\documentclass[12pt]{article}
\usepackage{amsthm}
\usepackage{amsfonts}
\usepackage{amsmath}
\usepackage{amssymb}
\usepackage[all]{xy}

\everymath{\displaystyle}
\CompileMatrices

\newtheorem{thm}{Theorem}[section]
\newtheorem{lemma}[thm]{Lemma}
\newtheorem{cor}[thm]{Corollary}
\newtheorem{prop}[thm]{Proposition}

\theoremstyle{definition}
\newtheorem{defn}[thm]{Definition}

\def\ve{\varepsilon}

\def\R{{\mathbb R}}

\def\C{{\mathbb C}}

\def\Z{{\mathbb Z}}

\def\so{_{\scriptscriptstyle O}}
\def\su{_{\scriptscriptstyle U}}
\def\st{_{\scriptscriptstyle T}}
\def\pr{^{\scriptstyle \R}}
\def\po{^{\scriptscriptstyle O}}
\def\pu{^{\scriptscriptstyle U}}
\def\pt{^{\scriptscriptstyle T}}
\def\sr{_{\scriptscriptstyle \R}}
\def\sc{_{\scriptscriptstyle \C}}
\def\crt{^{\scriptscriptstyle {\it CRT}}}
\def\scrt{_{\scriptscriptstyle {\it CRT}}}

\def\cp{{\mathcal{K}}}
\def\cpr{\cp\sr}

\def\top{{\mathcal{T}}}
\def\topr{\top\sr}
\def\topro{\top_{\scriptscriptstyle {\R, 0}}}

\def\id{\text {id}}

\def\im{\text {image} \,}
\def\coker{\text {coker}}
\def\tor{\text {Tor}}

\def\ac{A_{\scriptscriptstyle \C \,}}
\def\c{_{\scriptscriptstyle \C \,}}
\def\r{_{\scriptscriptstyle \R \,}}

\newcommand{\mapis}[1]{\buildrel {#1} \over \longrightarrow}
\newcommand{\sm}[4]{       \bigl( \begin{smallmatrix} 
					{#1} & {#2} \\ {#3} & {#4}
                   		    \end{smallmatrix} \bigr)         }
\newcommand{\mm}[4]{       \Bigl( \begin{smallmatrix} 
					{#1} & {#2} \\ {#3} & {#4}
                   		    \end{smallmatrix} \Bigr)         }

\newcommand{\smv}[2]{       \bigl( \begin{smallmatrix} 
					{#1} \\ {#2}
                   		    \end{smallmatrix} \bigr)         }

\newcommand{\smh}[2]{       \bigl( \begin{smallmatrix} 
					{#1} & {#2}
                   		    \end{smallmatrix} \bigr)         }
\newcommand{\smaaa}[3]{       \bigl( \begin{smallmatrix} 
					{#1} & {#2} & {#3}
                   		    \end{smallmatrix} \bigr)         }
\newcommand{\smabc}[3]{       \Bigl( \begin{smallmatrix} 
					{#1} \\ {#2} \\ {#3}
                   		    \end{smallmatrix} \Bigr)         }
\newcommand{\smaaabbb}[6]{       \bigl( \begin{smallmatrix} 
					{#1} & {#2} & {#3} \\
					{#4} & {#5} & {#6}
                   		    \end{smallmatrix} \bigr)         }
\newcommand{\mmaaabbb}[6]{       \Bigl( \begin{smallmatrix} 
					{#1} & {#2} & {#3} \\
					{#4} & {#5} & {#6}
                   		    \end{smallmatrix} \Bigr)         }
\newcommand{\smaabbcc}[6]{       \Bigl( \begin{smallmatrix} 
					{#1} & {#2} \\
					{#3} & {#4} \\
					{#5} & {#6}
                   		    \end{smallmatrix} \Bigr)         }

\begin{document}
\bibliographystyle{amsplain}


\title{Real C*-Algebras, United $K$-Theory, \\ and the K\"unneth Formula}
\date{\today}
\author{Jeffrey L. Boersema}

\maketitle

\vspace*{-0.525in}
\hspace*{4.71in}

\begin{center}
\quad \\[-10pt]
Department of Mathematics\\ 
Seattle University\\
Seattle, WA  98133\\
boersema@seattleu.edu
\end{center}

\setcounter{page}{1}

\thispagestyle{empty}

\vspace{0.2in}

\begin{abstract}  
\baselineskip15pt

We define united $K$-theory for real C*-algebras, generalizing Bousfield's 
topological united $K$-theory.  United $K$-theory incorporates three 
functors --- real $K$-theory,
complex $K$-theory, and self-conjugate $K$-theory --- and the natural transformations among
them.  The advantage of united $K$-theory over ordinary $K$-theory lies 
in its homological algebraic properties, which allow us to construct a 
K\"unneth-type, non-splitting, short exact
sequence whose middle term is the united $K$-theory of the tensor product of two
real C*-algebras $A$ and $B$ which holds as long as the complexification of $A$ is
in the bootstrap category $\mathcal{N}$.  Since united $K$-theory contains ordinary
$K$-theory, our sequence provides a way to compute the $K$-theory of the tensor
product of two real C*-algebras.  

As an application, we compute the united $K$-theory of the tensor product of two 
real Cuntz algebras.  Unlike in the complex case, it turns out that the isomorphism
class of the tensor product $\mathcal{O}\pr_{k+1} \otimes \mathcal{O}\pr_{l+1}$ is
not determined solely by the greatest common divisor of $k$ and $l$.  Hence, we have
examples of non-isomorphic, simple, purely infinite, real C*-algebras whose
complexifications are isomorphic.  \end{abstract}

\bigskip
\bigskip

\baselineskip18pt


\section*{Introduction}

In 1982, Claude Schochet (\cite{Schochet}) proved the existence of an exact 
K\"unneth-type sequence 
$$0 \rightarrow K_*(A) \otimes K_*(B) \rightarrow 
K_*(A \otimes B) \rightarrow 
\tor(K_*(A), K_*(B)) \rightarrow 0$$
which holds for complex C*-algebras $A$ and $B$ in which $A$ is in the bootstrap 
category $\mathcal{N}$.  Recall that $\mathcal{N}$ is the smallest subcategory of
complex, separable, nuclear C*-algebras which contains the separable type I
C*-algebras; which is closed under the operations of taking inductive limits, stable
isomorphisms, and crossed products by $\Z$ and $\R$; and which satisfies the two out
of three rule for short exact sequences (i.e. if $0 \rightarrow A \rightarrow B
\rightarrow C \rightarrow 0$ is exact and two of $A$, $B$, $C$ are in $\mathcal{N}$,
then the third is also in $\mathcal{N})$.

In this paper, we will develop an analogous K\"unneth sequence for real 
C*-algebras.  A real C*-algebra $A$ is a real Banach *-algebra which satisfies the C*-equation 
$\| a^*a \| = \|a\|^2$ as well as the axiom that $1+a^*a$ is invertible in the
unitization $A^+$ for all $a \in A$.  This is equivalent by \cite{Palmer} to saying
that $A$ is *-isometrically isomorphic to a norm-closed adjoint-closed algebra of
operators on a real Hilbert space.  

If $A$ is any real C*-algebra, it has a unique complexification $\C \otimes 
A$ (sometimes denoted $A\sc$ in this paper) which 
is a complex C*-algebra (see \cite{Goodearl} or \cite{Palmer}).  It is often useful
to complexify a real C*-algebra since complex C*-algebras are better understood.
However, the complexification by itself loses.  For 
example, the real C*-algebras $\R \oplus \R$ and $\C$ become isomorphic 
upon complexification.  Later in this paper we will describe examples of distinct real 
simple purely infinite C*-algebras whose complexifications are isomorphic.

Given two real C*-algebras $A$ and $B$, let $A \otimes\sr B$ 
be the real spacial tensor product.  The 
problem set before us is to express the $K$-theory of $A \otimes\sr B$ in terms of the 
$K$-theory of $A$ and of $B$ in the same way that Schochet's K\"unneth-type 
sequence expresses the $K$-theory of the tensor product of two complex 
C*-algebras in terms of the $K$-theory of the factors.

A key step in Schochet's proof (following Atiyah in \cite{Ati2}) is to build a free 
geometric resolution of $K_*(B)$ for any complex C*-algebra $B$.  This is a length
one resolution $$0 \rightarrow K_*(F_1) \xrightarrow{\mu_1} 
K_*(F_0) \xrightarrow{\mu_0} 
K_*(B) \rightarrow 0 $$
where $K_*(F_i)$ are free (in the category of $\Z$-graded 
$K_*(\C)$-modules) and the homomorphisms $\mu_i$ are induced by homomorphisms 
on the level of C*-algebras.  

Unfortunately, if $B$ is a real C*-algebra, $K_*(B)$ does not necessarily
have projective dimension one in the category of $K_*(\R)$-modules.  
Thus, it is impossible to build a
length-one free resolution of $K_*(B)$, so Schochet's proof does not adapt to give a
proof of a K\"unneth-type sequence for real C*-algebras.  In fact, counter-examples
showing that the analogous sequence is not exact can be found in \cite{Ati2} (for the
topological case) and \cite{Boer} (for a simpler example of C*-algebras).

United $K$-theory is the way around the problem of projective dimension.  
Following A.K. Bousfield's work in the topological case (see \cite{Bou}), we will
define $K\crt(A)$, the united $K$-theory of a real C*-algebra, as a triple of
$\Z$-graded abelian groups --- ordinary $K$-theory, complex $K$-theory, and
self-conjugate $K$-theory --- together with their respective module structures and
interrelating natural transformations.  This functor takes values in the category
{\it CRT} consisting of so-called {\it CRT}-modules.  

As in the topological case (and drawing
on a theorem of Bousfield), we find that $K\crt(B)$ has projective dimension
one in the category {\it CRT} and so the homological algebraic obstruction to
producing a K\"unneth Formula for real $K$-theory disappears.  Using this result, we
are able to produce free geometric resolutions and subsequently we prove the
following product formula.  

\vspace{.3cm}
{\bfseries Main Theorem} 
{\it Let $A$ and $B$ be real C*-algebras with $\C \otimes A \in \mathcal{N}$.  Then 
there is a short exact sequence} 
$$0 \rightarrow K\crt(A) \otimes\scrt K\crt(B) \rightarrow 
K\crt(A \otimes B) \rightarrow 
\tor\scrt(K\crt(A), K\crt(B)) \rightarrow 0\; .$$

This K\"unneth Sequence partially generalizes Schochet's result 
in \cite{Schochet} in the sense that by restricting attention to the complex part
of each of the {\it CRT}-modules in the exact sequence, we recover Schochet's K\"unneth
sequence for the complex C*-algebras $\C \otimes A$ and $\C \otimes B$.
 At the same time, it generalizes Bousfield's K\"unneth sequence for topological
united $K$-theory stated (but not proven) in \cite{Bou2}.  Indeed, Bousfield's
result can be recovered by taking $A = C(X, \R)$ and $B = C(Y, \R)$.  

It is also the case that $K\crt(B)$ has injective dimension one in the category 
{\it CRT}, a property that we intend to exploit in a subsequent paper to prove a
universal coefficient theorem for united $K$-theory, analogous to the universal
coefficient theorem for complex C*-algebras of \cite{RS}.

Since united $K$-theory incorporates ordinary $K$-theory, the main product theorem
provides a way of computing the ordinary $K$-theory of the tensor product of two
real C*-algebras.  As an example, we compute the united $K$-theory of the tensor
product of two real Cuntz algebras $\mathcal{O}\pr_{k+1}$ and
$\mathcal{O}\pr_{k+1}$.  In the case of complex Cuntz algebras, the $K$-theory of
the tensor product --- and therefore, the isomorphism class of the tensor product by
the classification theorems of Phillips and Kirchberg (\cite{Phillips} and
\cite{Kirchberg}) --- only depends on the greatest common divisor of $k$ and $l$. 
In the real case however, we have found that the greatest common divisor is not
enough to determine $K$-theory.  The $K$-theory of $\mathcal{O}\pr_{k+1} \otimes
\mathcal{O}\pr_{k+1}$ is a function of the greatest common divisor of $k$ and $l$ as
well as the values of $k$ and $l$ modulo 4.  As a particular example, it 
turns out that the real
C*-algebras $\mathcal{O}\pr_3 \otimes \mathcal{O}\pr_3$ and $\mathcal{O}\pr_3
\otimes \mathcal{O}\pr_5$ are not isomorphic, even though their complexifications
$\mathcal{O}_3 \otimes \mathcal{O}_3$ and $\mathcal{O}_3 \otimes \mathcal{O}_5$ are
isomorphic.

The organization of this paper is as follows.  In Section~\ref{ukt}, we define and 
develop united $K$-theory.  In Section~\ref{fgr}, we compute the united $K$-theory
for the real C*-algebras $\R$, $\C$, and $T$.  These are the prototypical free {\it
CRT}-modules and form the building blocks for the free geometric resolution which we
construct for an arbitrary real C*-algebra $B$.  Section~\ref{tp} is entirely
algebraic, wherein we define and develop the theory of tensor products in the
category {\it CRT} so that $K\crt(A) \otimes\scrt K\crt(B)$ and $\tor\scrt(K\crt(A),
K\crt(B))$ make sense in the statement of the main theorem.  Finally, in
Section~\ref{ks} we prove the main theorem and in Section~\ref{rca} we apply the
main theorem  to compute the united $K$-theory of the tensor product of two
arbitrary Cuntz algebras.

{\it Acknowledgments:}  The content of this paper represents work done for 
my dissertation at the University of Oregon.  I gratefully acknowledge my 
debt to my advisor, Chris Phillips, for his patience and helpfulness, both 
before and after graduation.


\section{United $K$-theory} \label{ukt}

We begin by recording the definition of united $K$-theory, even though the 
terms contained in it are as yet undefined.  It is the task of this chapter to make full 
sense of the definition.

\begin{defn} Let $A$ be a real C*-algebra.  The united $K$-theory of $A$ is the 
triple of $\Z$-graded abelian groups
$$K\crt(A) := \{KO_*(A), KU_*(A), KT_*(A)\} \, $$
together with the eight natural transformations 
$\{r, c, \ve, \zeta, \psi\su, \psi\st, \gamma, \tau \}$.
\end{defn}
The $K$-groups comprising united $K$-theory will be defined in 
Section~\ref{self-conjugate} and the natural transformations among the three graded groups 
will be described in Section~\ref{trans}.  In Section~\ref{crt-modules} we will record some important
properties of the category {\it CRT} which is the target category of United
$K$-theory.  Finally in Section~\ref{longexactsequences} we will prove the
existence of three long exact sequences involving the three graded groups comprising
united $K$-theory.

\subsection{Real, Complex, and Self-Conjugate $K$-theory} \label{self-conjugate}

Recall from Section~1.4 of \cite{Schroder} that the $K$-theory of a real
unital C*-algebra $A$ is defined to be the Grothendieck group of equivalence classes of
stable projections in $\cup_{k \in \Z_+} M_k(A)$.  The equivalence relation can be
taken in the sense of partial isometries, unitaries, or homotopies.  If $A$ is not
unital then $K$-theory is defined to be the kernel of the map $K(A^+) \rightarrow
K(\R)$.  For $n \geq 0$, we define $K_n(A) = K(S^n A)$ where $S^nA$ denotes the
$n$-fold suspension.  Finally, for $n < 0$, $K_n(A)$ is defined using the 8-fold
periodicity.  Note that in this paper, all $K$-groups are considered to be 
graded over $\Z$.

This definition of $K$-theory is identical with the familiar definition of
$K$-theory for complex C*-algebra (see \cite{Black}, \cite{Murphy}, or \cite{Wegge}) 
in the sense that if $A$ is a complex C*-algebra, the groups obtained would be the
same whether $A$ were considered a complex C*-algebra or just a real 
C*-algebra with a forgotten complex structure.

In this paper, we will use the phrase ``real $K$-theory'' to refer to the 
$K$-theory of a real C*-algebra as described above while the phrase 
``complex $K$-theory'' will refer to the $K$-theory of the complexification 
of a real C*-algebra, as in the following definitions.

\begin{defn} [Real $K$-theory] $KO_n(A) := K_n(A)$ for all $n \in \Z$.
\end{defn}

\begin{defn} [Complex $K$-theory] $KU_n(A) := K_n(\C \otimes A)$ for all $n \in \Z$.
\end{defn}

The real $K$-theory $KO_*(A)$ is a $\Z$-graded module over the ring $KO_*(\R)$, 
which according to \cite{Schroder} on p. 23, is given as below in gradings $0$
through $8$. $$
KO_*(\R) \hspace{.5cm} = \hspace{.5cm}
\Z \hspace{1cm} 
\Z_2 \hspace{1cm}
\Z_2 \hspace{1cm}
0 \hspace{1cm}
\Z \hspace{1cm}
0 \hspace{1cm}
0 \hspace{1cm}
0 \hspace{1cm}
\Z  $$
We recall that $KO_*(\R)$ is a graded skew-commutative ring generated by the unit 
$1\so \in KO_0(\R)$, $\eta\so \in KO_1(\R)$, $\xi \in KO_4(\R)$, and the Bott
element $\beta\so \in KO_8(\R)$.  These generators are subject to the relations $2
\eta\so = 0$, $(\eta\so)^3 = 0$, $\eta\so \cdot \xi = 0$, and $\xi^2 = 4  \beta\so$.
 Since the Bott element $\beta\so \in KO_8(\R)$ is invertible, multiplication by
$\beta\so$ induces a period 8 isomorphism of $KO_*(A)$.  

The complex $K$-theory $KU_*(A)$ is a $\Z$-graded module over the ring 
$KU_*(\R)$.
$$
KU_*(\R) \hspace{.5cm} = \hspace{.5cm}
\Z \hspace{1cm} 
0 \hspace{1cm}
\Z \hspace{1cm}
0 \hspace{1cm}
\Z \hspace{1cm}
0 \hspace{1cm}
\Z \hspace{1cm}
0 \hspace{1cm}
\Z  $$

The ring $KU_*(\R)$ is a free commutative polynomial ring over $\Z$ generated by 
the Bott element $\beta\su \in KU_2(\R)$ and its inverse $\beta\su^{-1} \in
KU_{-2}(\R)$.  In particular, the well-known Bott periodicity isomorphism for
complex $K$-theory is induced by multiplication by $\beta\su$.

Finally, we define self-conjugate $K$-theory.  Consider the real C*-algebra 
$$ T = \{f \in C([0,1], \C) \mid f(0) = \overline{f(1)}\} \, .$$  
If $A$ is any other real C*-algebra, we have 
$$T \otimes A \cong \{f \in C([0,1], \C \otimes A) \mid f(0) = \overline{f(1)} \}$$ 
where the conjugation of $\C \otimes A$ is defined by 
$\overline{\lambda \otimes a} = \bar{\lambda} \otimes a$.
  
\begin{defn} [Self-conjugate $K$-theory] $KT_n(A) := K_n(T \otimes A)$ for all $n \in \Z$.
\end{defn}

The motivation for this definition comes from the topological case.  Topological
 self-conjugate $K$-theory was first defined by D.W. Anderson in \cite{And} in
terms of self-conjugate vector bundles.  A self-conjugate vector bundle is a complex
vector bundle together with a given conjugate linear automorphism.  In the case of
real C*-algebras, such an object is replaced by a self-conjugate projection, which
is a projection $p$ in a matrix algebra over $\C \otimes A$ together with a given
equivalence to the conjugate projection $\overline{p}$.  If by ``equivalence'' we
mean homotopy equivalence, this is a continuous path of projections from $p$ to
$\overline{p}$, which is the same thing as a projection in a matrix algebra over $T
\otimes A$.

The rest of this section is devoted to developing the theory of self-conjugate 
$K$-theory for real C*-algebras.  A priori, $KT_*(A)$ has a periodicity of period 8 since $T \otimes A$
is a real C*-algebra, but we will eventually show that it has period 4
(see the discussion following Corollary~\ref{K(T)}).  Since $T$ is commutative there is a pairing $T
\otimes T \rightarrow T$ which makes $KT_*(\R)$ into a ring and makes $KT_*(A)$ into
a $KT_*(\R)$-module.  To compute the structure ring $KT_*(\R)$, we must develop the
relationship between self-conjugate and complex $K$-theory, particularly 
the long exact sequence relating the two.  

Let $\zeta$ be the C*-algebra homomorphism from $T$ to $\C$ defined by $\zeta(f) = f(0)$ 
for any $f \in T$.  If $A$ is any real C*-algebra, then the homomorphism
$$\zeta \otimes \id \colon T \otimes A \rightarrow \C \otimes A$$ 
induces a natural transformation $KT_n(A) \rightarrow KU_n(A)$, which we also call 
$\zeta$ or $\zeta_n$.  Similarly, the conjugation $\psi\su \colon \C \rightarrow
\C$ defines an involutive natural transformation $(\psi\su)_n \colon KU_n(A)
\rightarrow KU_n(A)$ for any C*-algebra $A$.

Finally, the inclusion $\gamma \colon S\C \rightarrow T$ defines a natural 
transformation $\gamma_{n} \colon KU_{n}(A) \rightarrow KT_{n-1}(A)$ as follows. 
For $n \geq 1$, the transformation $\gamma_n$ is induced by the homormophism
$$S^{n}(\C \otimes A) \rightarrow S^{n-1}(T \otimes A) $$ where by convention we say
that it is the outermost (or leftmost) suspension that gets used by $\gamma_n$ and
the innermost $n-1$ suspensions go along for the ride.  For $n < 0$, $\gamma_n$
is defined using the periodicity isomorphisms.

\begin{thm} \label{Seq1} There is a natural long exact sequence 
$$\cdots \longrightarrow  KU_{n+1}(A) \xrightarrow{\gamma} 
KT_n(A) \xrightarrow{\zeta} 
KU_n(A) \xrightarrow{1 - \psi\su}
KU_{n}(A) \longrightarrow \cdots.$$
\end{thm}

\begin{proof} The short exact sequence 
$$0 \longrightarrow 
S\C \otimes A \xrightarrow{\gamma} 
T \otimes A \xrightarrow{\zeta} \C \otimes A \longrightarrow 0$$
induces the long exact sequence
$$\cdots \longrightarrow  
K_{n}(S\C \otimes A)  \xrightarrow{\gamma} 
K_n(T \otimes A) \xrightarrow{\zeta}
K_n(\C \otimes A) \xrightarrow{\delta} 
K_{n-1}(S\C \otimes A) \longrightarrow \cdots.$$

Since $K_{n-1}(S\C \otimes A) \cong KU_n(A)$, this is the long exact sequence we want.  
It only remains to identify the index $\delta$.    We will prove that
$\delta = -(1 - \psi\su)$.  This is sufficient since the sign of the homomorphism
does not affect the exactness of the sequence.

We first consider the case $n = 1$.  Since the index map $\delta$ is most 
conveniently defined in terms of unitaries, we consider the unitary definition of
$K$-theory, denoted by $KU^{\mathcal{U}}_1(A)$,
and show that the following diagram commutes.  Here $\Theta$ is the usual
isomorphism from $KU^{\mathcal{U}}_1(A)$ to $KU_0(SA)$ as described in Section~7.2
of \cite{Wegge}. 
$$ \xymatrix{ KU^{\mathcal{U}}_{1}(A) \ar[dr]_{\delta}
\ar[r]^{\psi\su - 1}   & KU^{\mathcal{U}}_1(A)  \ar[d]^{\Theta}   \\      & 
KU_0(SA)  }$$

Let $u$ be a unitary element of $M_k(\C \otimes A)$.  Then 
$\delta([u]) = [V_t p_k V_t^*] - [p_k]$ where $V_t$ is a unitary lift of
$\sm{u}{0}{0}{u^*}$ in $M_{2k}(T \otimes A)$.  In other words, $V_t$ may be taken
to be any path in $\mathcal{U}_{2k}(\ac)$ from $\sm{u}{0}{0}{u^*}$ to
$\sm{\overline{u}}{0}{0}{\overline{u}^*}$.   The class $[V_t p_k V_t^*] - [p_k]$
is ostensibly an element of $KU_0(C(S^1, A))$; but since it is in the kernel of the
point-evaluation map $KU_0(C(S^1, A)) \rightarrow KU_0(A)$, it determines a unique
element of $KU_0(SA)$. 

On the other hand, $\Theta(\psi\su - 1)[u] = \Theta[\overline{u}u^*] = 
[W_t p_k W_t^*] - [p_k]$ where $W_t$ is any path in $\mathcal{U}_{2k}(A)$ from
$\sm{1}{0}{0}{1}$ to $\sm{\overline{u}u^*}{0}{0}{u \overline{u}^*}$.   Now,
$\overline{u}u^*$ is stably path-connected to $u^*\overline{u}$ through unitaries so
we can replace $W_t$ by any path $\widetilde{W_t}$ from  $\sm{1}{0}{0}{1}$ to
$\sm{u^*\overline{u}}{0}{0}{u \overline{u}^*}$.  Furthermore, the class
$[\widetilde{W_t} p_k \widetilde{W_t}^*]$ is equivalent to the class  $\bigl[
\sm{u}{0}{0}{u^*} \widetilde{W_t}  p_k \widetilde{W_t}^*  \sm{u}{0}{0}{u^*}^*
\bigr]$.   Therefore, in the definition of $\delta$, we take $V_t$ to be the path 
$\sm{u}{0}{0}{u^*} \widetilde{W_t}$  from $\sm{u}{0}{0}{u^*}$ to
$\sm{\overline{u}}{0}{0}{\overline{u}^*}$ and it follows that $\theta(\psi\su -
1)[u] = [V_t p_k V_t^*] - [p_k] = \delta([u])$. 

We have proven that $\delta = -(1 - \psi\su)$ for $n = 1$ and for all $A$.  
The general case follows using suspensions and periodicity.
\end{proof}

\begin{cor}  \label{K(T)}   
The self-conjugate $K$-theory of $\R$ in degrees zero through eight is     
$$KT_*(\R) \hspace{.5cm} = \hspace{.5cm}
\Z \hspace{1cm}
\Z_2 \hspace{1cm}
0 \hspace{1cm}
\Z \hspace{1cm}
\Z \hspace{1cm}
\Z_2 \hspace{1cm}
0 \hspace{1cm}
\Z \hspace{1cm}
\Z $$
\end{cor}

\begin{proof} Take $A = \R$ in the sequence of
Theorem~\ref{Seq1} and use 
\begin{align*}   
KU_n(\R) &= \begin{cases} \Z & n \equiv 0 \pmod 2 \\
	 	0 & n \equiv 2 \pmod 2~ \end{cases} \\
\intertext{and} 
(\psi\su)_n &=    \begin{cases} 1 & n \equiv 0 \pmod 4 \\
		-1 & n \equiv 2 \pmod 4~ \end{cases} \end{align*}
where $(\psi\su)_n \colon KU_n(\R) \rightarrow KU_n(\R)$ to compute 
$KT_*(\R)$.
\end{proof}

Let $1\st \in KT_0(\R)$ be the class of the unit.  Since $\zeta$ is a unital 
homomorphism, it sends $1\st$ to $1\su \in KU_0(\R)$.  The calculation above also
reveals that $\zeta \colon KT_4(\R) \rightarrow KU_4(\R)$ is an isomorphism.  We let
$\beta\st$ denote the generator of $KT_4(\R)$ satisfying $\zeta(\beta\st) =
\beta\su^2$.  We let $\eta\st = \gamma(\beta\su)$ denote the non-zero element of $KT_1(\R)$ 
and we let $\omega$ denote $\gamma(\beta\su^2)$, which is a generator of 
$KT_3(\R)$.  

Since $\zeta$ is a ring homomorphism and an isomorphism in degree $4k$ for 
all integers $k$, it follows that the element $\beta\st$ is invertible, 
where the inverse is given by the element $\zeta^{-1}(\beta\su^{-2})$.  
This fact and the relations $(\eta\st)^2 = 0$, $\omega^2 = 0$, and 
$\eta\st \cdot \omega = \omega \cdot \eta\st = 0$ (which are forced by 
dimension considerations) completely determine the ring structure of $KT_*(\R)$.

Since $\C$ and $T$ are commutative, there are external products in real, complex, and 
self-conjugate $K$-theory for any real C*-algebras $A$ and $B$.
\begin{align*}
{\alpha\so} \colon KO_m(A) \otimes KO_n(B) &\rightarrow KO_{m+n}(A \otimes B) \\
{\alpha\su} \colon KU_m(A) \otimes KU_n(B) &\rightarrow KU_{m+n}(A \otimes B) \\
{\alpha\st} \colon KT_m(A) \otimes KT_n(B) &\rightarrow KT_{m+n}(A \otimes B) \\
\end{align*}

We will use the notation $x \cdot\so y$, $x \cdot\su y$, and $x \cdot\st y$ to denote 
the real, complex, and self-conjugate products respectively.


\subsection{Natural Transformations} \label{trans}

In this subsection we define the rest of the eight natural transformations which 
interrelate the real, complex, and self-conjugate $K$-theories.  These
transformations are \begin{align*}   \label{operations}
c_n &\colon KO_n(A) \longrightarrow KU_n(A)  &
r_n &\colon KU_n(A) \longrightarrow KO_n(A) \\
\varepsilon_n &\colon KO_n(A) \longrightarrow KT_n(A) &
\zeta_n &\colon KT_n(A) \longrightarrow KU_n(A) \\
(\psi\su)_n &\colon KU_n(A) \longrightarrow KU_n(A)  &
(\psi\st)_n &\colon KT_n(A) \longrightarrow KT_n(A) \\
\gamma_n &\colon KU_n(A) \longrightarrow KT_{n-1}(A)  &
\tau_n &\colon KT_n(A) \longrightarrow KO_{n+1}(A)  \; .
\end{align*}

The transformations $\zeta$, $\psi\su$, and $\gamma$ have already been defined, 
induced by C*-algebra homomorphisms on the appropriate C*-algebras. In the same
way, $c$, $r$, $\varepsilon$, and $\psi\st$ are induced by C*-algebra homomorphisms 
 \begin{align*} c &\colon \R \longrightarrow \C  &
	r &\colon \C \longrightarrow M_2(\R) \\
\varepsilon &\colon \R \longrightarrow T  & 
	\psi\st &\colon T \longrightarrow T  \; 
\end{align*}
defined as follows.  We  take $c$ and $\varepsilon$ to be the unital 
inclusion maps.  
The realification $r$ is defined by $r(x + iy) = \sm{x}{-y}{y}{x}$ for all $x + iy
\in \C$.  Finally, let $f \in T$ be a path from $f(0)$ to $f(1) =
\overline{f(0)}$.  We define $\psi\st(f)\in T$ to be the reverse path from $f(1)$ to
$f(0)$ --- that is, $\psi\st(f)(s) = f(1-s)$ for all $s \in [0,1]$.

Defining $\tau$ will take a little more work.  First, we will define two homomorphisms 
$\sigma_1$ and $\sigma_2$ from $T$ to $C(S^1, M_4(\R))$.  
  
To define $\sigma_1$, let $U_s$ be any path in $\mathcal{U}_4(\R) = O(4)$ 
$$ \text{from} \qquad
\begin{pmatrix} 1 & 0 & 0 & 0 \\ 0 & 1 & 0 & 0 \\ 
	0 & 0 & 1 & 0  \\ 0 & 0 & 0 & 1 \end{pmatrix}
\qquad \text{to} \qquad
\begin{pmatrix} 0 & 1 & 0 & 0 \\ 1 & 0 & 0 & 0 \\
	 0 & 0 & 0 & 1  \\ 0 & 0 & 1 & 0  \end{pmatrix} \; .$$

Let $f$ be an element of $T$ and write $f(0) = x + iy \in \C$.  Then $f(1) = x - iy$.  
Now $U_s \cdot \sm{r(f(1))}{0}{0}{0} \cdot U_s^*$ is a path in $M_4(\R)$  
$$ \text{from} \qquad
\begin{pmatrix} x & y & 0 & 0\\ -y & x & 0 & 0 \\ 
	0 & 0 & 0 & 0 \\ 0 & 0 & 0 & 0  \end{pmatrix}
\qquad \text{to} \qquad
\begin{pmatrix} x & -y & 0 & 0 \\ y & x & 0 & 0 
	\\ 0 & 0 & 0 & 0 \\ 0 & 0 & 0 & 0   \end{pmatrix} \;  .$$

On the other hand, $r(f)$ is a path from $\sm{x}{-y}{y}{x}$ to $\sm{x}{y}{-y}{x}$.  We 
define $\sigma_1(f)$ to be the concatenation 
$$\sm{r(f)}{0}{0}{0} \; \# \; U_s \cdot \sm{r(f(1))}{0}{0}{0} \cdot U_s^* \, $$  
which is a loop in $M_4(\R)$ based at 
$$\sm{r(f(0))}{0}{0}{0} = \begin{pmatrix} x & -y & 0 & 0 \\ y & x & 0 & 0 
	\\ 0 & 0 & 0 & 0 \\ 0 & 0 & 0 & 0   \end{pmatrix} \; $$

We define $\sigma_2(f)$ to be simply the constant path with value 
$r(f(0)) = \sm{x}{-y}{y}{x} \in M_2(\R) \subset M_4(\R)$.

Now each of the homomorphisms $\sigma_1$ and $\sigma_2$ induces, for any real C*-algebra 
$A$, a natural homomorphism 
$$(\sigma_i)_* \colon KT_0(A) \longrightarrow 
		KO_0(C(S^1, M_4(A))) \cong KO_0(C(S^1, A)) \; .$$  
We define $\tau = (\sigma_1)_* - (\sigma_2)_*$.  This homomorphism takes values in 
$KO_0(C(S^1, A))$.  But elements in the image of $\tau$ are in the kernel of the
base-point evaluation map  $$KO_0(C(S^1, A)) \longrightarrow KO_0(A) $$ 
and hence determine elements of $KO_0(SA) = KO_1(A)$.  Therefore $\tau$ is a 
homomorphism $KT_0(A) \rightarrow KO_1(A)$.  By use of suspensions and periodicity,
we then have a natural transformation  $$\tau_n \colon KT_n(A) \rightarrow
KO_{n+1}(A)  $$  for all $n \in \Z$.

\begin{prop} \label{crtobject} 
For any real C*-algebra $A$, the natural transformations among $KO_*(A)$, $KU_*(A)$, 
and $KT_*(A)$ 
satisfy the following relations
\begin{align*}  
rc &= 2    & \psi\su \beta\su &= -\beta\su \psi\su & \xi &= r \beta\su^2 c \\
cr &= 1 + \psi\su & \psi\st \beta\st &= \beta\st \psi\st 
	&  \omega &= \beta\st \gamma \zeta \\ 
r &= \tau \gamma & \varepsilon \beta\so &= \beta\st^2 \varepsilon 
	& \beta\st \varepsilon \tau 
			&= \varepsilon \tau \beta\st + \eta\st \beta\st    \\
c &= \zeta \varepsilon & \zeta \beta\st &= \beta\su^2 \zeta
	 &  \varepsilon r \zeta &= 1 + \psi\st   \\
(\psi\su)^2 &= 1 & \gamma \beta\su^2 &= \beta\st \gamma       
	&  \gamma c \tau &= 1 - \psi\st \\
(\psi\st)^2 &= 1 & \tau \beta\st^2 &= \beta\so \tau & \tau &= -\tau \psi\st  \\
\psi\st \varepsilon &= \varepsilon & 
	\gamma &= \gamma \psi\su 
	\qquad & \tau \beta\st \varepsilon &= 0 \\
\zeta \gamma &= 0 & \eta\so &= \tau \varepsilon 
	& \varepsilon \xi &= 2 \beta\st \varepsilon \\
\zeta &= \psi\su \zeta & \eta\st &= \gamma \beta\su \zeta 
	& \xi \tau &= 2 \tau \beta\st \; .
\end{align*}
\end{prop}

In this section, we will prove the relations in the first two columns.  The
relations  in the third column will be proven in subsequent sections concluding in
Section~\ref{geosurjection}.  As we go, we will be careful to use only those
relations which have been proven.  

\begin{lemma} \label{mult2} Let $x_1 \in KO_i(A)$, $x_2 \in KO_j(B)$, 
$y_1 \in KU_k(C)$, $y_2 \in KU_l(D)$,
$z_1 \in KT_m(E)$, and $z_2 \in KT_n(F)$.  Then \begin{enumerate}
\item[{\rm(1)}] \label{p1}$\psi\su(y_1 \cdot\su y_2) 
		= \psi\su(y_1) \cdot\su \psi\su(y_2)$
\item[{\rm(2)}] \label{p2}$\zeta(z_1 \cdot\st z_2) = \zeta(z_1) \cdot\su \zeta(z_2)$
\item[{\rm(3)}] \label{p3}$\gamma(y_1 \cdot\su \zeta(z_1)) = \gamma(y_1) \cdot\st z_1$
\item[{\rm(4)}] \label{p4}$\gamma(\zeta(z_1) \cdot\su y_1) 
			= (-1)^m z_1 \cdot\st \gamma(y_1)$ 

\item[{\rm(5)}] \label{pt-3} $c(x_1 \cdot\so x_2) = c(x_1) \cdot\su c(x_2)$
\item[{\rm(6)}] \label{pt-2} $r(c(x_1) \cdot\su y) = x_1 \cdot\so r(y)$
\item[{\rm(7)}] \label{pt-1} $r(y \cdot\su c(x_1)) = r(y) \cdot\su x_1$
\item[{\rm(8)}] \label{pt1} $\varepsilon(x_1 \cdot\so x_2) 
	= \varepsilon(x_1) \cdot\st \varepsilon(x_2)$
\item[{\rm(9)}] \label{pt2} $\psi\st(z_1 \cdot\st z_2) 
	= \psi\st(z_1) \cdot\st \psi\st(z_2)$
\item[{\rm(10)}] \label{pt3} $\tau(z_1 \cdot\st \varepsilon(x_1)) 
	= \tau(z_1) \cdot\so x_1$
\item[{\rm(11)}] \label{pt4} $\tau(\varepsilon(x_1) \cdot\st z_1) 
	= (-1)^i x_1 \cdot\so \tau(z_1)$
\end{enumerate} \end{lemma}

\begin{proof} 

Parts~(1), (2), (5), (8), and (9) follow from the facts that $\psi\su$, 
$\zeta$, $c$, $\varepsilon$, and $\psi\st$ are ring homomorphisms on the 
level of real C*-algebras.

To show part~(3), we first show that the diagram
$$ \xymatrix{
S\C \otimes T \ar[d]^{\gamma \otimes 1} \ar[r]^{1 \otimes \zeta} 
& S\C \otimes \C \ar[r]^-{\cdot}
& S\C \ar[d]^{\gamma}  \\
T \otimes T \ar[rr]^{\cdot} & & T }$$
commutes up to homotopy, where $\cdot$ represents multiplication in $\C$ or $T$.  
Indeed, for $f \in S\C$, $g \in T$, and $s \in [0,1]$; let $F_s(f, g) \in T$ be
defined by $F_s(f, g)(t) = f(t) \cdot g(st)$.  Since $F_s(f, g)(0) = F_s(f, g)(1)=
0$, we see that $F_s(f, g)$ is in $T$ for each $s$.  Furthermore, 
$F_0(f,g) = \gamma(\alpha\su(f \otimes \zeta(g)))$ and 
$F_1(f,g) = \alpha\st(\gamma(f) \otimes g)$.

This implies that, up to homotopy, the diagram
$$ \xymatrix{S^k(\C \otimes C) \otimes S^m(T \otimes A) 
		\ar[d]^{\gamma \otimes 1} \ar[r]^{1 \otimes \zeta} 
& S^k(\C \otimes C) \otimes S^m(\C \otimes A)  \ar[r]^-{\cdot}
& S^{k+m}(\C \otimes C \otimes A)   \ar[d]^{\gamma}     \\
S^{k-1}(T \otimes C) \otimes S^{m}(T \otimes A)    \ar[rr]^{\cdot} 
& & S^{k+m-1}(T \otimes C \otimes A)		}$$
commutes
(It is at this point where we rely on our convention regarding 
suspensions in the definition of $\gamma_k$.) and hence in $K$-theory the formula
$\gamma(z_1 \cdot\su \zeta(y_1)) = \gamma(z_1) \cdot\st y_1$ holds.   
Part~(4) follows from part~(3) using skew-commutativity.

Part~(6) follows from the following commutative diagram which 
expresses that $r$ is 
an $\R$-module homomorphism:
$$ \xymatrix{
\R \otimes \C \ar[rr]^-{\id \otimes r} \ar[d]^{c \otimes \id}
& & \R \otimes M_2(\R) \ar[d]^{\cdot}  \\
\C \otimes \C \ar[r]^-{\cdot} & \C \ar[r]^-{r} & M_2(\R) \\ }$$

Now we move to part~(10).  For each $i$, the homomorphism 
$\sigma_i \colon T \rightarrow C(S^1, M_4(\R))$ is an $\R$-module
homomorphism.  Therefore, the following diagram commutes: $$ \xymatrix{
T \otimes \R \ar[rr]^-{\sigma_i \otimes \id} \ar[d]^{\id \otimes \varepsilon}
& & C(S^1, M_4(\R)) \ar[d]^{\cdot} \otimes \R \\
T \otimes T \ar[r]^-{\cdot} & T \ar[r]^-{\sigma_i} & C(S^1, M_4(\R)) \\ }$$
This establishes the formula 
$(\sigma_i)_*(z_1 \cdot\st \varepsilon(x_1))  = (\sigma_i)_*(z_1) \cdot\so x_1$ for 
each $i$.  Therefore, since $\tau = (\sigma_1)_* - (\sigma_2)_*$, it follows that
$\tau(z_1 \cdot\st \varepsilon(x_1))  = \tau(z_1) \cdot\so x_1$.

Finally, using skew-commutativity, parts (4), (7), and (11) follow from 
parts (3), (6), and (10).
\end{proof} 

\begin{proof}[Beginning of Proof of \ref{crtobject}]
The composition $r \circ c$ from $\R$ to $M_2(\R)$ is given by 
$x \mapsto \sm{x}{0}{0}{x}$ 
which induces multiplication by 2 on real $K$-theory.  The composition $c \circ r$
from $\C$ to $M_2(\C)$ is given by $x + iy \mapsto \sm{x}{-y}{y}{x}$.  This map is
unitarily equivalent (via the unitary $\tfrac{1}{\sqrt{2}}\sm{1}{i}{i}{1}$) to the
map $x + iy \mapsto \sm{x + iy}{0}{0}{x - iy}$\,, which induces $1 + \psi\su$ on
complex $K$-theory.  This proves the first two relations in the first column.

Next we prove the relation $r = \tau \gamma$.  Let $p$ be a projection in $S\C \otimes A$ 
and consider it (via the inclusion $\gamma$) as an element of $T \otimes A$.  Then 
$$\tau \gamma[p] = \sigma_1[p] - \sigma_2[p] 
	= [r(p) \; \# \; U_s \cdot r(p(1)) \cdot U_s^*] - [r(p(0))].$$  
But $r(p(0))$ and $U_s \cdot r(p(1)) \cdot U_s^*$ are both identically zero since 
$p(0) = p(1) = 0$.  
Therefore, $\tau \gamma[p] = [r(p)]$.  

The next five relations in the first column are satisfied directly on the level of 
C*-algebra homomorphisms.  
The relations $\zeta = \psi\su \zeta$
and $\gamma = \gamma \psi\su$ follow from the exact sequence of Theorem~\ref{Seq1}.

Now we prove the first six relations in the second column.  In preparation, we recall the 
formulas $\psi\su(\beta\su) = -\beta\su$ and $\zeta(\beta\st) = \beta\su^2$. 
Furthermore, we have the formulas $\varepsilon(\beta\so) = \beta\st^2$ and
$\psi\st(\beta\st) = \beta\st$.  Indeed, the first identity is established from
$\zeta \circ \varepsilon = c$ and $c(\beta\so) = \beta\su^4$ using the fact that
$\zeta_8$ is an isomorphism.  The second identity is established from the identities
$\psi\su(\beta\su^2) = \beta\su^2$ and $\zeta \psi\st = \psi\su \zeta$ 
(the latter 
holds on the level of C*-algebra homomorphisms) using the fact
that $\zeta_4$ is an isomorphism taking $\beta\st$ to $\beta\su^2$.

Let $x \in KO_n(A)$, $y \in KU_n(A)$, and $z \in KT_n(A)$.  Then using 
Lemma~\ref{mult2} we have
\begin{align*} \psi\su \beta\su(y) 
&= \psi\su(\beta\su \cdot\su y) = \psi\su(\beta\su) \cdot\su \psi\su (y) =
-\beta\su \cdot\su \psi\su(y) = -\beta\su \psi\su(y) ~,\\ \psi\st \beta\st (z) &=
\psi\st(\beta\st \cdot \st z) = \psi\st(\beta\st) \cdot\st \psi\st(z) = \beta\st
\cdot\st \psi\st(z) = \beta\st \psi \st(z) ~,\\ \varepsilon \beta\so(x) &=
\varepsilon(\beta\so \cdot \so x) = \varepsilon(\beta\so) \cdot\st \varepsilon(x) =
\beta\st^2 \cdot\st \varepsilon(x) = \beta\st^2 \varepsilon(x) ~,\\ \zeta
\beta\st(z) &= \zeta(\beta\st \cdot\st z) = \zeta(\beta\st) \cdot\su \zeta(z) =
\beta\su^2 \cdot\su \zeta(z) = \beta\su^2 \zeta(z)  ~,   \\ \gamma \beta\su^2(y) &=
\gamma(\beta\su^2 \cdot\su y) = \gamma(\zeta(\beta\st) \cdot\su y) = \beta\st
\cdot\su \gamma(y) = \beta\st \gamma(y)  ~, \\ \text{and} \quad \qquad \tau
\beta\st^2(z) &= \tau(\beta\st^2 \cdot\st z) = \tau(\varepsilon(\beta\so) \cdot\st
z) = \beta\so \cdot\so \tau(z) = \beta\so \tau(z)~.  \end{align*}

We end by proving the relations $\eta\so = \tau \varepsilon$ and $\eta\st 
= \gamma \beta\su \zeta$.  First we prove that $\tau \varepsilon(1\so) =
\eta\so~\in~KO_1(\R)$.  Now  $$\tau(1\st) = \sigma_1(1\st) - \sigma_2(1\st) 
= [1_2 \; \# \; U_s \cdot 1_2 \cdot U_s^*] -   [1_2] 
= [U_s \cdot \sm{1}{0}{0}{1} \cdot U_s^*] - [\sm{1}{0}{0}{1}].$$
But since $U_s$ is a path in $M_4(\R)$ from $\sm{1_2}{0}{0}{1_2}$ to 
$\sm{u}{0}{0}{u^*}$ where $u = \sm{0}{1}{1}{0}$, we have 
$$\tau \varepsilon(1\so) = \tau(1\st) = \Theta[\sm{0}{1}{1}{0}] 
= \Theta[\sm{-1}{0}{0}{1}] = \Theta[-1] = \eta\so \; .$$

This proves that $\tau \varepsilon (1\so) = \tau (1\st) = \eta\so$.  In general, 
we let $x \in KO_n(A)$ and we use Lemma~\ref{mult2}, parts (8) and (10):
\begin{align*}  
\tau \varepsilon(x) &= \tau \varepsilon(1\so \cdot\so x) 
= \tau(\varepsilon(1\so) \cdot\st \varepsilon(x)) \\
&= \tau \varepsilon(1\so) \cdot\so x = \eta\so \cdot\so x.
\end{align*}
Therefore $\eta\so = \tau \varepsilon$.

Finally, we compute 
$$\gamma \beta\su \zeta(1\st) = \gamma(\beta\su \cdot\su 1\su) = \gamma(\beta\su) 
= \eta\st\; \in KO_1(\R) \;.$$ 
In general, then, we have for any $z \in KT_n(A)$,
\begin{align*}
\gamma \beta\su \zeta(z) &= \gamma \beta\su \zeta(1\st \cdot\st z) 
= \gamma \beta\su(\zeta(1\st) \cdot\su \zeta(z))     \\
&= \gamma(\beta\su \zeta(1\st) \cdot\su \zeta(z)) 
= (\gamma \beta\su \zeta(1\st) \cdot\st z) 
= \eta\st \cdot\st z
\end{align*}
using Lemma~\ref{mult2}, parts (2) and (3).
\end{proof}


\subsection{{\it CRT}-modules} \label{crt-modules}

This section is a summary of sections 1, 2, and 3 of \cite{Bou} describing the category {\it CRT} 
which is the target category of united $K$-theory.

\begin{defn} (Section~2.1 of \cite{Bou}) A {\it CRT}-module is a triple
$M = \{ M\po, M\pu, M\pt \}$ of $\Z$-graded abelian groups where $M\po$ is a 
$KO_*(\R)$-module, $M\pu$ is a $KU_*(\R)$-module, and $M\pt$ is a
$KT_*(\R)$-module.  In addition to the operations implied by the module
structures, each object $M$ comes with {\it CRT}-operations  \begin{align*}    c_n
&\colon M\po_n \longrightarrow M\pu_n  & r_n &\colon M\pu_n \longrightarrow M\po_n \\
\varepsilon_n &\colon M\po_n \longrightarrow M\pt_n &
\zeta_n &\colon M\pt_n \longrightarrow M\pu_n \\
(\psi\su)_n &\colon M\pu_n \longrightarrow M\pu_n  &
(\psi\st)_n &\colon M\pt_n \longrightarrow M\pt_n \\
\gamma_n &\colon M\pu_n \longrightarrow M\pt_{n-1}  &
\tau_n &\colon M\pt_n \longrightarrow M\po_{n+1} 
\end{align*}
which are $KO_*(\R)$-module homomorphisms and satisfy all the relations of 
Proposition~\ref{crtobject}.
\end{defn}

\begin{defn}
A {\it CRT}-morphism $\phi \colon M \rightarrow N$ is a triple $\phi 
= \{\phi\po, \phi\pu, \phi\pt \}$ where $\phi\po \colon M\po \rightarrow N\po$ is a
graded $KO_*(\R)$-module homomorphism, $\phi\pu \colon M\pu \rightarrow N\pu$ is a
graded $KU_*(\R)$-module homomorphism, and $\phi\pt \colon M\pt \rightarrow N\pt$ is
a graded $KT_*(\R)$-module homomorphism.  Furthermore, $\phi$ must commute 
with the eight
transformations $r$, $c$, $\varepsilon$, $\zeta$, $\psi\su$, $\psi\st$, $\gamma$,
and $\tau$. \end{defn}

\begin{defn}
The abelian category {\it CRT} is the category whose objects are {\it CRT}-modules 
and whose morphisms are {\it CRT}-morphisms.
\end{defn}

\begin{thm}
United $K$-theory is a natural, exact, homotopy invariant, continuous, stable, 
covariant functor from the category of real C*-algebras to the category {\it CRT}.
\end{thm}

\begin{proof}
These properties all hold for ordinary real $K$-theory.  They hold then for 
for complex and self-conjugate $K$-theory because the operation of tensoring a
C*-algebra by $\C$ or by $T$ commutes with stabilization and is natural, exact,
homotopy invariant, continuous, and covariant.  Finally, since the homomorphisms
involved with these properties are natural, they commute with the internal {\it
CRT}-module homomorphisms and thus are {\it CRT}-module homomorphisms. \end{proof}

We now record some general results about the category {\it CRT} from \cite{Bou} 
which will be necessary in what follows.

\begin{defn} (Section~2.3 of \cite{Bou}) An object $M = \{M\po, M\pu, M\pt\}$ in 
the category {\it CRT} is said to be acyclic if the sequences  
\begin{equation} 
\cdots \longrightarrow  
M\pu_{n+1} \mapis{\gamma} 
M\pt_n \mapis{\zeta} 
M\pu_n \xrightarrow{1 - \psi\su}
M\pu_{n} \longrightarrow \cdots 
\end{equation} 
\begin{equation} 
\cdots \longrightarrow  
M\po_n \mapis{\eta\so} 
M\po_{n+1} \mapis{c} 
M\pu_{n+1} \mapis{r \beta\su^{-1}}
M\po_{n-1} \longrightarrow \cdots 
\end{equation}
\begin{equation} 
\cdots \longrightarrow  
M\po_{n} \mapis{\eta\so^2} 
M\po_{n+2} \mapis{\varepsilon} 
M\pt_{n+2} \xrightarrow{\tau \beta\st^{-1}}
M\po_{n-1} \longrightarrow \cdots 
\end{equation}
are exact.
\end{defn}

We will show in Section~\ref{longexactsequences} that the united $K$-theory of any real 
C*-algebra is acyclic.  The acyclicity condition imposes a high degree of rigidity
on the {\it CRT}-module, as demonstrated by the following propositions found in
Section~2.3 of \cite{Bou}.

\begin{prop} \label{CRTrigidity} Suppose 
$\phi \colon M \rightarrow N$ is a homomorphism of
acyclic objects of {\it CRT}.  If one of $\phi\po$, $\phi\pu$, and $\phi\pt$ is an
isomorphism, then the other two are also isomorphisms. \end{prop}

\begin{prop} Suppose $M = \{M\po, M\pu, M\pt \}$ is an acyclic object of {\it CRT}.  
If one of the groups $M\po$, $M\pu$, and $M\pt$ is trivial, then the other two are
also trivial. \end{prop}

The following two theorems of Bousfield will be used in 
Section~\ref{fgr} to form 
geometric, free, length-one resolutions of $K\crt(B)$ for any real C*-algebra $B$.

\begin{thm} (Theorem~3.4 of \cite{Bou}) \label{acyclic=pd1} An object $M$ in 
\it{CRT} has projective dimension at most 1 if and only if it is acyclic.
\end{thm}

\begin{thm} (Theorem~3.2 of \cite{Bou}) \label{free=pro} For an object 
$M \in {\it CRT}$, the following are equivalent
\begin{enumerate}
\item[{\rm(1)}] $M$ is projective
\item[{\rm(2)}] $M$ is free
\item[{\rm(3)}] $M$ is acyclic and $M\pu$ is a free abelian group
\end{enumerate}
\end{thm}

Furthermore, Bousfield explicitly describes the free {\it CRT}-modules which are 
generated by one element.  Arbitrary free {\it CRT}-modules are direct sums of
monogenic free {\it CRT}-modules.  We will come across examples of these monogenic
free {\it CRT}-modules in Section~\ref{computation} as the united $K$-theory of
$\R$, $\C$, and $T$.


\subsection{Long Exact Sequences} \label{longexactsequences}

\begin{thm} \label{acyclic} For any real C*-algebra $A$, the following sequences 
are exact:
\end{thm}
\begin{equation} \label{utseq} 
\cdots \longrightarrow  KU_{n+1}(A) \mapis{\gamma} 
KT_n(A) \mapis{\zeta} 
KU_n(A) \xrightarrow{1 - \psi\su}
KU_{n}(A) \longrightarrow \cdots \end{equation}
\begin{equation} \label{ouseq} 
\cdots \longrightarrow  KO_n(A) \mapis{\eta\so} 
KO_{n+1}(A) \mapis{c} 
KU_{n+1}(A) \mapis{r \beta\su^{-1}}
KO_{n-1}(A) \longrightarrow \cdots \end{equation}
\begin{equation} \label{otseq} 
\cdots \longrightarrow  KO_{n}(A) \mapis{\eta\so^2} 
KO_{n+2}(A) \mapis{\varepsilon} 
KT_{n+2}(A) \xrightarrow{\tau \beta\st^{-1}}
KO_{n-1}(A) \longrightarrow \cdots \end{equation}

In the language of \cite{Bou}, this theorem says that $K\crt(A)$ is acyclic.  It 
is an immediate corollary of Theorems~\ref{acyclic=pd1} and \ref{acyclic} that
$K\crt(A)$ has projective dimension at most $1$ in the category {\it CRT}.  
However, this corollary will not be justified until we finish the proof of
Proposition~\ref{crtobject}.  In this section, we will not use any of the unproven
relations of that proposition, nor will we use any consequences of them.

In Section~\ref{self-conjugate}, we proved the exactness of Sequence~\ref{utseq}.  
Sequence~\ref{ouseq} can be found in \cite{Schroder} or \cite{Karoubi}, but we will
give a new proof which avoids Clifford algebras.  Using similar techniques we will
also prove the exactness of Sequence~\ref{otseq} which is new in the 
non-commutative context; in the topological case, it is due to Anderson \cite{And}.  

We will use the following notation in this section.

\begin{defn} Let $A$ be any real C*-algebra, let $X$ be any locally compact
 topological space, and let $\tau$ be any involution of $X$.  Then $C(X,A)$ denotes
the algebra of continuous functions from $X$ to $A$ and $C_0(X,A)$ is the subalgebra
of those functions which vanish at infinity.  We further define \begin{align*}
C(X;\tau) &= \{f \in C(X,\C) \mid f(\tau(x)) = \overline{f(x)} \}   \\ C_0(X;\tau)
&= \{f \in C_0(X,\C) \mid f(\tau(x)) = \overline{f(x)} \}   \\ IA &= C([0,1], A) \\
SA &= C_0(\R, A) \\
S^{-1}A &= \{f \in C_0(\R, \C \otimes A) \mid f(-x) = \overline{f(x)} \} \\
S^{m,n} A &= S^m (S^{-1})^n A \; .
\end{align*}
\end{defn}

The operation $A \mapsto SA$ is the usual suspension operation.  The operation 
$A \mapsto S^{-1}A$ is an inverse suspension operation in the
sense of the following proposition.

\begin{prop} \label{KKequiv}
The real C*-algebras $\R$ and $S^{1,1}$ are $KK$-equivalent.  Hence, for any real
C*-algebra $A$, there is an isomorphism $K(S^{m,n} A) \cong K_{m-n}(A)$. 
 \end{prop}

In Section~1.5 of \cite{Schroder}, Schr\"oder adapts the Toeplitz algebra proof of 
Bott periodicity to the case of real C*-algebras to show that $K_n(A) \cong
K_{n}(S^{1,1}A)$.  Our proposition is somewhat stronger, but the proof takes the
same approach using the real Toeplitz algebra.  The proof
will follow the next lemma.  Recall that a short exact sequence is semisplit if it
has a linear completely positive contractive section.

\begin{lemma} \label{cpositive}
A short exact sequence of real C*-algebras is semisplit if the short exact sequence 
obtained by complexifying is semisplit.
\end{lemma}

\begin{proof}
Let $s$ be a section for the complexified short exact sequence
$$0 \rightarrow \C \otimes A 
	\rightarrow \C \otimes B
	\xrightarrow{1 \otimes \pi} \C \otimes C
	\rightarrow 0  \; $$
so that $(1 \otimes \pi) \circ s = 1$.  If we simply restrict $s$ to $C = 1 \otimes
C \subset \C \otimes C$, we cannot be guaranteed  that the image will lie in $B$. 
However, we claim that the linear projection $f \colon \C \otimes B \rightarrow B$
defined by $f(b_1 + i b_2) = b_1$ is contractive and completely positive.  In that
case the composition $f \circ s\rvert_{1 \otimes C} \colon C \rightarrow B$ is the
desired section for the short exact sequence 
$$0 \rightarrow A \rightarrow B \xrightarrow {\pi} C
\rightarrow 0 \; .$$

To prove the claim, assume that $B$ is an algebra of operators on a real 
Hilbert space $\mathcal{H}$.  Then $\C \otimes B$
is an algebra of operators on the complex Hilbert space $\C \otimes \mathcal{H}$.
    
We show that $f$ is positive using the condition (Theorem~VIII.3.8 in \cite{Conway}) 
that 
an operator $T$ on a Hilbert space is positive if and only if $\langle T x, x
\rangle \geq 0$ for all vectors $x$ in the Hilbert space.
Let $b_1 + i b_2$ be an arbitrary positive element in $\C \otimes B$.  Then
$b_1$ is positive since for all $x \in \mathcal{H}$, 
$$0 \leq \langle (b_1 + i b_2)x, x \rangle  	
	= \langle b_1 x, x \rangle + i \langle b_2 x, x \rangle  \; .$$ 
Since $\langle b_i x, x \rangle \in \R$, it follows that $\langle b_1 x, x \rangle
\geq 0$ and that $\langle b_2 x, x \rangle = 0$ for all $x$.  This shows that $f$
is positive.  To show $f$ is completely positive repeat the argument, replacing $B$
with $M_n(B)$. 

To show that $f$ contractive, we compute
\begin{align*}
\| b_1 + i b_2 \| 
&= \sup_{\|x_1 + ix_2 \| = 1} \|(b_1 + ib_2)(x_1 + ix_2) \|     \\
&\geq \sup_{\|x_1\| = 1}  \|(b_1 + ib_2)(x_1) \|               \\
&\geq \sup_{\|x_1\| = 1}  \|b_1(x_1) \|
= \|b_1\| \; .
\end{align*}
\end{proof}

\begin{proof} [Proof of Proposition~\ref{KKequiv}]
Let $\top$ denote the complex Toeplitz algebra generated by a single isometry $S$.  
Then as in Section~V.1 of \cite{Davidson} or Section~11.2 of \cite{Wegge}, $\top$
has an ideal isomorphic to the compact operators $\cp$ which is the kernel of the
homomorphism $\sigma \colon \top \rightarrow C(S^1, \C)$ which sends $S$ to the
identity function $z$.  Thus we have the short exact sequence \begin{equation}
\label{topseq} 0 \rightarrow \cp \rightarrow \top \xrightarrow{\sigma} C(S^1, \C)
\rightarrow 0 \; . \end{equation}

If we let $\topr$ be the real Toeplitz algebra generated by $S$, then $\topr$ contains 
the ideal $\topr \cap \cp = \cpr$, the algebra of real compact operators.  The
quotient is isomorphic to the real subalgebra of $C(S^1, \C)$ that is generated by
$z$, namely 
$$ C(S^1; \overline{\; \cdot \;})
= \{ f \in C(S^1, \C) \mid f(\overline{x}) = \overline{f(x)} \}    \; .$$ 
Therefore we have the short exact sequence 
\begin{equation} \label{toprseq} 0 \rightarrow \cpr \rightarrow
\topr \xrightarrow{\sigma} C(S^1; \overline{\; \cdot\;}) \rightarrow 0
\end{equation} as found in Section~1.5 of \cite{Schroder}.    

Now let $\pi \colon \topr \rightarrow \R$ be the composition of $\sigma$ with evaluation 
at $1 \in S^1$ and let $\topro$ be the kernel of $\pi$.  Then $\cpr$ is an ideal of
$\topro$ with quotient 
$$\{ f \in C(S^1; \overline{\; \cdot \;}) \mid f(1) = 0 \} \cong S^{-1} \; .$$   
Therefore we have our third short exact sequence
\begin{equation} \label{toproseq}
0 \rightarrow \cpr \rightarrow \topro \xrightarrow{\sigma} S^{-1} \rightarrow 0 \; . 
\end{equation}

Since $C(S^1, \C)$ is nuclear, the theorem of Choi and Effros recorded as
Theorem~15.8.3 of \cite{Black} implies
that Sequence~\ref{topseq} has a completely positive section.  Since
Sequence~\ref{topseq} is the complexification of Sequence~\ref{toprseq},
Lemma~\ref{cpositive} implies that Sequence~\ref{toprseq} also has a completely
positive section.  Finally, since $\topro$ is precisely the set of operators $T$ in
$\topr$ such that $\sigma(T)$ is in $S^{-1}$, this section restricts to a completely
positive section for Sequence~\ref{toproseq}.

Then by Proposition~2.5.6 of \cite{Schroder}, for any real separable C*-algebra $A$ 
the semisplit exact Sequence~\ref{toproseq} induces long exact sequences 
\begin{equation} \label{kk1}
\dots \rightarrow KK_{n+1}(A, \topro)
\rightarrow KK_n(A, S^{1,1}) 
\xrightarrow{\delta} KK_n(A, \R)
\rightarrow KK_n(A, \topro) \rightarrow \dots
\end{equation}
and
\begin{equation} \label{kk2}
\dots \rightarrow KK_n(\topro, A)
\rightarrow KK_n(\R, A)
\xrightarrow{\delta} KK_n(S^{1,1}, A)
\rightarrow KK_{n-1}(\topro, A) \rightarrow \dots
\end{equation}
where in both cases $\delta$ is given by the intersection product with an element 
$\delta \in KK_0(S^{1,1}, \R)$.  The element $\delta$ is one half of our $KK$-equivalence.

We claim that $KK_n(A, \topro) = KK_n(\topro, A) = 0$.  Using the long exact sequence 
on $K$-theory induced by
$$0 \rightarrow \topro \rightarrow \topr \xrightarrow{\pi} \R \rightarrow 0 \; ,$$
Schr\"oder shows (in Section 1.5 of \cite{Schroder}) that $K_*(\topro) = 0$ by showing 
that the homomorphism $\pi_*$ induced on $K$-theory has an inverse.  The same techniques 
immediately generalize to show
that the $KK$-theory homomorphisms $\pi_*$ and $\pi^*$ have inverses in the
sequences $$ \dots \rightarrow KK_n(A, \topro) \rightarrow KK_n(A, \topr)
\xrightarrow{\pi_*} KK_n(A, \R) \rightarrow KK_{n-1}(A, \topro) \rightarrow \dots
$$
and 
$$
\dots \rightarrow KK_{n+1}(\topro, A) \rightarrow KK_n(\R, A)
\xrightarrow{\pi^*} KK_n(\topr, A) \rightarrow KK_n(\topro, A) \rightarrow \dots
 \; .$$
Therefore $KK_n(A, \topro) = KK_n(\topro, A) = 0$.

It follows that $\delta$ is an isomorphism in the Sequences \ref{kk1} and \ref{kk2}.  
In particular, setting $A = \R$ in Sequence~\ref{kk1}, we have an isomorphism
$$KK_0(\R, S^{1,1}) \xrightarrow{\delta} KK_0(\R, \R) \cong \Z  \, .$$
To find the inverse $KK$-element to $\delta$, let $\ve \in KK_0(\R, S^{1,1})$ be the 
element such that $\ve \cap \delta = 1\sr \in KK_0(\R, \R)$.  It remains only to
show that $\delta \cap \ve = 1 \in KK_0(S^{1,1}, S^{1,1})$.

Applying Sequence~\ref{kk2} with $A = S^{1,1}$ gives an isomorphism
$$KK_0(\R, S^{1,1}) \rightarrow KK_0(S^{1,1}, S^{1,1})$$
and we conclude that $KK_0(S^{1,1}, S^{1,1}) \cong \Z$.
Now, we note that $\delta \cap \ve$ is an idempotent since 
$$(\delta \cap \ve)^2 = \delta \cap (\ve \cap \delta) \cap \ve  	
= \delta \cap 1 \cap \ve = \delta \cap \ve \; .$$ 
But $\delta \cap \ve$ is non-zero
since $\delta \cap \ve \cap \delta = \delta \neq 0$.  Therefore, $\delta 
\cap \ve$
must be the only non-zero idempotent in the ring $KK_0(S^{1,1}, S^{1,1}) \cong \Z$,
namely the identity. \end{proof}

We are now ready to prove the exactness of the sequences of 
Theorem~\ref{acyclic}.

\begin{proof}[Construction of Sequence~\ref{ouseq}] 
We start with the inclusion homomorphism $c \colon \R \rightarrow \C$ and construct 
the mapping cylinder and mapping cone.  See Section~2 of \cite{Schochet3} as a
basic reference for mapping cylinders and mapping cones.  In our case, the 
mapping cylinder is
\begin{align*} Zc &= \{ (x, f) \in \R \oplus I\C \mid f(0) = c(x) \}   \\
 &\cong \{ f \in I\C \mid f(0) \in \R \subset \C \}   \\
\intertext{and the mapping cone is}
Cc &= \{ (x, f) \in \R \oplus I\C \mid f(0) = c(x), f(1) = 0 \} \\
  &\cong \{ f \in I\C \mid f(0) \in \R, f(1) = 0 \} \\
 & \cong \{ f \in C_0(\R, \C) \mid f(-z) = \overline{f(z)} \}   \\
  & \cong S^{-1}
\end{align*}
Then there is a short exact sequence
\begin{equation} 
0 \rightarrow Cc \rightarrow Zc \rightarrow \C \rightarrow 0\, .
\end{equation}  
Since $Zc$ is homotopy equivalent to $\R$ and $Cc$ is isomorphic to $S^{-1}$, we
have up to homotopy equivalence a sequence
\begin{equation} \label{ouseq2}
0 \rightarrow S^{-1} \rightarrow \R \xrightarrow{c} \C \rightarrow 0
\end{equation}
which, upon tensoring on the right by $A$ and applying $K$-theory, 
gives rise to a long exact sequence
\begin{equation} \label{ouseq2'} 
\cdots \longrightarrow  KO_n(A) \xrightarrow{\alpha} 
KO_{n+1}(A) \xrightarrow{c} 
KU_{n+1}(A) \xrightarrow{\beta}
KO_{n-1}(A) \rightarrow \cdots \, . \end{equation}

Let $\alpha$ be the element of $KK_0(S^{-1}, \R) \cong KK_{1}(\R, \R)$ determined by 
the homomorphism $S^{-1} \rightarrow \R$ so that the $K$-theory homomorphism
$\alpha$ in the sequence above is realized by multiplication by the element
$\alpha$.  Now, the index map in the long exact sequence above is induced by a
homomorphism $S\C \rightarrow S^{-1}$ (see Proposition~2.3 of \cite{Schochet3}). 
Let $\beta \in KK(S\C, S^{-1}) \cong KK_{-2}(\C, \R)$ be the $KK$-element 
corresponding to this homomorphism.

We will show that $\alpha = \eta\so \in KK_1(\R, \R)$ and 
$\beta = \pm r \beta\su^{-1} \in KK_{-2}(\C, \R)$.  The first is easy since 
there are only two elements in $KK_1(\R, \R) \cong \Z_2$.  The element 
$\alpha$ cannot be zero since then Sequence~\ref{ouseq2'} would not be 
exact for $A = \R$.

To show that $\beta = \pm r \beta\su^{-1}$, we must first compute $KK_{-2}(\C, \R)$.  
For this, we use the long exact sequence induced by Sequence~\ref{ouseq2} on the
functor $KK(- , \R)$, $$KK_{-1}(\R, \R) \leftarrow KK_{-2}(\R, \R) \leftarrow
KK_{-2}(\C, \R) \leftarrow KK_0(\R, \R) \leftarrow KK_{-1}(\R, \R) \, $$
according to Proposition~2.5.4 of \cite{Schroder}.
Since $KK_{-2}(\R, \R) \cong KK_{-1}(\R, \R) = 0$ it follows that there is an 
isomorphism
$KK_{0}(\R, \R) \rightarrow KK_{-2}(\C, \R)$ given by the intersection product with 
$\beta$.  Therefore, $\beta$ is a generator of $KK_{-2}(\C, \R) \cong \Z$.  

Now, the inverse Bott element $\beta\su^{-1}$ is a generator of  $KK_{-2}(\C, \C) \cong \Z$ , 
so it remains only to prove that $r$ is a generator of $KK_0(\C, \R)
\cong \Z$.  For this, it is enough to show that $r \colon K_4(\C) \rightarrow
K_4(\R)$ is an isomorphism $\Z \rightarrow \Z$.  We use the short exact sequence 
$$0 \rightarrow K_4(\R) \xrightarrow{c} K_4(\C) \rightarrow K_4(\R)
\rightarrow 0 $$ 
induced by Sequence~\ref{ouseq2} and the fact that
$K_2(\R) = \Z_2$ to deduce that $c \colon K_4(\R) \rightarrow K_4(\C)$ is
multiplication by $\pm 2$ from $\Z$ to $\Z$.  Then the relation $cr = 2$ implies
that $r_4$ is an isomorphism.  \end{proof}

\begin{proof}[Construction of Sequence~\ref{otseq}]
To develop the third exact sequence, we follow the same procedure as above, 
starting with the homomorphism $\varepsilon \colon \R \rightarrow T$.

The mapping cylinder of $\ve$ is
\begin{align*}
Z\varepsilon 
	&= \{ (x, f) \in \R \oplus IT \mid f(0) = \varepsilon(x) \} \\
	&\cong \{ f \in I^2\C \mid 
	 f(s,0) = \overline{f(s,1)} \; \text{for all $s \in I$ and} \;
	 f(0,t_1) = f(0,t_2) \in \R \; \text{for all $t_1, t_2 \in I$} \}  \\
\intertext{and the mapping cone is}
C\varepsilon &= \{ (x, f) \in \R \oplus IT \mid f(0) = \varepsilon(x) \; \text{and} 
	\; f(1) = 0 \} \\
 &\cong \{ f \in I^2\C \mid 
	 f(s,0) = \overline{f(s,1)} \; \text{for all $s \in I$,} \\   
& \hspace{3cm} f(0,t_1) = f(0,t_2) \in \R \;
			 	\text{for all $t_1, t_2 \in I$, and} \\
& \hspace{3cm} f(1,t) = 0 \;	\text{for all $t \in I$} \}    \\
 &\cong \{ f \in I^2\C \mid	
	f(s,0) = f(s,1) = f(1,t) = 0 \; \text{and} \; 
	f(0,t) = \overline{f(0,1-t)} \}
\end{align*}
Call this last algebra $B$.  We will 
prove that $B$ is isomorphic to $S^{-2}$.  Recall that in previous proof, we
found that  
\begin{align*}
S^{-1} &\cong \{f \in I\C \mid f(0) \in \R, f(1) = 0 \} \, ; \\
\intertext{thus}
S^{-2} &\cong \{f \in I\C \mid f(0) \in \R, f(1) = 0 \}^{\otimes 2} \\
       &\cong \{f \in I^2 (\C \otimes\r \C) \mid 
	  f(1,t) = f(s,1) = 0 \; \text{for all $s$ and $t \in I$,}    \\
      & \hspace{4cm}
	  f(0,t) \in \R \otimes \C \; \text{for all $t$ in $I$, and} \\  
      & \hspace{4cm}
	f(s,0) \in \C \otimes \R \; \text{for all $s$ in $I$}
				\} \, .
\end{align*}

Now, $\C \oplus \C$ is isomorphic to $\C \otimes\r \C$ via the homomorphism
$$\alpha(x,y) = \tfrac{1}{2}(x + y) \otimes 1 + \tfrac{1}{2}(x - y)i \otimes i\, .$$ 
(This isomorphism is an isomorphism of complex C*-algebras when $\C \otimes\r \C$ 
receives the complex structure based on its first factor.)  Under this isomorphism
$\R \otimes\r \C$ is identified with $\{(x,\overline{x}) \mid x \in \C \} \subset \C
\oplus \C $ and $\C \otimes\r \R$ is identified with  $\{(x,x) \mid x \in \C \}$.

Hence $S^{-2}$ is isomorphic to the subalgebra $B'$ of $I^2(\C \oplus \C) 
	\cong I^2\C \oplus I^2\C$ consisting of pairs $(f, g)$ satisfying
\begin{align*}
f(s,1) &= f(1,t) = g(s,1) = g(1,t) = 0 \\ 
f(0,t) &= \overline{g(0,t)}, \text{and} \\
f(s,0) &= g(s,0) \; .
\end{align*}
Finally, the pasting map $(f, g) \mapsto h$ defined by
$$h(s,t) = 
\begin{cases}
	g(s,1-2t), &\text{if $t \leq \tfrac{1}{2}$.} \\
	f(s, 2t-1), &\text{if $t \geq \tfrac{1}{2}$} 
\end{cases} $$
gives an isomorphism from $B'$ to $B$.  We have shown that $C\varepsilon$ is isomorphic 
to $S^{-2}$.  

Then the short exact sequence 
$$0 \rightarrow C\varepsilon \rightarrow Z\varepsilon \rightarrow T \rightarrow 0$$
becomes, up to homotopy, the sequence
\begin{equation} \label{otseq2}
0 \rightarrow S^{-2} \rightarrow \R \xrightarrow{\varepsilon} T \rightarrow 0 
\end{equation}
which gives rise to a long exact sequence
\begin{equation} \label{otseq2'}
\cdots \longrightarrow  KO_n(A) \xrightarrow{\alpha} 
KO_{n+2}(A) \xrightarrow{\varepsilon} 
KT_{n+2}(A) \xrightarrow{\beta}
KO_{n-1}(A) \rightarrow \cdots \, . 
\end{equation}
Let $\alpha \in KK_2(\R, \R)$ and $\beta \in KK_{-3}(T, \R)$ be the $KK$-elements 
that implement the homomorphisms in the above sequence.  

We will show that $\alpha = \eta\so^2$ and $\beta = \pm\tau \beta\st^{-1}$.  The 
element $\alpha$ cannot be zero, since otherwise Sequence~\ref{otseq2'} above would
not be exact for $A=\R$.  Since $\eta\so^2$ is the only non-zero element of
$KK_2(\R, \R) \cong \Z_2$, we have $\alpha = \eta\so^2$.

Sequence~\ref{otseq2} gives rise to a long exact sequence
$$KK_{-1}(\R, \R) \leftarrow KK_{-3}(\R, \R) \leftarrow KK_{-3}(T, \R) 
	\xleftarrow{\beta} KK_0(\R, \R) \leftarrow KK_{-2}(\R, \R) \, .$$
Since $KK_{-2}(\R, \R) = KK_{-3}(\R, \R) = 0$ and since $KK_0(\R, \R) \cong \Z$, 
this implies that $KK_{-3}(T, \R) \cong \Z$ and the element $\beta$ is a generator.  

In the proof of the exactness of Sequence~\ref{ouseq}, we discovered that $r$ is an 
isomorphism from $KU_4(\R) \cong \Z$ to $KO_4(\R) \cong \Z$.  This isomorphism
passes through $KT_{3}(\R) \cong \Z$ via the factorization $r = \tau \gamma$.  So
$\tau$ must be a generator of $KK_1(T, \R) \cong \Z$.  Since $\beta\st$ is a
$KK$-equivalence, this implies $\tau \beta\st^{-1}$ is a generator of $KK_{-3}(T,
\R)$.

 \end{proof}


\section{Free Geometric Resolutions} \label{fgr}

In this section we show that given an arbitrary unital C*-algebra $A$, we can form 
a free geometric resolution of $K\crt(A)$.  In the first section, we compute the
united $K$-theory of the real C*-algebras $\R$, $\C$, and $T$.  These {\it
CRT}-modules give us geometric realizations of singly-generated free {\it
CRT}-modules and form the basic building blocks for the free geometric resolutions which
will be obtained in Section~\ref{geosurjection}


\subsection{Geometrically Realized Free {\it CRT}-Modules} \label{computation}

The united $K$-theory of the real C*-algebras $\R$, $\C$, and $T$ are given in 
Tables~\ref{crtR}, \ref{crtC}, and \ref{crtT}.  These tables give the groups
$KO_n(A)$, $KU_n(A)$, and $KT_n(A)$ up to isomorphism as well as the operations
$c_n$, $r_n$, $\varepsilon_n$, $\zeta_n$, $(\psi\su)_n$, $(\psi\st)_n$, $\gamma_n$,
and $\tau_n$ which are written as multiplication by a certain number or matrix in
terms of fixed ordered generators for each of the groups.  These preferred
generators are shown in the table for $K\crt(\R)$.  For $K\crt(\C)$ and $K\crt(T)$
the preferred generators are described in the course of the computation in the text.
 These preferred generators are chosen to be stable under the periodicity
isomorphisms.  The $KO_*(\R)$-, $KU_*(\R)$-, and $KT_*(\R)$-module structures are
given implicitly by the relations $\eta\so = \tau \varepsilon$, $\eta\st = \gamma
\beta\su \zeta$, $\xi = \tau \gamma \beta\su^2 c$, and $\omega = \beta\st \gamma
\zeta$.  These last two relations have not been proven in general yet.  However, our
calculation of $K\crt(\R)$ will imply that they hold for $A=\R$ and, following the
computation, we use this to prove the general case.

\begin{table}[h]
\caption{$K\crt(\R)$} \label{crtR}
$$\begin{array}{|c|c|c|c|c|c|c|c|c|c|}  
\hline \hline  
n & 0 & 1 & 2 & 3 & 4 & 5 & 6 & 7 & 8 \\
\hline  \hline
KO_n
& \Z \cdot 1\so & \Z_2 \cdot \eta\so & \Z_2 \cdot \eta\so^2 & 0 
	& \Z \cdot \xi & 0 & 0 & 0 & \Z \cdot \beta\so \\
\hline  
KU_n 
& \Z \cdot 1\su & 0 & \Z \cdot \beta\su & 0 & \Z \cdot \beta\su ^2 & 0 
	& \Z \cdot \beta\su ^3 & 0 & \Z \cdot \beta\su^4  \\
\hline  
KT_n 
& \Z \cdot 1\st & \Z_2 \cdot \eta\st & 0 & \Z \cdot \omega & \Z \cdot \beta\st 
     & \Z_2 \cdot \beta\st \eta\st & 0 & \Z \cdot \beta\st \omega 
		& \Z \cdot \beta\st^2 \\
\hline \hline
c_n & 1 & 0 & 0 & 0 & 2 & 0 & 0 & 0 & 1    \\
\hline
r_n & 2 & 0 & 1 & 0 & 1 & 0 & 0 & 0 & 2      \\
\hline
\varepsilon_n & 1 & 1 & 0 & 0 & 2 & 0 & 0 & 0 & 1   \\
\hline
\zeta_n &  1 & 0 & 0 & 0 &  1 & 0 & 0 & 0 & 1      \\
\hline
(\psi\su)_n & 1 & 0 & -1 & 0 & 1 & 0 & -1 & 0 & 1    \\
\hline
(\psi\st)_n & 1 & 1 & 0 & -1 & 1 & 1 & 0 & -1 & 1   \\
\hline
\gamma_n & 1 & 0 & 1 & 0 & 1 & 0 & 1 & 0 & 1     \\
\hline
\tau_n & 1 & 1 & 0 & 1 & 0 & 0 & 0 & 2 & 1        \\
\hline \hline
\end{array}$$
\end{table}

\begin{table} 
\caption{$K\crt(\C)$}  \label{crtC}
$$\begin{array}{|c|c|c|c|c|c|c|c|c|c|}  
\hline  \hline 
n & \makebox[1cm][c]{0} & \makebox[1cm][c]{1} & 
\makebox[1cm][c]{2} & \makebox[1cm][c]{3} 
& \makebox[1cm][c]{4} & \makebox[1cm][c]{5} 
& \makebox[1cm][c]{6} & \makebox[1cm][c]{7} 
& \makebox[1cm][c]{8} \\
\hline  \hline
KO_n 
& \Z & 0 & \Z & 0 & \Z & 0 & \Z & 0 & \Z  \\
\hline  
KU_n 
& \Z \oplus \Z & 0 & \Z \oplus \Z & 0 & \Z \oplus \Z & 0 
	& \Z \oplus \Z & 0 & \Z \oplus \Z   \\
\hline  
KT_n 
& \Z & \Z & \Z & \Z & \Z & \Z & \Z & \Z & \Z \\
\hline \hline
c_n & \smv{1}{1} & 0 & \smv{-1}{1} & 0 & \smv{1}{1} & 0 & \smv{-1}{1} & 0 & \smv{1}{1} 
	   \\
\hline
r_n & \smh{1}{1} & 0 & \smh{-1}{1} & 0 & \smh{1}{1} & 0 
	& \smh{-1}{1} & 0 & \smh{1}{1}      \\
\hline
\varepsilon_n & 1 & 0 & -1 & 0 & 1 & 0 & -1 & 0 & 1   \\
\hline
\zeta_n &  \smv{1}{1} & 0 & \smv{1}{-1} & 0 & \smv{1}{1} & 0 & \smv{1}{-1} & 0 
	& \smv{1}{1}      \\
\hline
(\psi\su)_n & \sm{0}{1}{1}{0} & 0 & \sm{0}{-1}{-1}{0} & 0 & 
	\sm{0}{1}{1}{0} & 0 & \sm{0}{-1}{-1}{0} & 0 & \sm{0}{1}{1}{0}    \\
\hline
(\psi\st)_n & 1 & -1 & 1 & -1 & 1 & -1 & 1 & -1 & 1   \\
\hline
\gamma_n & \smh{1}{1} & 0 & \smh{1}{-1} & 0 & 
	\smh{1}{1} & 0 & \smh{1}{-1} & 0 & \smh{1}{1}     \\
\hline
\tau_n & 0 & -1 & 0 & 1 & 0 & -1 & 0 & 1 & 0        \\
\hline \hline
\end{array}$$
\end{table}

\begin{table}  
\caption{$K\crt(T)$} \label{crtT}
$$\begin{array}{|c|c|c|c|c|c|c|c|c|c|}  
\hline  \hline 
n & \makebox[1cm][c]{0} & \makebox[1cm][c]{1} & 
\makebox[1cm][c]{2} & \makebox[1cm][c]{3} 
& \makebox[1cm][c]{4} & \makebox[1cm][c]{5} 
& \makebox[1cm][c]{6} & \makebox[1cm][c]{7} 
& \makebox[1cm][c]{8} \\
\hline  \hline
KO_n 
& \Z & \Z_2 & 0 & \Z & \Z & \Z_2 & 0 & \Z & \Z  \\
\hline  
KU_n 
& \Z & \Z & \Z & \Z & \Z & \Z & \Z & \Z & \Z   \\
\hline  
KT_n 
& \Z \oplus \Z_2 & \Z_2 & \Z & \Z \oplus \Z & \Z \oplus \Z_2 & \Z_2 & \Z 
& \Z \oplus \Z & \Z \oplus \Z_2  \\
\hline \hline
c_n & 1 & 0 & 0 & 2 & 1 & 0 & 0 & 2 & 1    \\
\hline
r_n & 2 & 1 & 0 & 1 & 2 & 1 & 0 & 1 & 2      \\
\hline
\varepsilon_n & \smv{1}{0} & 1 & 0 & \smv{2}{1} & \smv{1}{1} & 1 & 0 
	& \smv{2}{1} & \smv{1}{0}   \\
\hline
\zeta_n &  \smh{1}{0} & 0 & 0 & \smh{1}{0} & 
	\smh{1}{0} & 0 & 0 & \smh{1}{0} & \smh{1}{0}    \\
\hline
(\psi\su)_n & 1 & -1 & -1 & 1 & 1 & -1 & -1 & 1 & 1   \\
\hline
(\psi\st)_n & \sm{1}{0}{0}{1} & 1 & -1 & \sm{1}{0}{1}{-1}  
	   & \sm{1}{0}{0}{1} & 1 & -1 & \sm{1}{0}{1}{-1} & \sm{1}{0}{0}{1}  \\
\hline
\gamma_n & \smv{0}{1} & \smv{0}{1} & 1 & 1 & \smv{0}{1} 
			& \smv{0}{1} & 1 & 1 & \smv{0}{1}    \\
\hline
\tau_n & \smh{1}{1} & 0 & 1 & \smh{-1}{2} 
		& \smh{0}{1} & 0 & 1 & \smh{-1}{2} & \smh{1}{1}   \\
\hline \hline
\end{array}$$
\end{table}

Let $F(b, n, \R)$ (respectively $F(b, n, \C)$ and $F(b, n, T)$) denote the free 
{\it CRT}-module with a single generator $b$ in the real (respectively complex and
self-conjugate) part in degree $n$ (see Section~2 of \cite{Bou}).  Any free {\it
CRT}-module can be obtained as a direct sum of these three monogenic free {\it
CRT}-modules.

By direct comparison with Paragraph~2.4 of \cite{Bou}, we observe from our
computations that 
$K\crt(\R)$ is a free {\it CRT}-module with single generator $1\so \in KO_0(\R)$; 
$K\crt(\C)$ is a free {\it CRT}-module with generator $\kappa_1 = \smv{1}{0} \in
KU_0(\C)$;  and $K\crt(T)$ is a free {\it CRT}-module with generator $\chi =
\smv{-1}{0} \in KT_{-1}(T)$.  Therefore, any free {\it CRT}-module can be realized
as $K\crt(F)$, where $F$ is a direct sum of suspensions of the algebras $\R$, $\C$,
and $T$. 

\begin{proof}[Computation of Table~\ref{crtR}]
We have already established the group and module structures of $KO_*(\R)$, $KU_*(\R)$, 
and $KT_*(\R)$.  For the computation of $c$, $r$, and $\psi\su$, we refer to 
\cite{Karoubi} (III.5.19) or \cite{Bou}.  Alternatively, the values of $c$ and $r$ 
can be computed directly using Sequence~\ref{ouseq} and $\psi\su$ using the
relation $cr = 1+ \psi\su$.  We computed $\zeta$ and $\gamma$ in
Section~\ref{computation} using Sequence~\ref{utseq}.  The remaining homomorphisms
$\varepsilon$, $\psi\st$, and $\tau$ can be easily determined by the constraints $c
= \zeta \varepsilon$, $\eta\so = \tau \varepsilon$, $\varepsilon r \zeta = 1 +
\psi\st$, and $r = \tau \gamma$.  \end{proof}

At this point, we pause to prove two more of the relations stated in 
Proposition~\ref{crtobject} using the results of the computation of $K\crt(\R)$. 
These relations will be used in computing $K\crt(\C)$.  

\begin{proof}[Continuation of proof of Proposition~\ref{crtobject}]
First we prove $\xi = r \beta\su^2 c$.  Let $x \in KO_*(A)$.  Then using the relations 
$r = \tau \gamma$, $c = \zeta \varepsilon$ (from Proposition~\ref{crtobject});
Lemma~\ref{mult2}, parts (3), (5), and (10); and $r \beta\su^2
c(1\so) = \xi$ from Table~\ref{crtR} we have \begin{align*}  r \beta\su^2 c(x) 
		&= r \beta\su^2 c(1\so \cdot\so x)   \\
		&= r \beta\su^2 (1\su \cdot\su c(x))   \\
		&= \tau \gamma (\beta\su^2 \cdot\su \zeta \varepsilon(x))    \\
		&= \tau (\gamma \beta\su^2 \cdot\st \varepsilon(x))   \\
		&= \tau \gamma \beta\su^2 \cdot\so (x)       \\
		&= r \beta\su^2c(1\so) \cdot\so x   \\
		& = \xi \cdot\so x   \; .\end{align*}

Second, we prove $\omega = \beta\st \gamma \zeta$.  Let $z \in KT_*(A)$.  Then using 
$\gamma (1\su) = \beta\st^{-1} \omega$ (from Table~\ref{crtR}) and
Lemma~\ref{mult2} part (3), we have \begin{align*}
\beta\st \gamma \zeta(z) &= \beta\st \gamma (1\su \cdot\su \zeta(z))  \\
			&= \beta\st (\gamma(1\su) \cdot\st z)   \\
			&= \beta\st (\beta\st^{-1} \omega \cdot\st z)  \\
			&= \omega \cdot\st z \; . \end{align*}
\end{proof}

\begin{proof}[Computation of Table~\ref{crtC}]
We designate the preferred generators of $KO_n(\C)$ to be those elements which
correspond to the preferred generators of $KU_n(\R)$.

Now, as in the proof of Thoerem~\ref{acyclic}, $\C \oplus \C$ is isomorphic to 
$\C \otimes \sr \C$ via the isomorphism
$$\alpha(\lambda_1, \lambda_2) = 
\tfrac{1}{2}(\lambda_1 + \lambda_2) \otimes 1 + 
\tfrac{1}{2}(\lambda_1 - \lambda_2)i \otimes i \; .$$
Thus, $KU_n(\C) = K_n(\C \otimes \sr \C) \cong K_n(\C \oplus \C)$ is isomorphic to 
$\Z \oplus \Z$ if $n$ is even and is 0 otherwise.  We define the classes 
\begin{align*} 
\kappa_1 &= \alpha_*[(1,0)] 
	= [\tfrac{1}{2}(1 \otimes 1) + \tfrac{1}{2}(i \otimes i)]   \\
\kappa_2 &= \alpha_*[(0,1)]
	= [\tfrac{1}{2}(1 \otimes 1) - \tfrac{1}{2}(i \otimes i)] 
\end{align*}
in $KU_0(\C)$.  Then the preferred generators of $KU_{2n}(\C)$ are set to be 
$\beta\su^n \cdot \kappa_1$ and $\beta\su^n \cdot \kappa_2$ in that order.

Now the complexification operation $c \colon KO_n(\C) \rightarrow KU_n(\C)$ is computed 
by studying the composition 
$$\C \cong \R \otimes \sr \C \hookrightarrow \C \otimes \sr \C \cong \C \oplus \C$$ 
which is given by $\lambda \mapsto (\overline{\lambda}, \lambda).$  Similarly, the
operation $\psi\su \colon KU_n(\C) \rightarrow KU_n(\C)$ is computed by studying the
composition 
$$\C \oplus \C \cong \C \otimes \sr \C \xrightarrow{~\psi\su \otimes 1~}
\C \otimes \sr \C \cong \C \oplus \C$$  
which is given by $(\lambda_1, \lambda_2)
\mapsto (\overline{\lambda_2}, \overline{\lambda_1})$.

To compute $KT_*(\C)$, we use the exact sequence 
$$0 \rightarrow KT_n(\C) \mapis{\zeta} 
KU_n(\C) \xrightarrow{1 - \psi\su}
KU_n(\C) \mapis{\gamma}
KT_{n-1}(\C) \rightarrow 0$$
for $n$ even. 

Now $1 - (\psi\su)_n = \sm{1}{-1}{-1}{1}$ for $n \equiv 0 \pmod 4$ and 
$1 - (\psi\su)_n = \sm{1}{1}{1}{1}$ for $n \equiv 2 \pmod 4$.  So for any even $n$,
the kernel and the cokernel of $1 - (\psi\su)_n$ are both isomorphic to $\Z$.  Thus
$KT_n(\C) \cong \Z$ for all $n$.  For each $n$, choose the preferred generators of
$KT_n(\C)$ such that $\gamma_n$ and $\zeta_n$ are as described in the table.  Since
$\zeta \beta\st = \beta\su^2 \zeta$ and $\gamma \beta\su^2 = \beta\st \gamma$, the
collection of these preferred generators is stable under the periodicity isomorphism
$\beta\st$.

Finally, the operations $r$, $\varepsilon$, $\psi\st$, and $\tau$ are determined by the 
relations $cr=1+\psi\su$, $c= \zeta \varepsilon$, $\ve r \zeta = 1 + \psi\st$, and
$r=\tau \gamma$.     \end{proof} 

\begin{proof}[Computation of Table~\ref{crtT}]
The structure of real $K$-theory is immediate, since $KO_n(T) \cong KT_n(\R)$.  We 
specify generators of $KO_*(T)$ which correspond to the previously specified
generators of $KT_*(\R)$.  The $KO_*(\R)$-module structure of $KO_*(T)$ is inferred
from the $KT_*(\R)$-module structure on $KT_*(\R)$ via the inclusion $\ve 
\colon \R \hookrightarrow T$.  For example, for the element $\xi \in KO_4(\R)$ 
and any $z \in KT_k(\R) = KO_k(T)$ we have
$$\xi \cdot z = \ve(\xi) \cdot z = 2 \beta\st \cdot z \; ,$$
a fact which we shall use shortly. 

For all $n$, we have $KU_n(T) \cong K_n(\C \otimes T) \cong KT_n(\C) \cong 
\Z$.  Then use the sequence
\begin{equation} \label{ouseqc} 
KO_{n-1}(T) \mapis{\eta\so} 
KO_n(T) \mapis{c} 
KU_n(T) \mapis{r\beta\su^{-1}}
KO_{n-2}(T) \mapis{\eta\so} 
KO_{n-1}(T) \end{equation}
from Theorem~\ref{acyclic} to compute $c$ and $r$.

Let $y_0 = c({1}\st)$ be the preferred generator of $KU_0(T)$ and let $y_1$ be 
the generator of $KU_1(T)$ such that $c({\omega}) = 2 \beta\su \cdot y_1$. 
Then we take $\beta\su^n \cdot y_0$ to be the preferred generator of $KU_{2n}(T)$
and we take $\beta\su^n \cdot y_1$ to be the preferred generator of $KU_{2n+1}(T)$. 

From these choices of generators and from Sequence~\ref{ouseqc}, along with the
relations $\psi\su c = c$, $\psi\su \beta\su = \beta\su \psi\su$, $rc = 2$ and $cr = 1 + \psi\su$, we can mostly compute the operations of
$r$, $c$, and $\psi\su$.  However, the values of $r_4$, $c_4$, $r_7$, and $c_7$ are
only determined up to sign.  The following calculations show that $r_4 = 2$ and $r_7
= 1$ and it follows that $c_4 = 1$ and $r_7 = 2$:
\begin{align*}
r(\beta\su^2 \cdot y_0) &= r(\beta\su^2 \cdot c(1\st)) \\
	&= r(\beta\su^2) \cdot 1\st \\
	&= \xi \cdot 1\st  \\
	&= 2 \beta\st \\
\intertext{and}
2 r(\beta\su^3 \cdot y_1) &= r(\beta\su^2 \cdot 2 \beta\su y) \\
	&= r(\beta\su^2 \cdot c(\omega) ) \\
	&= r(\beta\su^2) \cdot \omega \\
	&= \xi \cdot \omega \\
	&= 2 \beta\st \omega 
\end{align*}

The groups $KT_*(T)$ can be computed immediately from the sequence
\begin{equation} \label{utseqt} 
KU_{n+1}(T) \xrightarrow{1 - \psi\su}
KU_{n+1}(T) \mapis{\gamma}
KT_{n}(T) \mapis{\zeta}
KU_n(T) \xrightarrow{1 - \psi\su}
KU_n(T) \end{equation}
using what we know about $KU_*(T)$ and $\psi\su$.  Choose generators for $KT_n(T)$ so 
that $\gamma$ and $\zeta$ are as specified in the table (this choice is not
unique).  

To compute $\varepsilon$ we use the sequence
\begin{equation}  \label{otseqt}
KO_{n-2}(T) \mapis{\eta\so^2} 
KO_n(T) \mapis{\varepsilon} 
KT_n(T) \mapis{\tau \beta\st^{-1}}
KO_{n-3}(T) \mapis{\eta\so^2} 
KO_{n-1}(T) \; . \end{equation}
Note that in this sequence, $\eta\so^2 = 0$ for all $n$.  It is immediate that 
$\varepsilon_1 = \varepsilon_5 = 1$.

Looking at Sequence~\ref{otseqt} with $n = 0$, we see that the image of $\varepsilon_0$ 
is a subgroup of index 2 in $KT_0(T) \cong \Z \oplus \Z_2$.  Since $\zeta_0
\varepsilon_0 = c_0 = 1$, we conclude that $\varepsilon_0$ is either $\smv{1}{0}$ or
$\smv{1}{1}$.  By changing the basis of $KT_0(T)$ we may assume that $\varepsilon_0
= \smv{1}{0}$.  (If a change is necessary, the new generators will be $\smv{1}{1}$
and $\smv{0}{1}$ given in terms of the old generators.  This basis change does not
disturb our previously computed representations of $\gamma$ and $\zeta$.)

From $\zeta_3 \varepsilon_3 = c_3 = 2$, we see that $\varepsilon_3$ is $\smv{2}{a}$ 
where $a$ is some integer.  But the exact sequence above tells us that the image of
$\varepsilon_3$ is a free summand of $KT_3(T) \cong \Z \oplus \Z$.  Thus $a$ must be
odd.  Again we may change the basis $KT_3(T)$ so that $\varepsilon_3 = \smv{2}{1}$. 
(The new generators will be $\smv{1}{\tfrac{a-1}{2}}$ and $\smv{0}{1}$.  Again this
change does not disturb our representations of $\gamma$ and $\zeta$.)

Using the same reasoning that we used in the computation of $\varepsilon_0$ above, we know 
that $\varepsilon_4$ is either $\smv{1}{0}$ or $\smv{1}{1}$.  However at this point
we do not have the freedom to change generators.  (The generators of $KT_4(T)$ are to
be the suspensions of the generators of $KT_0(T)$.)  Instead,  we compute 
$\tau_0 = \smh{1}{1}$ using the relations $\eta\so = \tau \varepsilon$ and $r = \tau
\gamma$.  Then we use the relation $\tau \beta\st^{-1} \varepsilon = 0$ from
Sequence~\ref{otseq} to see that $\varepsilon_4 = \smv{1}{1}$.

The computation of $\varepsilon_7$ is similar to that of $\varepsilon_4$.  First we use 
the relation $\zeta_7 \varepsilon_7 = c_7 = 2$ and the exact Sequence~\ref{otseqt}
to deduce that $\varepsilon_7 = \smv{2}{a}$ for some odd $a$.  Independently, we
compute $\tau_3 = \smh{-1}{2}$ using the relations $\eta\so = \tau \varepsilon$ and
$r = \tau \gamma$.  Finally, we use the relation $\tau \beta\st^{-1} \varepsilon =
0$ to see that $a = 1$.

It remains to compute $\psi\st$ and $\tau$, which can easily be done using the relations 
$\ve r \zeta = 1 + \psi\st$, $r = \tau \gamma$, and $\eta\so = \tau \ve$.
\end{proof}  


\subsection{Geometrically Realized Surjections and Resolutions} \label{geosurjection} 

In \cite{Schochet}, Schochet produces a geometric resolution by first producing for any 
unital complex C*-algebra $A$ a homomorphism $\mu \colon F \rightarrow \cp(\mathcal{H})
\otimes A$ that is surjective on both $K_0$ and $K_1$ where $F$ is a 
C*-algebra such that $K_*(F)$ is free
(Lemma~3.1 of \cite{Schochet}).  Roughly, he used homomorphisms from $\R$ in order
to hit each projection (representing an element in $K_0(A)$) and a homomorphism from
$S\R$ in order to hit each unitary (representing an element in $K_1(A)$).  Then he
embedded $A$ in $\cp \otimes A$ in order to have room to combine the 
homomorphisms.  Once $\mu$ is consructed, the resolution is completed using 
a mapping cone construction.

In what follows we will also produce a geometric resolution by first producing a 
homomorphism which is surjective on united $K$-theory.  Now there are eight
distinct groups in real $K$-theory and only two of them have simple intrinsic descriptions
(projections for $K_0$ and unitaries for $K_1$).  Therefore, we have to resort to
suspensions, representing elements of $K_i(A)$ by projections in $S_iA$.  Therefore
the target of our homomorphism will be the stabilization of a certain eightfold
suspension of $A$ described below.

We will make essential use of Schochet's stabilization trick in order to have room to 
paste together the many homomorphisms that are necessary to hit the many elements
of all the different groups of united $K$-theory.  Finally, we will use a mapping
cone construction similar to Schochet's and the fact that $K\crt(A)$ has projective
dimension 1 in {\it CRT} to turn a geometric surjective homomorphism into a
geometric resolution of length 1.

Let $A$ be a real unital C*-algebra.  Let 
$$S^+A = \{ f \in C([0,1], A) \mid f(0) = f(1) \in \R \cdot 1_A \} $$ 
be the unitized suspension of $A$.  We let $\mathfrak{S}A$ denote the unital algebra 
obtained by repeating this process eight times:
$$\mathfrak{S}A := (S^+)^8A\; .$$
There is a split exact sequence
$$0 \rightarrow S^8A \rightarrow \mathfrak{S}A \rightarrow (S^+)^7 \R \rightarrow 0 \; $$
 yielding the direct sum decomposition
$$K\crt(\mathfrak{S}A) \cong K\crt(S^8A) \oplus K\crt((S^+)^7 \R) \; .$$ 
There is also an isomorphism $K\crt(A) \cong K\crt(S^8A)$ so $K\crt(\mathfrak{S}A)$ 
carries the same essential information as $K\crt(A)$.  In our geometric resolution,
we will replace $K\crt(A)$ by $K\crt(\mathcal{K}(\mathcal{H}) \otimes
\mathfrak{S}A)$.  The main theorems of this section are the following.

\begin{prop} \label{surjection} Let $A$ be a real unital C*-algebra.  Then there is a 
real C*-algebra $F$, a real Hilbert space $\mathcal{H}$, and a C*-algebra
homomorphism  $\mu \colon F \rightarrow \mathcal{K}(\mathcal{H}) \otimes
\mathfrak{S}A$  such that $K\crt(F)$ is a free {\it CRT}-module and the induced
homomorphism $$\mu_* \colon K\crt(F) \rightarrow K\crt(\mathcal{K}(\mathcal{H})
\otimes \mathfrak{S}A) \cong K\crt(\mathfrak{S}A)$$ is a surjection.
\end{prop}

\begin{thm} \label{resolution} Let $A$ be a real unital C*-algebra.  Then there are real
 C*-algebras $F_1$ and $F_0$, C*-algebra homomorphisms $\mu_1$ and $\mu_2$, and a
Hilbert space $\mathcal{H}$ making the following sequence exact $$0 \rightarrow
K\crt(F_1) \xrightarrow{(\mu_1)_*} K\crt(F_0) \xrightarrow{(\mu_0)_*}
K\crt(\mathcal{K}(\mathcal{H}) \otimes \mathfrak{S}A)
\rightarrow 0 \; .$$
Furthermore, $K\crt(F_0)$ and $K\crt(F_1)$ are free {\it CRT}-modules.
\end{thm}

We tackle the proof of Proposition~\ref{surjection} by considering the real, complex, 
and self-conjugate cases independently in the next lemmas.  Once it is proven, it
will be used with Lemma~\ref{crtobject} to prove Theorem~\ref{resolution}.

\begin{lemma} \label{realsurjection} Let $A$ be a real unital C*-algebra.  Then there is 
a real C*-algebra $F$, a real Hilbert space $\mathcal{H}$, and a C*-algebra
homomorphism  $\mu \colon F \rightarrow \mathcal{K}(\mathcal{H}) \otimes A$ 
such that $K\crt(F)$ is a free {\it CRT}-module and the induced homomorphism
$$\mu_* \colon KO_0(F) \rightarrow KO_0(\mathcal{K}(\mathcal{H}) \otimes A) 
		\cong KO_0(A)$$
is a surjection.
\end{lemma}

\begin{proof} 
Let $\{p_s\}_{s \in S}$ be a collection of projections $p_s \in
\mathcal{K}(\mathcal{H}_s) \otimes A$ such that $\{[p_s]\}_{s \in S}$ generates 
$KO_0(A)$.  Here $S$ is a possibly uncountable index set $S$, and for
each $s \in S$, $\mathcal{H}_s$ is a finite dimensional Hilbert space.  Then
construct homomorphisms $$\mu_s \colon \R \rightarrow \mathcal{K}(\mathcal{H}_s)
\otimes A$$  defined by $\mu_s(\lambda) = \lambda p_s$ so that $[p_s]$ is in the
image of $(\mu_s)_*$.

Now form the Hilbert space direct sum $\mathcal{H} = \bigoplus_{s \in S} \mathcal{H}_s$ 
and the C*-algebra direct sum $F = \bigoplus_{s \in S} \R$.  Then we patch together
the homomorphisms $\mu_s$ to form a homomorphism $$\mu \colon F \rightarrow
\mathcal{K}(\mathcal{H}) \otimes A $$ by defining $\mu$ to be the composition 
	$$\R \mapis{\mu_s} \mathcal{K}(\mathcal{H}_s) \otimes A 
	\hookrightarrow \mathcal{K}(\mathcal{H}) \otimes A$$
on each summand $\R$ of $F$.
\end{proof}

Recall that $\kappa_1 \in KU_0(\C)$ is the free generator of $K\scrt(\C)$.  From
Table~\ref{crtC} we have $r(\kappa_1) = 1\su$.

\begin{lemma} \label{complexprojectiontarget}
Let $A$ be a real unital C*-algebra.  Let $p$ be a projection in $\C \otimes A$.  Then 
there is a C*-algebra homomorphism
$\mu \colon \C \rightarrow M_2(A)$
such that the induced homomorphism 
$$\mu_* \colon KU_0(\C) \rightarrow KU_0(M_2(A)) \cong KU_0(A)$$ 
carries $\kappa_1$ to $[p]$ modulo the image of $\beta\su^{-1} c$. 
\end{lemma}

\begin{proof}
Define $\mu$ by $\mu(\lambda) = r(\lambda p) \in M_2(A)$.  Then
$$r(\mu_*(\kappa_1) - [p]) = \mu_*(r(\kappa_1)) - r[p] = \mu_*[1] - r[p] = 0 \; \; .$$  
Therefore $\mu_*(\kappa_1) - [p]$ is in the kernel of $r$ which is the image of 
$\beta\su^{-1} c$.
\end{proof}

\begin{lemma} \label{complexsurjection} Let $A$ be a real unital C*-algebra.  Then there 
is a real C*-algebra $F$, a real Hilbert space $\mathcal{H}$, and a C*-algebra
homomorphism  $\mu \colon F \rightarrow \mathcal{K}(\mathcal{H}) \otimes A$ 
such that $K\crt(F)$ is a free {\it CRT}-module and the induced homomorphism
$$\mu_* \colon KU_0(F) \rightarrow KU_0(\mathcal{K}(\mathcal{H}) \otimes A) 
			\cong KU_0(A)$$
is a surjection modulo the image of $\beta\su^{-1}c$.
\end{lemma}

\begin{proof}
Again, choose a family of projections $p_s \in \mathcal{K}(\mathcal{H}_s) \otimes \C
\otimes A$ such that the classes $[p_s]$ generate $KU_0(A)$.  By
Lemma~\ref{complexprojectiontarget}, we have homomorphisms
$$\mu_s \colon \C \rightarrow M_2(\mathcal{K}(\mathcal{H}_s) \otimes A) 
			\cong \mathcal{K}(\mathcal{H}'_s) \otimes A \; $$
such that $\mu_s(\kappa_1) = [p_s]$ modulo $\beta\su^{-1} c$.
Here $\mathcal{H}'_s$ is taken to be the Hilbert space tensor product 
$\R^2 \otimes \mathcal{H}_s$ so that 
$M_2(\R) \otimes \mathcal{K}(\mathcal{H}_s) \cong 
\mathcal{K}(\R^2) \otimes \mathcal{K}(\mathcal{H}_s) \cong
\mathcal{K}(\mathcal{H}'_s)$.

Then the family {$\mu_s$} is patched together as in the proof of 
Lemma~\ref{realsurjection} to form
$$\mu_s \colon F \rightarrow \mathcal{K}(\mathcal{H}) \otimes A$$
where  $F = \bigoplus_{s \in S} \C$ and 
$\mathcal{H} = \bigoplus_{s \in S} \mathcal{H}'_s$.
\end{proof}

We now turn to the self-conjugate case.  In this case the target of our homomorphism 
will not be a matrix algebra over $A$, but a matrix algebra over $S^+A$.  With
respect to the natural decomposition  	$$K\crt(S^+A) \cong K\crt(SA) \oplus
K\crt(\R) \; , $$ 
let $\pi \colon K\crt(S^+A) \rightarrow K\crt(SA)$ be the projection
onto the first summand.  For the source of our homomorphism, we will use the algebra
$T$.  Recall that $K\crt(T)$ is a free {\it CRT}-module generated by $\chi =
\smv{-1}{0} \in KT_{-1}(T)$.  The only other fact we need to remember about $\chi$
is that $\tau(\chi) = 1\st \in KO_0(T)$ (see Table~\ref{crtT}).

\begin{lemma} \label{selfconjugateprojectiontarget} Let $A$ be unital and let $p$ be 
a projection in $T \otimes A$.  Then there is a C*-algebra homomorphism  
$\mu \colon T \rightarrow M_{16}(S^+A)$
such that the composition 
$$\pi \circ \mu_* \colon KT_{-1}(T) \rightarrow KT_{-1}(M_{16}(S^+A)) 
\rightarrow KT_{-1}(SA) \cong KT_0(A)$$ 
carries $\chi$ to $[p]$ modulo the image of $\beta\st^{-1} \varepsilon$.  
\end{lemma}

\begin{proof} 
Fix a unitary $u$ in $M_8(A)$ such that
$$u 
\begin{pmatrix}  r(p(0)) & 0 & 0 & 0 \\ 0 & 0 & 0 & 0 \\
		0 & 0 & 1_2 - r(p(0)) & 0 \\ 0 & 0 & 0 & 0
\end{pmatrix}  
u^* = 
\begin{pmatrix}  1_2 & 0 & 0 & 0 \\ 0 & 0 & 0 & 0 \\
		0 & 0 & 0 & 0 \\ 0 & 0 & 0 & 0
\end{pmatrix} \; .$$
Recall from Section~\ref{trans} that for $i = 1, 2$, we have homomorphisms
$$\sigma_i \colon T \otimes A \rightarrow C(S^1, M_4(A)) \; $$
such that $\tau = (\sigma_1)_* - (\sigma_2)_*$.  Then let 
$P' \in C(S^1, M_8(A))$ be the 
projection defined by 
$$P' = u
\begin{pmatrix} \sigma_1(p) & 0 \\ 0 & \sigma_2(1 - p) \end{pmatrix}
u^* \; .$$
Note that $P'$ is in $M_8(S^+A)$ since 
$$P'(0) = P'(1) = u \sm{\sigma_1(p)(0)}{0}{0}{\sigma_2(1-p)(0)} u^*
	 = u \sm{r(p(0))}{0}{0}{1 - r(p(0))} u^* 
	= \sm{1_2}{0}{0}{0}  \; .$$
Finally, let $P = \sm{P'}{0}{0}{0} \in M_{16}(S^+A)$.  

We will construct a homomorphism $\mu \colon T \rightarrow M_{16}(S^+A)$ that sends 
the unit to $P$.  For this we first define the C*-algebra homomorphism $\mu' \colon
T \rightarrow M_{16}(\R)$ as follows.  Given $f \in T$, let $\zeta(f) = f(0) = x +
iy$.  Then define  $$\mu'(f) = 	\begin{pmatrix} 		 	x & 0 & 0 &
0 & -y & 0 & 0 & 0  \\ 			0 & x & 0 & 0 & 0 & -y & 0 & 0  \\
			0 & 0 & x & 0 & 0 & 0 & -y & 0  \\
			0 & 0 & 0 & x & 0 & 0 & 0 & -y   \\
			y & 0 & 0 & 0 & x & 0 & 0 & 0  \\
			0 & y & 0 & 0 & 0 & x & 0 & 0  \\
			0 & 0 & y & 0 & 0 & 0 & x & 0  \\
			0 & 0 & 0 & y & 0 & 0 & 0 & x  
		\end{pmatrix} $$

Now $\mu'(f)$ commutes with the projection $P$ for all $f$ so the rule
$$\mu(f) = \mu'(f) \cdot P$$
defines a C*-algebra homomorphism which sends the unit to $P$.

We claim that $\mu_*$ carries ${1\st} \in KO_0(T)$ to 
$\tau[p] + 2 \cdot 1\so \in KO_0(M_8(\R) \otimes S^+A) \cong KO_1(A) \oplus
KO_0(\R)$.  Indeed,  \begin{align*} 
\mu_*({1\st}) 
&= [P]    \\
&= [\sigma_1(p)] + [\sigma_2(1-p)]      \\
&= [\sigma_1(p)] + [\sigma_2(1)] - [\sigma_2(p)]     \\
&= \tau[p] + [r(1)] \\
&= \tau[p] + 2 \cdot 1\so \; .
\end{align*}

Then 
$$ \tau (\pi \circ \mu_*)(\chi) 
= (\pi \circ \mu_*) \tau (\chi)
= (\pi \circ \mu_*) (1\st)
= \pi(\tau[p] + 2\so)
= \tau[p] \; . $$

Therefore, $(\pi \circ \mu_*)(\chi - [p])$ is in the kernel of $\tau$ which is the image 
of $\beta\st^{-1} \varepsilon$, according to Sequence~\ref{otseq}.   
\end{proof}

\begin{lemma} \label{selfconjugatesurjection} Let $A$ be a real unital C*-algebra.  Then 
there is a C*-algebra $F$, a Hilbert space $\mathcal{H}$, and a C*-algebra
homomorphism  $\mu \colon F \rightarrow \mathcal{K}(\mathcal{H}) \otimes S^+A$ 
such that $K\crt(F)$ is a free {\it CRT}-module and the induced homomorphism
$$\mu_* \colon KT_{-1}(F) \rightarrow KT_{-1}(\mathcal{K}(\mathcal{H}) \otimes S^+A) 
\cong KT_{-1}(S^+A) \cong KT_0(A) \oplus KT_{-1}(\R)$$
is a surjection modulo the image of $\beta\st^{-1} \ve$.
\end{lemma}

\begin{proof} 
Let $\mu_0 \colon \R \rightarrow S^+A = \mathcal{K}(\mathcal{H}_0) \otimes S^+A$ be the 
unital inclusion where $\mathcal{H}_0$ is a one-dimensional real Hilbert space. 
The homomorphism $\mu_0$ is designed to insure that the image of $\mu_*$ includes
$KT_{-1}(\R) \subset KT_{-1}(S^+(A))$.  Now let $p_s \in \mathcal{K}(\mathcal{H}_s)
\otimes T \otimes A$ be a family of projections indexed by $S$ so that $\{ [p_s]
\}_{s \in S}$ and the image of $\beta\st^{-1} \ve$ generate $KT_0(A)$.    For each $s \in S$, let  
$$\mu_s \colon T
\longrightarrow M_{16}(\mathcal{K}(\mathcal{H}_s) \otimes S^+A)  \cong
\mathcal{K}(\mathcal{H}'_s) \otimes S^+A$$ 
be given according to
Lemma~\ref{selfconjugateprojectiontarget} such that $(\pi \circ (\mu_s)_*)\chi =
[p_s]$ modulo $\beta\st^{-1} \varepsilon$.

Then, taking $F = \R \oplus \bigoplus_{s \in S} T$ and $\mathcal{H} = \mathcal{H}_0 
\oplus \bigoplus_{s \in S} \mathcal{H}_s'$, we form 
$$\mu \colon F \longrightarrow \mathcal{K}(\mathcal{H}) \otimes S^+A$$ by patching 
together the homomorphisms $\mu_0$ and $\{ \mu_s \}_{s \in S}$.
\end{proof}

\begin{proof}[Proof of Proposition~\ref{surjection}]
Let $A$ be a given real unital C*-algebra.  Using the now familiar patching technique, 
it suffices to produce enough maps of the form
$$\mu_s \colon F_s \rightarrow \mathcal{K}(\mathcal{H}_s) \otimes \mathfrak{S}A$$
(where $K\crt(F_s)$ is a free {\it CRT}-module for each $s$) so that the images of 
$(\mu_s)_*$ generate $K\crt(\mathfrak{S}A)$.  There will be 25 maps in all.

Recall that we have the decomposition
$$K\crt(\mathfrak{S}A) \cong K\crt(S^8A) \oplus K\crt((S^+)^7 \R) \; $$ 
arising from the split exact sequence
$$0 \rightarrow S^8A \rightarrow \mathfrak{S}A \rightarrow (S^+)^7 \R 
		\rightarrow 0 \; . $$

The first map in our collection is taken to be the section
$$s \colon (S^+)^7 \R \rightarrow \mathfrak{S}A$$
whose image in united $K$-theory is the the second summand in the decomposition above.

Now for each integer $i$ satisfying $0 \leq i \leq 7$, we let 
$$\mu_i\po \colon F_i\po \rightarrow \mathcal{K}(\mathcal{H}_i\po) \otimes S^iA$$
be given as in Lemma~\ref{realsurjection} so that the image of $(\mu_i\po)_*$ contains 
$KO_0(S^iA) \cong KO_i(A)$.  Then we suspend to get a map with the correct codomain: 
	$$S^{8-i} \mu_i\po \colon S^{8-i}(F_i\po) \rightarrow 
	S^8(\mathcal{K}(\mathcal{H}_i\po) \otimes A) 
	\hookrightarrow \mathcal{K}(\mathcal{H}_i\po) \otimes \mathfrak{S}A \; .$$  
The image of the maps $\{S^{8-i}\mu_i\po \}_{s = 0}^{7}$ induced on real $K$-theory is 
all of $KO_*(S^8A) \subset KO_*(\mathfrak{S}A)$.  It also follows, therefore, that
the image of the maps induced on complex and self-conjugate $K$-theory contains the
image of $c$ (and therefore of $\beta\su^{-1} c$) in $KU_*(S^8A)$ and the image 
of $\ve$ (and therefore of $\beta\st^{-1} \ve$) in 
$KT_*(S^8A)$.   

Similarly, we use Lemma~\ref{complexsurjection} to build maps
$$S^{8-i} \mu_i\pu \colon S^{8-i}(F_i\pu) 
\rightarrow S^8(\mathcal{K}(\mathcal{H}_i\pu) \otimes A) 
\hookrightarrow
\mathcal{K}(\mathcal{H}_i\pu) \otimes \mathfrak{S}A \; $$ 
so that the image of the
corresponding maps induced on $K$-theory contains $KU_*(S^8A)$ modulo 
$\beta\su^{-1} c$.

Finally, for each integer $i$, Lemma~\ref{selfconjugatesurjection} guarantees the 
existence of a map
$$\mu_i\pt \colon F_i\pt \rightarrow \mathcal{K}(\mathcal{H}_i\pt) \otimes S^+(S^iA)$$
such that the image of the map induced on self-conjugate $K$-theory, in degree $-1$, 
is surjective modulo $\beta\st^{-1} \ve$.  Then we suspend to get maps
$$S^{8-i} \mu_i\pt \colon S^{8-i} F_i\pt \rightarrow \mathcal{K}(\mathcal{H}_i\pt) 
\otimes S^{8-i}(S^+S^iA)
\hookrightarrow \mathcal{K}(\mathcal{H}_i\pt) \otimes \mathfrak{S}A \; .$$
The image of the corresponding maps induced on self-conjugate $K$-theory lies in 
$KT_*(\mathfrak{S}A)$.  But modulo the image of $s_*$ and the image of
$\beta\st^{-1} \varepsilon$ we know that the image of $S^{8-i} \mu_i\pt$ contains
all of $KT_*(S^8 A)$.  

Therefore, the images of the maps on united $K$-theory induced by the homomorphisms $s$, 
$\mu_i\po$, $\mu_i\pu$, and $\mu_i\pt$ contain all of $K\crt(\mathfrak{S}A)$ and
hence these maps can be used to patch together a map $\mu$ that is surjective on
united $K$-theory.  \end{proof}

\begin{proof}[Completion of Proof of Proposition~\ref{crtobject}]
Given a real unital C*-algebra $A$, let $F$ and $\mu$ be given as in the 
Proposition~\ref{surjection}.  Since $K\crt(F)$ is a free {\it CRT}-module, it
follows in particular that the {\it CRT}-relations are satisfied.  Therefore, the
{\it CRT}-relations hold for the image of $\mu_*$, which contains $K\crt(S^8A) \cong
K\crt(A)$.  So $K\crt(A)$ is a {\it CRT}-module for all real unital C*-algebras and
thus, using a unitization argument, for all real C*-algebras.    \end{proof}

Now that Proposition~\ref{crtobject} is finally proven, it is also established that 
$K\crt(A)$ has projective dimension at most $1$ by Theorem~\ref{acyclic=pd1}.  This
fact and Proposition~\ref{surjection} will both be used in the proof of
Theorem~\ref{resolution}.

\begin{proof}[Proof of Theorem \ref{resolution}] Given $A$, we take $F$, $\mu$, and 
$\mathcal{H}$ as in Proposition~\ref{surjection} so that 
$$K\crt(F) \mapis{\mu_*} K\crt(\mathcal{K}(\mathcal{H}) \otimes \mathfrak{S}A)$$
is a surjection.

Now let $F_0$ be the mapping cylinder of $\mu$ and let $F_1$ be the mapping cone of 
$\mu$:
\begin{align*}
F_0 &= \{(f, g) \in F \oplus C([0,1], \mathcal{K}(\mathcal{H}) \otimes \mathfrak{S}A) 
	\mid g(0) = \mu(f) \}    \\
F_1 &= \{(f, g) \in F \oplus C([0,1], \mathcal{K}(\mathcal{H}) \otimes \mathfrak{S}A) 
	\mid g(0) = \mu(f), g(1) = 0 \} \; .
\end{align*}

The homomorphism $\pi \colon F_0 \rightarrow F$ defined by $\pi(f, g) = f$ is a homotopy 
equivalence.  Furthermore, the homomorphism $\mu_0 \colon F_0 \rightarrow
\mathcal{K}(\mathcal{H}) \otimes \mathfrak{S}A$ defined by $\mu_0(f, g) =  g(1)$ is
equivalent to $\mu$ in the sense that the diagram $$ \xymatrix{ F_0 \ar[r]^-{\mu_0}
\ar[d]_{\pi} & K(H) \otimes \mathfrak{S}A   \\
F \ar[ur]_{\mu}
}  $$
commutes up to homotopy.  Indeed, there is a homotopy from $\mu_0$ to the 
homomorphism $(f,g) \mapsto g(0)$ which would make the diagram commute on 
the nose.

The homomorphism $\mu_0$ is surjective and its kernel is $F_1$.  Let $\mu_1 \colon F_1 
\rightarrow F_0$ be the inclusion.  Then we have the short exact sequence of
C*-algebras $$0 \rightarrow 
F_1 \mapis {\mu_1} 
F_0 \mapis {\mu_0} 
\mathcal{K}(\mathcal{H}) \otimes \mathfrak{S}A \rightarrow 0 \; . $$

Since $(\mu_0)_*$ is surjective on united $K$-theory, the long exact sequence 
on united $K$-theory breaks 
down into a resolution of $K\crt(\mathfrak{S}A)$.  Furthermore, $K\crt(F_0) \cong
K\crt(F)$ is a free {\it CRT}-module.  Since $K\crt(\mathfrak{S}A)$ has projective
dimension at most 1 in {\it CRT}, the {\it CRT}-module $K\crt(F_1)$ must also be
projective (see Theorem~9.5 of \cite{Rotman}), hence free (by
Theorem~\ref{free=pro}).    
\end{proof}


\section{Tensor Products} \label{tp}

Before we can state and prove the K\"unneth Theorem for united $K$-theory, we must 
discuss the tensor product in the category {\it CRT}.  I am grateful to A.K.
Bousfield for sharing with me some valuable unpublished notes on this topic (\cite
{Bou3}).  The following sections are largely drawn from those notes.

In Section~\ref{tp1}, we define the tensor product in the category {\it CRT}.  Given 
two {\it CRT}-modules $M$ and $N$, their tensor product $M \otimes\scrt N$ will
also be a {\it CRT}-module.  In Section~\ref{acyclictp} we prove some results about
the behavior of tensor products when one of the factors is acyclic or when one of
the factors is free.  These results are important both for the proof of our K\"unneth
formula and for purposes of computation (see Section~\ref{rca}).


\subsection{Tensor Products in \it{CRT}} \label{tp1}

Let $M = \{M\po, M\pu, M\pt\}$ be a {\it CRT}-module.  In particular $M\po$ is a 
left module over $KO_*(\R)$.  We give $M\po$ the structure of a right $KO_*(\R)$
module by the rule $m \cdot x = (-1)^{|m|\cdot|x|} x \cdot m$ for all $m \in M\po$
and $x \in KO_*(\R)$ where $|m|$ and $|x|$ denote the graded degree of those
elements.  Similarly, $M\pu$ and $M\pt$ are given right $KU_*(\R)$- and
$KT_*(\R)$-module structures, respectively.

Now if $N = \{N\po, N\pu, N\pt \}$ is also a {\it CRT}-module we have the tensor products
 $M\po \otimes_{KO_*(\R)} N\po$, $M\pu \otimes_{KU_*(\R)} N\pu$, and $M\pt
\otimes_{KT_*(\R)} N\pt$ using the right module structures on $M$ described above
and the natural left module structures on $N$.  These tensor products are again
$KO_*(\R)$-, $KU_*(\R)$-, and $KT_*(\R)$-modules, respectively.

\begin{defn} \label{CRT-pairing} Let $M$, $N$, and $P$ be {\it CRT}-modules.  A 
{\it CRT}-pairing $ {\alpha} \colon (M, N) \longrightarrow P$ consists of three
homomorphisms: \begin{enumerate}
\item[{\rm(1)}] a $KO_*(\R)$-module homomorphism 
	${\alpha}\po \colon 
		M\po \otimes_{KO_*(\R)} N\po \longrightarrow P\po \; $ 
\item[{\rm(2)}] a $KU_*(\R)$-module homomorphism  
	${\alpha}\pu \colon 
		M\pu \otimes_{KU_*(\R)} N\pu \longrightarrow P\pu \;$  
\item[{\rm(3)}] a $KT_*(\R)$-module homomorphism
	${\alpha}\pt \colon
		M\pt \otimes_{KT_*(\R)} N\pt \longrightarrow P\pt \; .$
\end{enumerate}
We use the notation ${\alpha}\po(m\so \otimes n\so)~=~m\so \cdot\so n\so$, 
${\alpha}\pu(m\su \otimes n\su)~=~m\su \cdot\su n\su$, and ${\alpha}\pt(m\st
\otimes n\st)~=~m\st \cdot\st n\st$ to express these products.  Furthermore, in
order to be a {\it CRT}-pairing, the homomorphisms $\alpha\po$, $\alpha\pu$, and $\alpha\pt$ must
satisfy the following properties for all $m\so \in M\po$, $n\so \in N\po$, $m\su \in
M\pu$, $n\su \in N\pu$, $m\st \in M\pt$, and $n\st \in N\pt$: 
\begin{enumerate}
\item[{\rm(1)}] 
	$c(m\so \cdot\so n\so) = c(m\so) \cdot\so c(n\so)$ 
\item[{\rm(2)}]
	$r(c(m\so) \cdot\su n\su) = m\so \cdot\so r(n\su)$ 
\item[{\rm(3)}] 
	$r(m\su \cdot\su c(n\so)) = r(m\su) \cdot\so n\so$ 
\item[{\rm(4)}] 
	$\varepsilon(m\so \cdot\so n\so) 
		= \varepsilon(m\so) \cdot\st \varepsilon(n\so)$ 
\item[{\rm(5)}]
	$\zeta(m\st \cdot\st n\st) = \zeta(m\st) \cdot\su \zeta(n\st)$ 
\item[{\rm(6)}]
	$\psi\su(m\su \cdot\su n\su) = \psi\su(m\su) \cdot\su \psi\su(n\su)$ 
\item[{\rm(7)}]
	$\psi\st(m\st \cdot\st n\st) = \psi\st(m\st) \cdot\st \psi\st(n\st)$ 
\item[{\rm(8)}]
	$\gamma(m\su \cdot\su \zeta(n\st)) = \gamma(m\su) \cdot\st n\st$ 
\item[{\rm(9)}]
	$\gamma(\zeta(m\st) \cdot\su n\su) 
		= (-1)^{|m\st|} m\st \cdot\st \gamma(n\su)$ 
\item[{\rm(10)}] 
	$\tau(m\st \cdot\st \varepsilon(n\so)) = \tau(m\st) \cdot\so n\so$ 
\item[{\rm(11)}]
	$\tau(\varepsilon(m\so) \cdot\st n\st) 
		= (-1)^{|m\so|} m\so \cdot\so \tau(n\st)$  
\item[{\rm(12)}] $\varepsilon \tau(m\st \cdot\st n\st) 
		= \varepsilon \tau(m\st) \cdot\st n\st 		
 		+ (-1)^{|m\st|} m\st \cdot\st \varepsilon \tau (n\st)   
		+ \eta\st(m\st \cdot\st n\st) \; .$ 
\end{enumerate} \end{defn}

In the above definition, the first three equations are redundant in the sense that they 
can be derived from the others using the relations $r = \tau \gamma$ and $c = \zeta
\varepsilon$.  We include them explicitly on our list simply out of convenience,
remembering that it is enough to check properties 4 through 12 to have a {\it
CRT}-pairing.

The tensor product in the category of {\it CRT}-modules is defined as the solution of a 
universal property described in the statement of the following proposition.

\begin{prop} \label{crttensorproduct}
Let $M$ and $N$ be {\it CRT}-modules.  Then there exists a unique {\it CRT}-module, 
denoted $M \otimes\scrt N$, and {\it CRT}-pairing 
$\iota \colon (M, N) \longrightarrow M \otimes\scrt N$ such that 
for every {\it CRT}-pairing ${\alpha} \colon (M, N)
\longrightarrow P$ there exists a unique {\it CRT}-module homomorphism $\alpha
\colon M \otimes\scrt N \longrightarrow P$ making the diagram below commute.
$$\xymatrix{ (M, N) \ar[r]^-{\iota} \ar[dr]_{{\alpha}}  	&  M \otimes\scrt N
\ar[d]^{\alpha}   \\ 	& P   }$$
\end{prop}

\begin{proof} We construct $M \otimes\scrt N$ as follows.  Let $R$ be the free 
{\it CRT}-module generated by the elements of $M\po \times N\po$, $M\pu \times
N\pu$, and $M\pt \times N\pt$.  These generators are denoted as pure tensors $m\so
\otimes n\so$, $m\su \otimes n\su$, and $m\st \otimes n\st$.  Let $Q$ be the
smallest {\it CRT}-submodule of $R$ such that the composition pairing 
$(M,N) \rightarrow R \rightarrow R/Q$ 
is a {\it CRT}-pairing.  In other words, $Q$
is the {\it CRT}-submodule generated by elements of $R$ corresponding to the {\it
CRT}-pairing relations in Definition~\ref{CRT-pairing} \end{proof}

\begin{defn} 
Given a {\it CRT}-module $M$, we define the suspension of $M$ by shifting the degrees 
such that
$$(\Sigma M)_k = M_{k+1}$$
and the desuspension by
$$(\Sigma^{-1} M)_k = M_{k-1}$$
\end{defn}
These are defined such that $K\crt(S^n A) = \Sigma^n K\crt(A)$ 
for any integer $n$.

The tensor product in the category of {\it CRT}-modules enjoys the usual 
properties associated with tensor products as stated in the following 
proposition.  The proof uses standard techniques of homological algebra, 
the details of which can be found in Section IV.1 of 
\cite{Boer}.

\begin{prop} \label{tensorproperties}
The tensor product functor is a symmetric bifunctor which in each argument is covariant,
 continuous, right exact, and commutes with suspensions and direct sums.
\end{prop}

Since the tensor product is right exact and since there are enough projective modules 
in {\it CRT} to form projective resolutions, we also have  
derived functors $\tor^i\scrt(M,N)$ for all {\it CRT}-modules $M$ and $N$ (see
\cite{Rotman}, page 121).  

The next three propositions characterize tensor products with monogenic free objects 
of {\it CRT}.  They are due to Bousfield (\cite{Bou3}) and proofs can be found in
\cite{Boer}.

\begin{prop} \label{tensorR} Let $N$ be an arbitrary {\it CRT}-module.  Then
\begin{align*}
F(b, 0, \R) \otimes\scrt N  
& \cong \{ \,  
\{b \otimes n\so \mid n\so \in N\po \}, \\  
& \hspace{2cm} \{cb \otimes n\su \mid n\su \in N\pu \},  \\
& \hspace{3cm} \{\varepsilon b \otimes n\st \mid n\st \in N\pt \} \, \}    \\
& \cong N  \; . 
\end{align*}
\end{prop} 

\begin{prop} \label{tensorC} Let $N$ be an arbitrary {\it CRT}-module.  Then
\begin{align*}
F(b, 0, \C) \otimes\scrt N 
&\cong \{ \,
\{r(b \otimes n\su) \mid n\su \in N\pu\},  \\
& \hspace{2cm} \{b \otimes n\su \mid n\su \in N\pu\} \oplus 
\{\psi\su b \otimes n\su \mid n\su \in N\pu\},  \\
& \hspace{3cm} \{\gamma (b \otimes n\su) \mid n\su \in N\pu \}  \oplus 
\{\varepsilon r(b \otimes n\su) \mid n\su \in N\pu \}
\, \}    \\
&\cong \{N\pu, N\pu \oplus N\pu, \Sigma N\pu \oplus N\pu \} \; . 
\end{align*}
\end{prop} 

\begin{prop} \label{tensorT} Let $N$ be an arbitrary {\it CRT}-module.  Then
\begin{align*}
F(b, 0, T) \otimes\scrt N 
&\cong \{ \,
\{ \tau(b \otimes n\st) \mid n\st \in N\pt \},  \\
& \hspace{2cm} \{ \zeta b \otimes n\su \mid n\su \in N\pu \} \oplus
	\{ c \tau b \otimes n\su \mid n\su \in N\pu \},  \\
& \hspace{3cm} \{ b \otimes n\st \mid n\st \in N\pt \} \oplus 
	\{ \varepsilon \tau b \otimes n\st \mid n\st \in N\pt \}
\, \} \\
&\cong \{\Sigma^{-1} N\pt, N\pu \oplus \Sigma^{-1} N\pu, N\pt \oplus \Sigma^{-1} N\pt \}   \; . \end{align*}
\end{prop} 

We say that a {\it CRT}-module $N$ is flat if tensoring on the right by $N$ 
preserves exact sequences of {\it CRT}-modules.

\begin{cor} Free objects in {\it CRT} are flat.
\end{cor}

\begin{proof}
It follows immediately from Propositions~\ref{tensorR}, \ref{tensorC}, and \ref{tensorT} 
that $F(b, 0, \R)$, $F(b, 0, \C)$, and $F(b, 0, T)$ are flat.  Since all free {\it
CRT}-modules are directs sums of suspensions and desuspensions of these three
monogenic free objects, the corollary follows using Lemma~\ref{tensorproperties}.
\end{proof}

\begin{cor} \label{complextensor}
Suppose $M$ and $N$ are {\it CRT}-modules such that $M$ is acyclic.  Then the natural 
homomorphism 
$$M\pu \otimes_{KU_*(\R)} N\pu \rightarrow (M \otimes\scrt N)\pu  $$
is an isomorphism.  
\end{cor}

\begin{proof}
Again, this follows from Propositions~\ref{tensorR}, \ref{tensorC}, and \ref{tensorT} 
in case $M$ is free.  Otherwise, $M$ has a free resolution 
$$0 \rightarrow F_1 \rightarrow F_0 \rightarrow M \rightarrow 0$$ 
of length one by Theorem~\ref{acyclic=pd1}.  Then consider the commutative diagram
$$\xymatrix{
0    					\ar[d]
& 0     				\ar[d]  \\
\tor_{KU_*(\R)}(M\pu, N\pu)		\ar[d]  \ar[r]
& (\tor\scrt(M, N))\pu			\ar[d]  \\
(F_1)\pu \otimes_{KU_*(\R)} N\pu   	\ar[d]  \ar[r]
& (F_1 \otimes\scrt N)\pu		\ar[d]  \\
(F_0)\pu \otimes_{KU_*(\R)} N\pu   	\ar[d]  \ar[r]
& (F_0 \otimes\scrt N)\pu	   	\ar[d]  \\
M\pu \otimes_{KU_*(\R)} N\pu  	 	\ar[d]  \ar[r]
& (M \otimes\scrt N)\pu		   	\ar[d]  \\
0				   
& 0
} $$
in which the vertical sequences are exact and the middle two horizontal homomorphisms 
are isomorphisms.  It follows from the five lemma that the other horizontal
homomorphisms are isomorphisms. \end{proof}

On the other hand, in the real and self-conjugate cases, there are homomorphisms 
$$M\po \otimes_{KO_*(\R)} N\po \rightarrow (M \otimes\scrt N)\po  $$
and 
$$M\pt \otimes_{KT_*(\R)} N\pt \rightarrow (M \otimes\scrt N)\pt $$
which are not isomorphisms in general.


\subsection{Acyclic Objects and Tensor Products} \label{acyclictp}

In this section we will prove the following proposition.

\begin{prop} \label{acyclicflat} Let $M$ be free and $N$ be acyclic in {\it CRT}.  
Then $M \otimes\scrt N$ is acyclic.
\end{prop}

On the other hand, we will see in Table~\ref{tensorA} that the tensor product 
of two arbitrary 
acyclic {\it CRT}-modules is not necessarily acyclic.  

To prove Proposition~\ref{acyclicflat} it suffices to assume that $M$ is a free 
monogenic {\it CRT}-module.  In fact, since tensor products commute with the
operation of suspension, it suffices to assume that $M$ is either $F(b, 0, \R)$,
$F(b, 0, \C)$, or $F(b, 0, T)$ (since suspensions and products preserve acyclicity).  
The case $M = F(b, 0, \R)$ is immediate since $F(b, 0, \R) \otimes\scrt N \cong N$.  
The other two cases are addressed in the following lemmas.

\begin{lemma} $F(b, 0, \C) \otimes\scrt N$ is acyclic for any {\it CRT}-module $N$.
\end{lemma}  

\begin{proof}  Recall that 
\begin{align*}
F(b, 0, \C) \otimes\scrt N 
&\cong \{ \,
\{r(b \otimes n\su) \mid n\su \in N\pu\},  \\
& \hspace{2cm} \{b \otimes n_1 + \psi\su b \otimes n_2 \mid n_i \in N\pu\},  \\
& \hspace{3cm} \{\gamma (b \otimes n_1) 
		+ \varepsilon r(b \otimes n_2) \mid n_i \in N\pu \}  \}
\, \}    \end{align*}

To show $F(b, 0, \C) \otimes\scrt N$ is acyclic, there are three sequences that must 
be shown to be exact.  This can be done directly, using the representation of $F(b,
0, \C) \otimes\scrt N$ above.  We illustrate by showing that the image of $\gamma$
is the kernel of $\zeta$ (which shows exactness at one point of the exact sequence
relating the complex and self-conjugate parts).

Let $b \otimes n_1 + \psi\su(b) \otimes n_2$ be an arbitrary element of 
$(F(b,0,\C) \otimes\scrt N)\pu$.  Then
$$\gamma(b \otimes n_1 + \psi\su(b) \otimes n_2) 
= \gamma(b \otimes n_1) + \gamma \psi\su(b \otimes \psi\su(n_2))   
	= \gamma(b \otimes (n_1 + \psi\su(n_2)))  \, .$$
On the other hand, if $\gamma(b \otimes n_1) + \varepsilon r(b \otimes n_2)$ is an 
arbitrary element of $(F(b,0,\C) \otimes\scrt N)\pt$ then
$$\zeta(\gamma(b \otimes n_1) + \varepsilon r(b \otimes n_2))
	= cr(b \otimes n_2)    
	= (b \otimes n_2) + (\psi\su(b) \otimes \psi\su(n_2)) \; .$$
Therefore, 
\begin{align*}
\ker \zeta &= \{ \gamma(b \otimes n_1) + \varepsilon r(b \otimes n_2) 
			\mid n_1 \in N\pu \text{~and~} n_2 = 0 \}	  \\
	&= \{\gamma(b \otimes n) \mid n \in N\pu \} \\
	&= \im \gamma   \, .
\end{align*}
\end{proof}		

Before proving that $F(b,0,T) \otimes N$ is acyclic, we need the following lemma 
concerning acyclic objects (compare \cite{Bou}, Paragraph~2.3).  
 
\begin{lemma} If $M$ is acyclic, then the sequence 
\begin{equation} \label{sequencettu} 
\cdots \longrightarrow  
M\pt_{n} \mapis{\eta\st} 
M\pt_{n+2} \xrightarrow{\smv{\zeta}{c \tau}} 
M\pu_{n+2} \oplus M\pu_{n+3}
		\xrightarrow{\smh{-\varepsilon r \beta\su}{\gamma \beta\su}}
M\pt_{n+2} \longrightarrow \cdots 
\end{equation} is exact.
\end{lemma}

\begin{proof} Consider the following diagram.
$$ \xymatrix {
{\cdots} \ar[r] &
M\pu \ar[rr]^{\gamma} \ar[d]^{\beta\su^{-1}} & &
M\pt \ar[r]^\zeta \ar[d]^{\zeta \varepsilon \tau}  &
M\pu \ar[r]^{1 - \psi\su} \ar[d]^{\varepsilon r \beta\su^{-1}}	&
M\pu \ar[rr]^\gamma \ar[d]^{\beta\su^{-1}} & &
M\pt \ar[r] \ar[d]^{\zeta \varepsilon \tau}  &
{\cdots} \\
{\cdots} \ar[r] &
M\pu \ar[rr]^{\beta\su(1 - \psi\su)}	& &
M\pu \ar[r]^{\gamma \beta\su^{-1}}	&
M\pt \ar[r]^{\zeta}	& 
M\pu \ar[rr]^{\beta\su(1 - \psi\su)} & &
M\pu \ar[r] &
{\cdots} } $$
In this diagram, the rows are exact since $M$ is acyclic.  We check that the squares 
commute:
\begin{align*}
\beta\su (1 - \psi\su) \beta\su^{-1} 
		&= 1 + \psi\su 
		= cr 
		= \zeta \varepsilon \tau \gamma  
					\displaybreak[0]  \\
\gamma \beta\su^{-1} \zeta \varepsilon \tau 
		&= \beta\st^{-1} \gamma \beta\su \zeta \varepsilon \tau
		= \beta\st^{-1} \eta\st \varepsilon \tau 
		= \beta\st^{-1} \varepsilon \tau \eta\st \\
		&~~~~= \varepsilon \tau \beta\st^{-1} \eta\st
			+ \eta\st \beta\st^{-1} \eta\st   
		= \ve \tau \beta\st^{-1} \gamma \beta\su \zeta  
		= \varepsilon \tau \gamma \beta\su^{-1} \zeta 
		= \varepsilon r \beta\su^{-1} \zeta   
					\displaybreak[0]  \\
\zeta \varepsilon r \beta\su^{-1} 
		&= cr \beta\su^{-1} 
		= (1 + \psi\su) \beta\su^{-1} 
		= \beta\su^{-1} (1 - \psi\su) \; .
\end{align*}
In this computation the relation $\eta\st \varepsilon \tau = \varepsilon \tau \eta\st$ 
follows from the fact that $\varepsilon$ and $\tau$ are $KO_*(\R)$-module
homomorphisms and the relation $\varepsilon(\eta\so) = \eta\st$.  

Now the exact Sequence~\ref{sequencettu} is easily derived from this commutative 
ladder above by a diagram chase, recalling that $\eta\st = \gamma \beta\su \zeta$.
\end{proof}

\begin{lemma} If $N$ is acyclic then $F(b, 0, T) \otimes\scrt N$ is acyclic.
\end{lemma}  

\begin{proof}
Recall that 
\begin{align*}
F(b, 0, T) \otimes\scrt N 
&\cong \{ \,
\{\tau(b \otimes n) \mid n \in N\pt \},  \\
& \hspace{2cm} \{\zeta b \otimes n_1 + c \tau b \otimes n_2 \mid n_i \in N\pu\},  \\
& \hspace{3cm} \{b \otimes n_1 
		+ \varepsilon \tau b \otimes n_2 \mid n_i \in N\pt \}  
\, \} \; . \end{align*}

Again, it can easily be shown by hand that the three sequences that need to be exact 
are exact.  For this, the two exact sequences
\begin{equation} \label{Nuut}
\cdots \longrightarrow  
N\pu_{n+1} \mapis{\gamma} 
N\pt_{n} \mapis{\zeta} 
N\pu_{n} \xrightarrow{1 - \psi\su}
N\pu_{n} \longrightarrow \cdots 
\end{equation}
and
\begin{equation} \label{Nttu}
\cdots \longrightarrow  
N\pt_{n} \mapis{\eta\st} 
N\pt_{n+2} \xrightarrow{\smv{\zeta}{\zeta \varepsilon \tau}} 
N\pu_{n+2} \oplus N\pu_{n+3}
		\xrightarrow{\smh{-\varepsilon r \beta\su}{\gamma \beta\su}}
N\pt_{n+2} \longrightarrow \cdots \; 
\end{equation}
must be used.  

For example, to show that the image of $\gamma$ is the kernel of $\zeta$ in 
$F(b,0,T) \otimes\scrt N$, we compute 
$$\gamma(\zeta b \otimes n_1 + c \tau b \otimes n_2) 
	= b \otimes \gamma(n_1) + \varepsilon \tau b \otimes \gamma(n_2)$$
and
$$\zeta(b \otimes n_1 + \varepsilon \tau b \otimes n_2) 
	= \zeta b \otimes \zeta(n_1) + c \tau b \otimes \zeta(n_2) \; .$$
Therefore, 
\begin{align*}
\im \gamma 
&= \{b \otimes n_1 + \varepsilon \tau b \otimes n_2 \mid n_1, n_2 \in \im \gamma \}   \\
&= \{b \otimes n_1 + \varepsilon \tau b \otimes n_2 \mid n_1, n_2 \in \ker \zeta \}    \\
&= \ker \zeta \; .  
\end{align*}

We will also show that the image of $\eta\so$ is the kernel of $c$ to illustrate how 
Sequence~\ref{Nttu} is used:
$$\eta\so(\tau(b \otimes n)) = \tau(b \otimes \eta\st n)$$
and
\begin{align*}
c(\tau(b \otimes n)) &= \zeta \varepsilon \tau(b \otimes n) \\
	&= \zeta(\varepsilon \tau b \otimes n + b \otimes \varepsilon \tau(n) 		
		+ \eta\st b \otimes n) \\
	&= c \tau b \otimes \zeta(n) + \zeta b \otimes c \tau(n)
\end{align*}
Hence
\begin{align*}
\im \eta\so &= \{ \tau(b \otimes n) \mid n \in \im \eta\st \} \\
	&= \{ \tau(b \otimes n) \mid n 
			\in \ker \zeta \cap \ker c \tau \} \\
	&= \ker c \, .   
\end{align*}
\end{proof}

This completes the proof of Proposition~\ref{acyclicflat}


\section{The K\"unneth Formula} \label{ks}


\subsection{The Statement of the K\"unneth Formula} \label{Kstatement}

In this section we will prove that the usual pairings in real, complex, and
self-conjugate $K$-theory form a {\it CRT}-pairing in united $K$-theory.  This
pairing forms one of the maps of our K\"unneth exact sequence, which is stated as a
theorem near the end of this section.

\begin{prop}
Let $A$ and $B$ be C*-algebras.   Then there is a {\it CRT}-module homomorphism
$$\alpha \colon K\crt(A) \otimes\scrt K\crt(B) \longrightarrow K\crt(A \otimes B) \; .$$
\end{prop}

\begin{proof}
By Lemma~\ref{mult2}, the homomorphisms
\begin{align*} 
{\alpha\so} \colon KO_*(A) \otimes_{KO_*(\R)} KO_*(B) 
	&\longrightarrow KO_*(A \otimes B)        \\
{\alpha\su} \colon KU_*(A) \otimes_{KU_*(\R)} KU_*(B) 
	&\longrightarrow KU_*(A \otimes B)        \\
{\alpha\st} \colon KT_*(A) \otimes_{KT_*(\R)} KT_*(B) 
	&\longrightarrow KT_*(A \otimes B)        
\end{align*}
are known to satisfy all of the {\it CRT}-pairing relations of 
Definition~\ref{CRT-pairing} except the last.  We must show that \begin{equation} \label{theformula}
\varepsilon \tau(z_1 \cdot\st z_2) = \varepsilon \tau(z_1) \cdot\st z_2 			
	+ (-1)^{|z_1|} z_1 \cdot\st \varepsilon \tau (z_2)  
	+ \eta\st(z_1 \cdot\st z_2)
\end{equation}
for any $z_1 \in KT_*(A)$ and $z_2 \in KT_*(B)$.   

First, we suppose $z_1 = \varepsilon(x_1)$ for some $x_1 \in KO_*(A)$.  In this 
case,
\begin{align*} \varepsilon \tau(z_1 \cdot\st z_2) 
&= \varepsilon \tau(\varepsilon(x_1) \cdot\st z_2) \\
&= (-1)^{|x_1|} \varepsilon(x_1 \cdot\so \tau(z_2)) \\
&= (-1)^{|x_1|} \varepsilon(x_1) \cdot\st \varepsilon \tau(z_2) \\
&= (-1)^{|z_1|} z_1 \cdot\st \varepsilon \tau(z_2)   \; .
\end{align*}
We also have 
\begin{align*}
\varepsilon \tau (z_1) \cdot\st z_2 
&= \varepsilon \tau \varepsilon(x_1) \cdot\st z_2  \\
&= \varepsilon(\eta\so \cdot\so x_1) \cdot\st z_2  \\
&= \eta\st(\varepsilon(x_1) \cdot\st z_2) \\
&= - \eta\st(z_1 \cdot\st z_2)
\end{align*}
(since $2 \eta\st = 0$) so Equation~\ref{theformula} holds in this case.    

Similarly, if $z_2 = \varepsilon(x_2)$ for some $x_2 \in KO_*(B)$ we have 
$\varepsilon \tau(z_1 \cdot\st z_2) = \varepsilon \tau(z_1) \cdot\st z_2$ and
$(-1)^{|z_1|} z_1 \cdot\st \varepsilon \tau(z_2) = -\eta\st(z_1 \cdot\st z_2)$ and
so Equation~\ref{theformula} holds in this case, too.

Next, we prove the claim that Equation~\ref{theformula} holds for a pair of 
elements $z_1$ and $z_2$ if and only if it holds for the pair $\beta\st z_1$ and
$z_2$.  Indeed, suppose that it holds for $z_1$ and $z_2$.  Then \begin{align*}
\varepsilon \tau(\beta\st z_1 \cdot z_2) 
&= \varepsilon \tau \beta\st(z_1 \cdot z_2)  \\
&= \beta\st \varepsilon \tau(z_1 \cdot z_2) + \eta\st \beta\st (z_1 \cdot z_2) \\
&= \beta\st(\varepsilon \tau z_1 \cdot z_2 
	+ (-1)^{|z_1|} z_1 \cdot \varepsilon \tau z_2 + \eta\st z_1 \cdot z_2) +
\beta\st \eta\st (z_1 \cdot z_2) \\ &= \varepsilon \tau (\beta\st z_1) \cdot z_2 
	+ \eta\st (\beta\st z_1 \cdot z_2) 
	+ (-1)^{|z_1|} \beta\st z_1 \cdot \varepsilon \tau (z_2)
\end{align*}
using the relations $\beta\st \varepsilon \tau = \varepsilon \tau \beta\st 
	+ \eta\st \beta\st$ and $2 \eta\st = 0$.  Conversely, suppose that
Equation~\ref{theformula} holds for $\beta\st z_1$ and $z_2$.  Then repeating the
above calculation, we find that it holds for $\beta\st^2 z_1$ and $z_2$.  But
$\beta\st^2 = \beta\so$ and all operations commute with $\beta\so$.  Since
$\beta\so$ is invertible, Equation~\ref{theformula} holds for $z_1$ and $z_2$.

Furthermore, because of the symmetry, Equation~\ref{theformula} holds for a pair $z_1$ 
and $z_2$ if and only if it holds for the pair $z_2$ and $z_1$. 

We are now ready to prove Equation~\ref{theformula} in general.  It suffices to assume 
that $z_1$ and $z_2$ are both elements in degree zero and that $A$ and $B$ are
unital.  Furthermore, we may assume that $z_1$ is represented by a projection in $T
\otimes A$ and that $z_2$ is represented by a projection in $T \otimes B$ (rather
than differences of projections in matrix algebras over $T \otimes A$ and $T \otimes
B$).  We will prove Equation~\ref{theformula} for the pair $\beta\st z_1$ and
$\beta\st z_2$.  Using Lemma~\ref{selfconjugateprojectiontarget}, let  $\mu_1 \colon
T \rightarrow M_{16}(S^+A)$ and $\mu_2 \colon T \rightarrow M_{16}(S^+B)$ be
homomorphisms such that the equations $(\mu_1)_*(\chi) = z_1 \in KT_{-1}(S^+A)$ and
$(\mu_2)_*(\chi) = z_2 \in KT_{-1}(S^+B)$, hold modulo $\beta\st^{-1} \varepsilon$. 
(Recall that $\chi \in KT_{-1}(T)$ is the free generator of $K\crt(T)$.)   Then the
equations $(\mu_1)_*(\beta\st \chi) = \beta\st z_1$ and $(\mu_2)_*(\beta\st \chi) =
\beta\st z_2$ hold modulo $\varepsilon$.

Using these maps, it suffices to consider the case $A = B = T$ with
the elements $z_1, z_2 \in KT_{3}(T)$.  Using the periodicity isomorphism 
$\beta\st$ and the claim above, it is
enough to consider $z_1, z_2 \in KT_{-1}(T)$.  In fact, we will establish
the formula for any elements $z_1$ and $z_2 \in KT_{-1}(T)$.

We first prove that $KT_{-1}(T)$ is generated by the image of 
$\varepsilon \colon KO_{-1}(T) \rightarrow KT_{-1}(T)$ and the image of $i_* \colon
KT_{-1}(\R) \rightarrow KT_{-1}(T)$ where $i$ is the inclusion $\R \hookrightarrow
T$.  From Table~\ref{crtT} we know that the image of $\varepsilon$ is the $+1$
eigenspace of $\psi\st$, forming a free summand of $KT_{-1}(T) \cong \Z \oplus \Z$. 
On the other hand, $\psi\st$ is multiplication by $-1$ on $KT_{-1}(\R) \cong \Z$ by
Table~\ref{crtR}; and therefore, the image of $i_*$ is contained in the $-1$
eigenspace of $\psi\st$.

Now consider the involution $\iota \colon T \otimes T \rightarrow T \otimes T$ which 
interchanges the two factors of the tensor product.  The map induced by $\iota$ on
$KT_{-1}(T) \cong K_{-1}(T \otimes T)$ interchanges the image of $\varepsilon$ and
the image of $i_*$.  Therefore $KT_{-1}(T)$ can be written as the direct sum of the
image of $\varepsilon$ and the image of $i_*$, each of which is isomorphic to $\Z$.

Since Equation~\ref{theformula} holds if either $z_1$ or $z_2$ is in the image of 
$\varepsilon$, it suffices to suppose that both $z_1$ and $z_2$ are in the image of
$i_*$.  Therefore, it suffices to demonstrate that Equation~\ref{theformula} holds
for $KT_{-1}(\R) \cong \Z$.

Now $KT_{-1}(\R)$ is generated by $\beta\st^{-1} \omega$.  Referring to 
Table~\ref{crtR} we have $(\beta\st^{-1} \omega)^2~=~0$, $\eta\st \beta\st = 0$,
and $\varepsilon \tau(\beta\st^{-1} \omega) = 2 \cdot 1\st$.  Thus \begin{align*} 
\varepsilon \tau(\beta\st^{-1} \omega) \cdot\st \beta\st^{-1} \omega 			
  + - \beta\st^{-1} \omega \cdot\st \varepsilon \tau (\beta\st^{-1} \omega)   
   &+ \eta\st(\beta\st^{-1} \omega \cdot\st \beta\st^{-1} \omega)    \\
&= 2 \cdot 1\st \cdot\st \beta\st^{-1} \omega 
	- \beta\st^{-1} \omega \cdot\st 2 \cdot 1\st + 0   \\
&= 0  \\
&= \varepsilon \tau (\beta\st^{-1} \omega \cdot\st \beta\st^{-1} \omega)
\end{align*} 
proving that Equation~\ref{theformula} holds for $z_1 = z_2 = \beta\st^{-1} \omega$.  
\end{proof}

Now that the homomorphism $\alpha$ is defined, we are finally prepared to understand 
the statement of the main theorem.  Recall that $\mathcal{N}$ is the bootstrap
category of complex C*-algebras described in the introduction.

\begin{thm}[K\"unneth Theorem] \label{KTheorem} Let $A$ and $B$ be real C*-algebras 
such that $\C \otimes A$ is in $\mathcal{N}$.  Then there is a natural short exact
sequence \begin{multline} \label{KF}
0 \longrightarrow 
K\crt(A) \otimes\scrt K\crt(B) \xrightarrow{~~\alpha~~}
K\crt(A \otimes B)   \\
\xrightarrow{~~\beta~~}
\tor\scrt(K\crt(A), K\crt(B)) \longrightarrow 0 \;
\end{multline}
which does not necessarily split.  Here $\alpha$ is the pairing described above and 
$\beta$ is a {\it CRT} homomorphism of degree $-1$.
\end{thm}

We will prove this theorem in the next section.  In Section~\ref{cuntzproducts} we 
will come across examples which show that the sequence does not split.  

The following corollary states that the subcategory of C*-algebras $A$ for which 
the K\"unneth sequence holds for all $B$ is at least as large as we would expect. 
Let $\mathcal{N}\r$ be the bootstrap category of real C*-algebras:  the smallest
subcategory of real separable C*-algebras whose complexification is nuclear which
contains real separable type I C*-algebras; which is closed under the operations of
taking inductive limits, stable isomorphisms, and crossed products by $\Z$ and $\R$;
and which satisfies the two out of three rule for short exact sequences. 

\begin{cor} 
Let $A$ and $B$ be real C*-algebras with $A \in \mathcal{N}\r$.  Then the K\"unneth 
sequence is exact for $A$ and $B$.
\end{cor} 

\begin{proof}
Let $A\c$ denote the complexification $\C \otimes A$ of $A$.  It suffices to show 
that the collection of real C*-algebras $A$ such that $A\c \in \mathcal{N}$
contains real separable type I C*-algebras; is closed under the operations of taking
inductive limits, stable isomorphisms, and crossed products by $\Z$ and $\R$; and
satisfies the two out of three rule for short exact sequences.

Suppose that $A$ is a real, separable, type I C*-algebra.  Then $A\c$ is certainly 
separable.  We show that $A\c$ is type I.  Let $\pi$ be a non-zero, irreducible
representation of $A\c$ on a complex Hilbert space $\mathcal{H}$.  Then the
restriction $\pi_A$ of $\pi$ to $A$ is a non-zero representation of $A$ on the
Hilbert space $\mathcal{H}$ (thought of as a real Hilbert space).  If $\pi_A$ is
irreducible, then $\pi_A(A)$ would contain a compact operator (since $A$ is type I)
and then so must $\pi(A\c)$.

On the other hand, if $\pi_A$ is not irreducible, then there is a
non-zero  $\xi \in H$ such that $E_1 = \pi(A) \xi \subsetneq  H$.  Let $E_2 = \pi(i
A) \xi = i E_1$.  Then $E_1 + E_2 = \pi(A\c) \xi = H$.  Now, $E_1 \cap E_2$ is a
complex $A\c$-invariant subspace of $H$.  Since $\pi$ is irreducible, either $E_1
\cap E_2 = 0$ or $E_1 \cap E_2 = H$.  However, $E_1 \cap E_2 = H$ implies $E_1 = H$,
contrary to hypothesis.  Hence $E_1 \cap E_2 = 0$, so $H = E_1 \oplus E_2$.  
In this case, the representation 
$\pi$ of $A\c$ on $H$ is isomorphic to the complexification of the representation
$\pi$ of $A$ on $E_1$.  Since the former is irreducible, so is the latter.  Therefore
the latter contains a compact operator, and so does the former.

Either way, $\pi(A\c)$ contains a compact operator and therefore $A\c$ is a type I
C*-algebra. 

The operation of complexification commutes with the operations of taking inductive 
limits, of stabilizing, and of taking crossed products by by $\Z$ and by $\R$. 
Therefore the category in question is closed under these four operations.

Finally, since complexification preserves exact sequences, the collection of C*-algebras 
whose complexification is in $\mathcal{N}$ satisfies the two out of three rule.
\end{proof} 


\subsection{The Proof of the K\"unneth Formula} \label{Kproof}

We first consider the case that $K\crt(B)$ is a free {\it CRT}-module.  In this case, 
the {\it CRT}-module $\tor\scrt(K\crt(A), K\crt(B))$ vanishes and so the claim that
Sequence~\ref{KF} is exact collapses to the statement that $\alpha$ is an
isomorphism.  This is proven below in Proposition~\ref{specialcase}.  Once that is
accomplished, we will prove Theorem~\ref{KTheorem} in the general case by using a
geometric resolution of $B$.

\begin{prop} \label{specialcase} Let $A$ and $B$ be real C*-algebras such that 
$K\crt(B)$ is free and $\C \otimes A \in \mathcal{N}$.  Then $\alpha(A, B)$ is an
isomorphism. \end{prop}

\begin{proof} By Theorem~\ref{acyclic} and Proposition~\ref{acyclicflat}, 
$K\crt(A) \otimes\scrt K\crt(B)$ and $K\crt(A \otimes B)$ are both acyclic objects.
 By Proposition~\ref{CRTrigidity}, it suffices to show that  
$$\alpha\pu \colon (K\crt(A) \otimes\scrt K\crt(B))\pu \longrightarrow 
(K\crt(A \otimes B))\pu$$  
is an isomorphism. Using Lemma~\ref{complextensor}, this homomorphism can be rewritten as
$$\alpha\pu \colon KU_*(A) \otimes_{KU_*(\R)} KU_*(B) \longrightarrow 
KU_*(A \otimes B) \; .$$
Now since $K\crt(B)$ is a free {\it CRT}-module, $KU_*(B)$ is a free abelian group by 
Theorem~3.2 of \cite{Bou}.  Furthermore $\C \otimes A$ is in the complex bootstrap
category $\mathcal{N}$, so the K\"unneth formula for complex C*-algebras
(Theorem~2.14 of \cite{Schochet}) shows that $\alpha\su$ is an isomorphism.
\end{proof}

The next proposition and the lemmas which follow will complete the proof of the 
K\"unneth Theorem.

\begin{prop} \label{gencase} Suppose $A$ is a real C*-algebra such that 
$\C \otimes A$ is nuclear and $\alpha(A, B)$ is an isomorphism whenever $K\crt(B)$ is a 
free {\it CRT}-module. 
Then for any C*-algebra $B$, the K\"unneth formula holds for the pair $(A, B)$.
\end{prop}

\begin{proof}
Let $A$ be as specified in the hypothesis.  Let $B$ be an arbitrary unital C*-algebra.  
Then, using Theorem~\ref{resolution}, let 
\begin{equation} \label{geometry} 
0 \longrightarrow F_1 \mapis{\mu_1} F_0 \mapis{\mu_0} 
\mathcal{K}(\mathcal{H}) \otimes \mathfrak{S}B \longrightarrow 0 \end{equation}
be a sequence of C*-algebras such that
\begin{equation} \label{algebra}
0 \longrightarrow K\crt(F_1) \xrightarrow{(\mu_1)_*} 
K\crt(F_0) \xrightarrow{(\mu_0)_*} K\crt(\mathfrak{S}B) \longrightarrow 0 \end{equation}
is a free resolution of $K\crt(\mathfrak{S}B)$.

Now, we take Sequence~\ref{geometry} and tensor everything on the left by the 
C*-algebra A obtaining  
\begin{equation}   \label{geometrytensor}
0 \longrightarrow A \otimes F_1 \xrightarrow{1 \otimes \mu_1} 
A \otimes F_0 \xrightarrow{1 \otimes \mu_0} 
A \otimes \mathcal{K}(\mathcal{H}) \otimes \mathfrak{S}B \longrightarrow 0 \; .
\end{equation}
The condition that $\C \otimes A$ is nuclear insures that tensoring by $A$ 
preserves exact sequences.  We unsplice the induced long exact 
sequence on united $K$-theory, forming
$$0 \longrightarrow 
\coker((1 \otimes \mu_1)_*) \xrightarrow{(1 \otimes \mu_0)_*}
K\crt(A \otimes \mathfrak{S}B) \mapis{\delta} 
\ker((1 \otimes \mu_1)_*) \longrightarrow 
0 $$
where the connecting homomorphism $\delta$ is a map of degree $-1$.

It remains to make the identifications
\begin{align*}
\ker((1 \otimes \mu_1)_*) &= \tor \scrt(K\crt(A), K\crt(\mathfrak{S}B)) \\
\coker((1 \otimes \mu_1)_*) &= K\crt(A) \otimes \scrt K\crt(\mathfrak{S}B) \;
\end{align*}
which is done using the following commutative diagram.
\begin{equation} \label{kunnethdiagram}
\xymatrix{    
0 \ar[d]   \\
\tor\scrt(K\crt(A), K\crt(\mathfrak{S}B)) \ar[d]  		  &&
K\crt(A \otimes \mathfrak{S}B) \ar[d]^{\delta}      \\
K\crt(A) \otimes\scrt K\crt(F_1) 
	\ar[rr]^-{\alpha(A, F_1)} \ar[d]^{1 \otimes (\mu_1)_*}    &&
K\crt(A \otimes F_1) \ar[d]^{(1 \otimes \mu_1)_*}   \\
K\crt(A) \otimes\scrt K\crt(F_0) 
	\ar[rr]^-{\alpha(A, F_0)} \ar[d]^{1 \otimes (\mu_0)_*}    &&
K\crt(A \otimes F_0) \ar[d]^{(1 \otimes \mu_0)_*}   \\
K\crt(A) \otimes\scrt K\crt(\mathfrak{S}B) 
	\ar[rr]^-{\alpha(A, \mathfrak{S}B)} \ar[d]    &&
K\crt(A \otimes \mathfrak{S}B)   		    \\
0   	}
\end{equation}

The vertical sequence on the left is exact, derived from the resolution 
given by Sequence~\ref{algebra}.  The vertical sequence on the right is 
also exact, induced by Sequence~\ref{geometrytensor}.  
The horizontal maps $\alpha(A, F_1)$ and $\alpha(A, F_0)$ are
isomorphisms since $K\crt(F_0)$ and $K\crt(F_1)$ are free.  The diagram 
then shows that the kernel of $(1 \otimes \mu_1)_*$ is isomorphic to the 
kernel of $1 \otimes (\mu_1)_*$ which, in turn, is isomorphic to $\tor\scrt(K\crt(A),
K\crt(\mathfrak{S}B))$.  Composing $\delta$ with this isomorphism, we obtain the homomorphism 
$\beta \colon K\crt(A \otimes \mathfrak{S}B) 
\rightarrow \tor\scrt(K\crt(A), K\crt(\mathfrak{S}B))$ 
of the K\"unneth formula.

Similarly, Diagram~\ref{kunnethdiagram} shows that the cokernel of $(1 \otimes \mu_1)_*$ is
isomorphic to the cokernel of $1 \otimes (\mu_1)_*$ which, in turn, is isomorphic to
$K\crt(A) \otimes\scrt K\crt(\mathfrak{S}B)$.  Furthermore, in this last
identification, the composition $$K\crt(A) \otimes K\crt(\mathfrak{S}B)
\longrightarrow  \coker(1 \otimes \mu_1)_* \longrightarrow  K\crt(A \otimes
\mathfrak{S}B)$$ is the homomorphism $\alpha(A, \mathfrak{S}B)$.

This proves that the K\"unneth formula holds for the pair $(A, \mathfrak{S}B)$.  Recall 
that $\mathfrak{S}B$ is the eight-fold unital suspension of $B$.  By repeatedly
applying Lemmas~\ref{suspension} and \ref{unitization} below, it follows that the
K\"unneth formula holds for the pair $(A, B)$.  

A further application of Lemma~\ref{unitization} shows that the K\"unneth formula 
holds if $B$ is not unital.    
\end{proof}

\begin{lemma} \label{suspension} The K\"unneth formula is satisfied for the pair 
$(A, B)$ if and only if it is satisfied for the pair $(A, SB)$.
\end{lemma}

\begin{proof}
Consider the following diagram.
$$\xymatrix{
0 \ar[d] & & 0 \ar[d]   \\
K\crt(A) \otimes\scrt K\crt(B)   \ar[d]^{\alpha(A, B)} \ar[rr]     & &
K\crt(A) \otimes\scrt K\crt(SB)   \ar[d]^{\alpha(A, SB)}  \\
K\crt(A \otimes B)    \ar[d]^{\beta(A, B)} \ar[rr]    & &
K\crt(A \otimes SB)   \ar[d]^{\beta(A, SB)}     \\
\tor\scrt(K\crt(A), K\crt(B))   \ar[d]  \ar[rr]   & &
\tor\scrt(K\crt(A), K\crt(SB))   \ar[d]   \\
0 & &  0  } $$

The horizontal maps in this diagram are all {\it CRT}-module isomorphisms with 
degree $-1$.  Therefore, there exists a homomorphism $\beta(A, B)$ making the left
vertical sequence exact if and only if there exists a homomorphism $\beta(A, SB)$
making the right vertical sequence exact. \end{proof}

\begin{lemma} \label{unitization} Let $A$ and $B$ be real C*-algebras 
with $\C \otimes A$ nuclear.  
Suppose that the K\"unneth Formula holds for the pair $(A, B^+)$.  Then it holds
for the pair $(A, B)$. \end{lemma}

\begin{proof}
We consider the short exact sequence
$$0 \longrightarrow
B \mapis{i}
B^+ \mapis{\pi}
\R \longrightarrow 0 \; .$$
Because of the presence of a section $s \colon \R \rightarrow B^+$ such that 
$\pi \circ s = 1$ (and because $\C \otimes A$ is nuclear) the vertical sequences in the
diagram below are split exact. $$\xymatrix{
& 0 \ar[d] & 0 \ar[d] & 0 \ar[d]    \\
0 \ar[r] &  K\crt(A) \otimes\scrt K\crt(B)  
	\ar[r]^-{\alpha(A, B)} \ar[d]^{1 \otimes i_*}   &
K\crt(A \otimes B) 
	\ar[r]^-{\beta(A, B)} \ar[d]^{(1 \otimes i)_*} &
\tor(K\crt(A), K\crt(B)) 
	\ar[d]^{\tor(1, i)}  \ar[r]	&
0     \\
0 \ar[r] &
K\crt(A) \otimes\scrt K\crt(B^+)  
	\ar[r]^-{\alpha(A, B^+)} \ar[d]^{1 \otimes \pi_*}   &
K\crt(A \otimes B^+) 
	\ar[r]^-{\beta(A, B^+)} \ar[d]^{(1 \otimes \pi)_*}	&
\tor(K\crt(A), K\crt(B^+)) \ar[r] \ar[d]    &
0   \\
0 \ar[r]  &
K\crt(A) \otimes\scrt K\crt(\R) 
	\ar[r]^-{\alpha(A, \R)} \ar[d]	&
K\crt(A \otimes \R) 
	\ar[r] \ar[d]	&
0   \\
& 0 & 0
} $$

Furthermore, the second and third rows are exact by hypothesis.  In the first row, 
we define $\beta(A, B)$ to be the composition $\tor(1, i)^{-1} \circ \beta(A, B^+)
\circ (1 \otimes i)_*$.  This makes the first row a complex and by the snake lemma
(Exercise~6.16 of \cite{Rotman}) it is exact.        
\end{proof}

The proof of Theorem~\ref{KTheorem} is completed by the following Proposition.

\begin{prop} The K\"unneth Sequence is natural with respect to 
homomorphisms of real C*-algebras in either variable.
\end{prop}

\begin{proof}
It is clear that $\alpha$ is natural.  So it suffices to show that 
$\beta$ is natural with respect to either argument.  First, we show that 
the homomorphism $\beta$ is independent of the choice of geometric 
resolution.  Indeed, suppose that
$$0 \rightarrow E_1 \xrightarrow{\lambda_1} E_0 \xrightarrow{\lambda_0} 
		\mathcal{K}(\mathcal{H}) \otimes \mathfrak{S}B \rightarrow 0$$
and
$$0 \rightarrow F_1 \xrightarrow{\mu_1} F_0 \xrightarrow{\mu_0} 
		\mathcal{K}(\mathcal{H}) \otimes \mathfrak{S}B \rightarrow 0$$
are two geometric resolutions as in Theorem~\ref{resolution}.  Then we can 
form a third geometric resolution
$$0 \rightarrow D_1 \xrightarrow{\kappa_1} D_0 \xrightarrow{\kappa_0} 
		\mathcal{K}(\mathcal{H}) \otimes \mathfrak{S}B \rightarrow 0$$
by defining $D_0 = E_0 \oplus F_0$ and $D_1$ to be the kernel of
$$\kappa_0 = \smh{\mu_0}{\lambda_0} \colon D_0 \rightarrow 
\mathcal{K}(\mathcal{H}) \otimes \mathfrak{S}B \; .$$
Then we obtain a commutative diagram 
$$\xymatrix{
0 \ar[r] &
F_1 \ar[r] \ar[d]  &
F_0 \ar[r] \ar[d]  &
 {\mathcal{K}(\mathcal{H}) \otimes \mathfrak{S}B} \ar[r] \ar[d]^=  &
0   \\
0 \ar[r] &
D_1 \ar[r]   &
D_0 \ar[r] &
 {\mathcal{K}(\mathcal{H}) \otimes \mathfrak{S}B} \ar[r]   &
0   	} \; $$
where $F_0 \rightarrow D_0$ is the inclusion onto the second summand of 
$D_0 = E_0 \oplus F_0$.  This homomorphism of 
resolutions induces a homomorphism from each group of 
Diagram~\ref{kunnethdiagram} to the corresponding group of a similar diagram 
based on 
the resolution with $D_1$ and $D_0$.  The resulting commutative diagram 
shows that the homomorphism $\beta$ based on the resolution $D_1$ and $D_0$ 
is the same as the homomorphism $\beta$ based on the resolution with $F_1$ and 
$F_0$.  Similarly, a diagram based on the inclusion onto the first summand shows that this 
homomorphism $\beta$ is the same as the one based on the resolution with 
$E_1$ and $E_0$.  Therefore $\beta$ does not depend on the particular 
geometric resolution.

Now, if $\phi \colon A \rightarrow A'$ is a homomorphism of real C*-algebras, 
then we choose a fixed resolution of $\mathcal{K}(\mathcal{H}) \otimes 
\mathfrak{S}B$ to compute $\beta(A, B)$ and $\beta(A', B)$.  Because 
$\alpha$ is natural, there is a homomorphism induced by $\phi$ from each group of 
Diagram~\ref{kunnethdiagram} to the corresponding groups of the same diagram 
with $A$ replaced by $A'$.  Since $\beta$ is defined in terms of $\alpha$ 
by way of this commutative diagram, the diagram shows that $\beta$ 
is natural with respect to homomorphisms in the first argument.

Finally, suppose that we have a real C*-algebra homomorphism 
$\phi \colon B \rightarrow B'$.  Let 
$\mu \colon F \rightarrow \mathcal{K}(\mathcal{H}) \otimes \mathfrak{S}B$ and 
$\mu' \colon F' \rightarrow \mathcal{K}(\mathcal{H}) \otimes \mathfrak{S}B'$
both be given by Proposition~\ref{surjection}.  Let $\phi$ also denote the 
induced homomorphism from $\mathcal{K}(\mathcal{H}) \otimes \mathfrak{S}B$ to 
$\mathcal{K}(\mathcal{H}) \otimes \mathfrak{S}B'$.  Then 
there is a commutative diagram
$$\xymatrix{
F \ar[rr]^-{\mu} \ar[d]^{\smv{1}{0}}
&& {\mathcal{K}(\mathcal{H}) \otimes \mathfrak{S}B \ar[d]^{\phi}} \\
F \oplus F' \ar[rr]^-{\mu''}
&& {\mathcal{K}(\mathcal{H}) \otimes \mathfrak{S}B'}
			}$$
where $\mu'' = \smh{\phi \circ \mu}{\mu'}$.
The naturality of the mapping cone construction produces the following 
commutative diagram
$$\xymatrix{
0 \ar[r]
& F_1  \ar[r] \ar[d]
& F_0  \ar[r] \ar[d]
& {\mathcal{K}(\mathcal{H}) \otimes \mathfrak{S}B} \ar[r] \ar[d]^{\phi}
& 0   \\
0 \ar[r]
& F''_1 \ar[r]
& F''_0 \ar[r]
& {\mathcal{K}(\mathcal{H}) \otimes \mathfrak{S}B'} \ar[r]
& 0		}$$
where $F_1$ and $F_0$ are the mapping cone and mapping cylinder of $\mu$ 
while $F''_1$ and $F''_0$ are the mapping cone and mapping cylinder of 
$\mu''$.  This is a homomorphism from one geometric resolution 
to the other and it induces a family of homomorphisms from the groups of
Diagram~\ref{kunnethdiagram} to the groups of the same diagram with $B$ 
replaced by $B'$.  
Since $\beta$ does not depend on the particular geometric 
resolution, it follows that $\beta$ commutes with $\phi$.

\end{proof}


\section{Application:  Real Cuntz Algebras} \label{rca}


\subsection{The United $K$-Theory of the Real Cuntz Algebras} \label{cuntz}

For $k \geq 1$, the real Cuntz algebra $\mathcal{O}\pr_{k+1}$ is the universal real 
C*-algebra generated by $k+1$ isometries $S_1, S_2, \dots, S_{k+1}$ subject to the
relation $\sum_{i = 1}^{k+1} S_i S_i^* = 1$ (see \cite{Schroder}, page 4).  The
complexification of $\mathcal{O}\pr_{k+1}$ is the complex Cuntz algebra
$\mathcal{O}_{k+1}$ introduced in \cite{Cun1}.  

The $K$-theory of the complex Cuntz algebras was computed by Cuntz in 
\cite{Cun2} while the 
ordinary $K$-theory of the real Cuntz algebras was computed by Schr\"oder 
(Theorem~1.6.8 in \cite{Schroder}).  However, this last computation is not entirely
correct, as it involves a mistaken solution to a certain extension problem.  In this
section we compute the united $K$-theory of the real Cuntz algebras.  In particular,
we find the correct real $K$-groups.  In fact, it is the structure of united
$K$-theory --- namely, the relationship between real and complex $K$-theory --- that
allowed us to detect Schr\"oder's error and correctly solve the extension problem.

The form that $K\crt(\mathcal{O}\pr_{k+1})$ takes depends on the congruence class of 
$k$ modulo 4.  Tables~\ref{cuntzA}, \ref{cuntzB}, and \ref{cuntzC} below show the
united $K$-theory in the case that $k$ is odd, the case that $k$ is congruent to 2
modulo 4, and the case that $k$ is congruent to 0 modulo 4, respectively.

\begin{table}[h]   
\caption{$K\crt(\mathcal{O}\pr_{k+1})$ for $k$ odd} \label{cuntzA}
$$\begin{array}{|c|c|c|c|c|c|c|c|c|c|}  
\hline  \hline 
n & \makebox[1cm][c]{0} & \makebox[1cm][c]{1} & 
\makebox[1cm][c]{2} & \makebox[1cm][c]{3} 
& \makebox[1cm][c]{4} & \makebox[1cm][c]{5} 
& \makebox[1cm][c]{6} & \makebox[1cm][c]{7} 
& \makebox[1cm][c]{8} \\
\hline  \hline 
KO_n  &
\Z_k & 0 & 0 & 0 & 
\Z_k & 0 & 0 & 0 & 
\Z_k         \\
\hline  
KU_n   &
\Z_k & 0 & 
\Z_k & 0 & 
\Z_k & 0 & 
\Z_k & 0 & 
\Z_k         \\
\hline  
KT_n   &
\Z_k & 0 & 0 & \Z_k &  
\Z_k & 0 & 0 & \Z_k &
\Z_k         \\
\hline \hline
c_n & 1 & 0 & 0 & 0 & 2 & 0 & 0 & 0 & 1 \\
\hline
r_n & 2 & 0 & 0 & 0 & 1 & 0 & 0 & 0 & 2  \\
\hline
\varepsilon_n & 1 & 0 & 0 & 0 & 2 & 0 & 0 & 0 & 1 \\
\hline
\zeta_n & 1 & 0 & 0 & 0 & 1 & 0 & 0 & 0 & 1 \\
\hline
(\psi\su)_n & 1 & 0 & -1 & 0 & 1 & 0 & -1 & 0 & 1 \\
\hline
(\psi\st)_n & 1 & 0 & 0 & -1 & 1 & 0 & 0 & -1 & 1 \\
\hline
\gamma_n & 1 & 0 & 0 & 0 & 1 & 0 & 0 & 0 & 1 \\
\hline
\tau_n & 0 & 0 & 0 & 1 & 0 & 0 & 0 & 2 & 0  \\
\hline \hline
\end{array}$$
\end{table}

\begin{table}   
\caption{$K\crt(\mathcal{O}\pr_{k+1})$ for $k \equiv 2 \pmod 4$} \label{cuntzB}
$$\begin{array}{|c|c|c|c|c|c|c|c|c|c|}  
\hline  \hline  
n & \makebox[1cm][c]{0} & \makebox[1cm][c]{1} & 
\makebox[1cm][c]{2} & \makebox[1cm][c]{3} 
& \makebox[1cm][c]{4} & \makebox[1cm][c]{5} 
& \makebox[1cm][c]{6} & \makebox[1cm][c]{7} 
& \makebox[1cm][c]{8} \\
\hline  \hline 
KO_n  &
\Z_k & \Z_2 & \Z_4 & \Z_2 & \Z_k &
0 & 0 & 0 & \Z_k         \\
\hline  
KU_n   &
\Z_k & 0 & 
\Z_k & 0 & 
\Z_k & 0 & 
\Z_k & 0 & 
\Z_k         \\
\hline  
KT_n   &
\Z_k & \Z_2 & \Z_2 & \Z_k &  
\Z_k & \Z_2 & \Z_2 & \Z_k &
\Z_k         \\
\hline \hline
c_n & 1 & 0 & \tfrac{k}{2} & 0 & 2 & 0 & 0 & 0 & 1 \\
\hline
r_n & 2 & 0 & 2 & 0 & 1 & 0 & 0 & 0 & 2  \\
\hline
\varepsilon_n & 1 & 1 & 1 & \tfrac{k}{2} & 2 & 0 & 0 & 0 & 1 \\
\hline
\zeta_n & 1 & 0 & \tfrac{k}{2} & 0 & 1 & 0 & \tfrac{k}{2} & 0 & 1 \\
\hline
(\psi\su)_n & 1 & 0 & -1 & 0 & 1 & 0 & -1 & 0 & 1 \\
\hline
(\psi\st)_n & 1 & 1 & 1 & -1 & 1 & 1 & 1 & -1 & 1 \\
\hline
\gamma_n & 1 & 0 & 1 & 0 & 1 & 0 & 1 & 0 & 1 \\
\hline
\tau_n & 1 & 2 & 1 & 1 & 0 & 0 & 0 & 2 & 1  \\
\hline \hline
\end{array}$$
\end{table}

\begin{table}   
\caption{$K\crt(\mathcal{O}\pr_{k+1})$ for $k \equiv 0 \pmod 4$} \label{cuntzC}
$$\begin{array}{|c|c|c|c|c|c|c|c|c|c|}  
\hline  \hline  
n & \makebox[1cm][c]{0} & \makebox[1cm][c]{1} & 
\makebox[1cm][c]{2} & \makebox[1cm][c]{3} 
& \makebox[1cm][c]{4} & \makebox[1cm][c]{5} 
& \makebox[1cm][c]{6} & \makebox[1cm][c]{7} 
& \makebox[1cm][c]{8} \\
\hline  \hline 
KO_n  &
\Z_k & \Z_2 & \Z_2^2 & \Z_2 & \Z_k &
0 & 0 & 0 & \Z_k         \\
\hline  
KU_n   &
\Z_k & 0 & 
\Z_k & 0 & 
\Z_k & 0 & 
\Z_k & 0 & 
\Z_k         \\
\hline  
KT_n   &
\Z_k & \Z_2 & \Z_2 & \Z_k &  
\Z_k & \Z_2 & \Z_2 & \Z_k &
\Z_k         \\
\hline \hline
c_n & 1 & 0 & \smh{0}{\tfrac{k}{2}} & 0 & 2 & 0 & 0 & 0 & 1 \\
\hline
r_n & 2 & 0 & \smv{1}{0} & 0 & 1 & 0 & 0 & 0 & 2  \\
\hline
\varepsilon_n & 1 & 1 & \smh{0}{1} & \tfrac{k}{2} & 2 & 0 & 0 & 0 & 1 \\
\hline
\zeta_n & 1 & 0 & \tfrac{k}{2} & 0 & 1 & 0 & \tfrac{k}{2} & 0 & 1 \\
\hline
(\psi\su)_n & 1 & 0 & -1 & 0 & 1 & 0 & -1 & 0 & 1 \\
\hline
(\psi\st)_n & 1 & 1 & 1 & -1 & 1 & 1 & 1 & -1 & 1 \\
\hline
\gamma_n & 1 & 0 & 1 & 0 & 1 & 0 & 1 & 0 & 1 \\
\hline
\tau_n & 1 & \smv{1}{0} & 1 & 1 & 0 & 0 & 0 & 2 & 1  \\
\hline \hline
\end{array}$$
\end{table}

\begin{proof}[Computation of Tables~\ref{cuntzA}, \ref{cuntzB}, and \ref{cuntzC}]
The complex $K$-theory of the real Cuntz algebras is the same as the ordinary $K$-theory 
of the complex Cuntz algebras.  Hence, we refer to \cite{Cun2} for the values of
$KU_*(\mathcal{O}\pr_{k+1})$.  Furthermore, since $KU_0(\mathcal{O}\pr_{k+1}) \cong
\Z_k$ is generated by the class of the identity in $\mathcal{O}\pr_{k+1} \subseteq
\C \otimes \mathcal{O}\pr_{k+1}$, we know that $KU_0(\mathcal{O}\pr_{k+1})$ is in
the image of $c_0$.  It follows that $(\psi\su)_n$ is the trivial involution for $n
\equiv 0 \pmod 4$ and multiplication by $-1$ for $n \equiv 2 \pmod 4$.  

Let $\phi_i \colon M_{(k+1)^i} \rightarrow M_{(k+1)^{i+1}}$ be the matrix algebra 
embedding of multiplicity $k+1$ and let $A_{{(k+1)}^\infty}$ be the real UHF
algebra obtained by taking the direct limit of the system $\{M_{(k+1)^i}, \phi_i\}$.
 We find the real $K$-theory of $A_{{(k+1)}^\infty}$ in degrees 0 through 8 using
the fact that $K$-theory passes through direct limits: $$\begin{array}{cccccccccc}  
  KO_*(A_{{(k+1)}^\infty})=  & \makebox[.85cm][c]{$\Z[\tfrac{1}{k+1}]$}  &
\makebox[.85cm][c]{$\Z_2$}  & \makebox[.85cm][c]{$\Z_2$} 
& \makebox[.85cm][c]{0} 
& \makebox[.85cm][c]{$\Z[\tfrac{1}{k+1}]$} 
& \makebox[.85cm][c]{0} 
& \makebox[.85cm][c]{0} 
& \makebox[.85cm][c]{0} 
& \makebox[.85cm][c]{$\Z[\tfrac{1}{k+1}]$} \\
&&&&&&&& \multicolumn{2}{l}{\text{for $k+1$ odd}}  \\
KO_*(A_{{(k+1)}^\infty})= 
& \makebox[.85cm][c]{$\Z[\tfrac{1}{k+1}]$} 
& \makebox[.85cm][c]{0} 
& \makebox[.85cm][c]{0} 
& \makebox[.85cm][c]{0} 
& \makebox[.85cm][c]{$\Z[\tfrac{1}{k+1}]$} 
& \makebox[.85cm][c]{0} 
& \makebox[.85cm][c]{0} 
& \makebox[.85cm][c]{0} 
& \makebox[.85cm][c]{$\Z[\tfrac{1}{k+1}]$}   \\
&&&&&&&& \multicolumn{2}{l}{\text{for $k+1$ even}}  
\end{array}$$

There is an automorphism $\alpha$ of $\mathcal{K} \otimes A_{{(k+1)}^\infty}$ which 
induces multiplication by $\tfrac{1}{k+1}$ on real $K$-theory such that there is an
isomorphism  $$\mathcal{K} \otimes \mathcal{O}\pr_{k+1} \cong (\mathcal{K} \otimes
A_{{(k+1)}^\infty}) \times_\alpha \Z \; . $$ (See page 51 of \cite{Schroder},
following \cite{Cun1} in the complex case.)

We compute the $K$-theory of this crossed product, and hence of $\mathcal{O}\pr_{k+1}$, 
using the real $K$-theory version of the Pimsner-Voiculescu sequence (found as
Theorem~1.5.5 in \cite{Schroder}, following \cite{PV} in the complex case): 
$$KO_{n+1}(\mathcal{O}\pr_{k+1}) \longrightarrow KO_n(A_{{(k+1)}^\infty})
\xrightarrow{1 - \alpha_*}  KO_n(A_{{(k+1)}^\infty}) \longrightarrow 
KO_n(\mathcal{O}\pr_{k+1}) \longrightarrow
$$

In the cases where $KO_n(A_{{(k+1)}^\infty}) \cong \Z[\tfrac{1}{k+1}]$, the homomorphism 
$1 - \alpha_* = 1 - \tfrac{1}{k+1} = \tfrac{k}{k+1}$ is injective with cokernel
isomorphic to $\Z_k$.  In the cases where $k+1$ is odd and $n$ is congruent to 1 or
2 (mod 8), the homomorphism $1 - \alpha_*$ is zero.   

Using this, the Pimsner-Voiculescu sequence almost completely determines 
$KO_*(\mathcal{O}\pr_{k+1})$; the exception is $KO_2(\mathcal{O}\pr_{k+1})$ which
is an extension of $\Z_2$ by $\Z_2$ in case $k+1$ is odd.  We show that it is $\Z_4$
for $k \equiv 2 \pmod 4$ and it is $\Z_2^2$ for $k \equiv 0 \pmod 4$.

We will use the exact sequence
\begin{equation} \label{oouOk}
\cdots \longrightarrow  KO_n(\mathcal{O}\pr_{k+1}) \mapis{\eta\so} 
KO_{n+1}(\mathcal{O}\pr_{k+1}) \mapis{c} 
KU_{n+1}(\mathcal{O}\pr_{k+1}) \mapis{r \beta\su^{-1}}
KO_{n-1}(\mathcal{O}\pr_{k+1}) \longrightarrow \cdots \; . 
\end{equation}
Since $KU_3(\mathcal{O}\pr_{k+1}) \cong 0$, we deduce that $(\eta\so)_1$ must be 
injective and $(\eta\so)_2$ must be surjective.

Suppose first that $KO_2(\mathcal{O}\pr_{k+1}) \cong \Z_4$.  Then $(\eta\so)_1 = 2$ and 
$(\eta\so)_2 = 1$.  Now the image of $(\eta\so)_1$ is the kernel of $c_2$, so we
conclude that $c_2 = \tfrac{k}{2}$.  Also, since the kernel of $(\eta\so)_2$ is the
image of $r_2$, we conclude that $r_2 = 2$.  Therefore, the relation $rc = 2$ forces
$k \equiv 2 \pmod 4$.

On the other hand, suppose that $KO_2(\mathcal{O}\pr_{k+1}) \cong \Z_2^2$.  Since 
$(\eta\so)_1$ is injective we can write $(\eta\so)_1 = \smv{1}{0} \colon \Z_2
\rightarrow \Z_2^2$ (written in terms of appropriately chosen generators of
$KO_2(\mathcal{O}\pr_{k+1})$).  Now $KO_1(\mathcal{O}\pr_{k+1})$ is in the image of
$\eta\so$.  Since $\eta\so^3 = 0$ we know that $(\eta\so)_2$ vanishes on the image
of $(\eta\so)_1$.  But $(\eta\so)_2$ is non-zero, so we conclude that $(\eta\so)_2 =
\smh{0}{1}$.

Using Sequence~\ref{oouOk} again, we determine $r_2 = \smv{1}{0}$ and 
$c_2 = \smh{0}{\tfrac{k}{2}}$.  Then $rc = 0$ implies $\tfrac{k}{2} \equiv 0 \pmod
2$, which implies $k \equiv 0 \pmod 4$.

This completes the computation of $KO_*(\mathcal{O}\pr_{k+1})$.  Now the computation 
of $KT_*(\mathcal{O}\pr_{k+1})$ and the remaining operations proceeds without
difficulty using the exact sequence $$\cdots \longrightarrow 
KU_{n+1}(\mathcal{O}\pr_{k+1}) \mapis{\gamma}  KT_n(\mathcal{O}\pr_{k+1})
\mapis{\zeta}  KU_n(\mathcal{O}\pr_{k+1}) \xrightarrow{1 - \psi\su}
KU_{n}(\mathcal{O}\pr_{k+1}) \longrightarrow \cdots \; ,$$
and the known values of $\psi\su$, as well as the relations $rc = 2$, 
$c = \zeta \varepsilon$, $r = \tau \gamma$, and $\ve r \zeta = 1 + \psi\st$.
\end{proof}


\subsection{The United $K$-Theory of Tensor Products of Real Cuntz Algebras} 
\label{cuntzproducts}

In this section, we compute the united $K$-theory of 
$\mathcal{O}\pr_{k+1} \otimes \mathcal{O}\pr_{l+1}$ where $k$ and $l$ are positive
integers.  Throughout, let $g = \text{gcd}(k,l)$.  Also if $k$ and $l$ are both even,
let $k', l' \in \Z_2$ be defined by \begin{align*} k' &= \begin{cases}
	0 & \text{if~ $\tfrac{k}{2} \equiv 0 \pmod g$} \\
	1 & \text{if~ $\tfrac{k}{2} \equiv \tfrac{g}{2} \pmod g$} 
	\end{cases}  \\
\intertext{and}	
l' &= \begin{cases}
	0 & \text{if~ $\tfrac{l}{2} \equiv 0 \pmod g$} \\
	1 & \text{if~ $\tfrac{l}{2} \equiv \tfrac{g}{2} \pmod g$} \; . 
	\end{cases}
\end{align*}	

The simplest case --- to describe and to compute --- is when one of the indices is odd. 
 Table~\ref{cuntzproductA} shows the united $K$-theory of $\mathcal{O}\pr_{k+1}
\otimes \mathcal{O}\pr_{l+1}$ if $k$ is odd.  

If $k$ and $l$ are both even, the united $K$-theory of 
$\mathcal{O}\pr_{k+1} \otimes \mathcal{O}\pr_{l+1}$ falls into three cases:  
both $k$ and $l$ are congruent to 2 modulo 4;
both $k$ and $l$ are congruent to 0 modulo 4;
and $k$ is congruent to 2 and $l$ is congruent to 0 modulo 4.  The results of these 
computations are shown in Tables~\ref{cuntzproductB}, \ref{cuntzproductC}, and
\ref{cuntzproductD} respectively.

\begin{table}[h]  
\caption{$K\crt(\mathcal{O}\pr_{k+1} \otimes \mathcal{O}\pr_{l+1})$ for $k$ odd} 
\label{cuntzproductA}
$$\begin{array}{|c|c|c|c|c|c|c|c|c|c|}  
\hline  \hline  
n & \makebox[1cm][c]{0} & \makebox[1cm][c]{1} & 
\makebox[1cm][c]{2} & \makebox[1cm][c]{3} 
& \makebox[1cm][c]{4} & \makebox[1cm][c]{5} 
& \makebox[1cm][c]{6} & \makebox[1cm][c]{7} 
& \makebox[1cm][c]{8} \\
\hline  \hline 
KO_n  &
\Z_g & \Z_g & 0 & 0 & 
\Z_g & \Z_g & 0 & 0 & 
\Z_g         \\
\hline  
KU_n   &
\Z_g & \Z_g & 
\Z_g & \Z_g & 
\Z_g & \Z_g & 
\Z_g & \Z_g & 
\Z_g         \\
\hline  
KT_n   &
\Z_g^2 & \Z_g & 0 & \Z_g &  
\Z_g^2 & \Z_g & 0 & \Z_g &
\Z_g^2         \\
\hline \hline
c_n & 1 & 1 & 0 & 0 & 2 & 2 & 0 & 0 & 1 \\
\hline
r_n & 2 & 2 & 0 & 0 & 1 & 1 & 0 & 0 & 2  \\
\hline
\varepsilon_n & \smv{1}{0} & 1 & 0 & 0 & \smv{2}{0} & 2 & 0 & 0 & \smv{1}{0} \\
\hline
\zeta_n & \smh{1}{0} & 1 & 0 & 0 & \smh{1}{0} & 1 & 0 & 0 & \smh{1}{0} \\
\hline
(\psi\su)_n & 1 & 1 & -1 & -1 & 1 & 1 & -1 & -1 & 1 \\
\hline
(\psi\st)_n &   \sm{1}{0}{0}{-1} & 1 & 0 & -1 & 
		\sm{1}{0}{0}{-1} & 1 & 0 & -1 &	
		\sm{1}{0}{0}{-1} \\
\hline
\gamma_n & 1 & \smv{0}{1} & 0 & 0 & 1 & \smv{0}{1} & 0 & 0 & 1 \\
\hline
\tau_n & \smh{0}{2} & 0 & 0 & 1 & \smh{0}{1} & 0 & 0 & 2 & \smh{0}{2}  \\
\hline \hline
\end{array}$$
\end{table}

\begin{table}  
\caption{$K\crt(\mathcal{O}\pr_{k+1} \otimes \mathcal{O}\pr_{l+1})$ for 
$k \equiv l \equiv 2 \pmod 4$} \label{cuntzproductB}
$$\begin{array}{|c|c|c|c|c|c|c|c|c|c|}  
\hline  \hline  
n & \makebox[1cm][c]{0} & \makebox[1cm][c]{1} & 
\makebox[1cm][c]{2} & \makebox[1cm][c]{3} 
& \makebox[1cm][c]{4} & \makebox[1cm][c]{5} 
& \makebox[1cm][c]{6} & \makebox[1cm][c]{7} 
& \makebox[1cm][c]{8} \\
\hline  \hline 
KO_n  &
\Z_g & \Z_{2g} & \Z_2^2 & \Z_2^2 & 
\Z_{2g} & \Z_g & 0 & 0 & 
\Z_g         \\
\hline  
KU_n   &
\Z_g & \Z_g & 
\Z_g & \Z_g & 
\Z_g & \Z_g & 
\Z_g & \Z_g & 
\Z_g         \\
\hline  
KT_n   &
\Z_g^2 & \Z_{2g} & \Z_2^2 & \Z_{2g} &  
\Z_g^2 & \Z_{2g} & \Z_2^2 & \Z_{2g} &
\Z_g^2         \\
\hline \hline
c_n & 1 & 1 & \smh{\tfrac{g}{2}}{\tfrac{g}{2}} & 0 & 
	1 & 2 & 0 & 0 & 1 \\
\hline
r_n & 2 & 2 & 0 & \smv{1}{1} & 
	2 & 1 & 0 & 0 & 2  \\
\hline
\varepsilon_n & \smv{0}{1} & 1 & \sm{1}{0}{0}{1} & \smh{g}{g} & 
\left( \begin{smallmatrix} {\tfrac{g}{2}} \\ {\tfrac{g}{2}+1} 
       \end{smallmatrix} \right)	 
			& 2 & 0 & 0 & \smv{0}{1} \\
\hline
\zeta_n & \smh{\tfrac{g}{2}}{1} & 1 & 
	\smh{\tfrac{g}{2}}{\tfrac{g}{2}} & \tfrac{g}{2} &
	\smh{\tfrac{g}{2}}{1} & 1 & 
	\smh{\tfrac{g}{2}}{\tfrac{g}{2}} & \tfrac{g}{2} &
	\smh{\tfrac{g}{2}}{1}   \\
\hline
(\psi\su)_n & 1 & 1 & -1 & -1 & 1 & 1 & -1 & -1 & 1 \\
\hline
(\psi\st)_n &   \sm{-1}{0}{0}{1} & 1 & \sm{1}{0}{0}{1} & -1 & 
		\sm{-1}{0}{0}{1} & 1 & \sm{1}{0}{0}{1} & -1 &	
		\sm{-1}{0}{0}{1} \\
\hline
\gamma_n & 2 & 
	\left( \begin{smallmatrix} {1} \\ {\tfrac{g}{2}} 
		\end{smallmatrix} \right)
	& g & \smv{1}{1} & 
	   2 & 
	 \left( \begin{smallmatrix} {1} \\ {\tfrac{g}{2}} 
		\end{smallmatrix} \right)
	& g & \smv{1}{1} & 2 \\
\hline
\tau_n & \smh{g+2}{g} & \smv{1}{1} & \sm{1}{0}{0}{1} & 1 & 
	 \smh{1}{0} & 0 & 0 & 1 & 
	\smh{g+2}{g}  \\
\hline \hline
\end{array}$$
\end{table}

\begin{table}  
\caption{$K\crt(\mathcal{O}\pr_{k+1} \otimes \mathcal{O}\pr_{l+1})$ for 
$k \equiv l \equiv 0 \pmod 4$} \label{cuntzproductC}
$$\begin{array}{|c|c|c|c|c|c|c|c|c|c|}  
\hline  \hline 
n & \makebox[1cm][c]{0} & \makebox[1cm][c]{1} & 
\makebox[1cm][c]{2} & \makebox[1cm][c]{3} 
& \makebox[1cm][c]{4} & \makebox[1cm][c]{5} 
& \makebox[1cm][c]{6} & \makebox[1cm][c]{7} 
& \makebox[1cm][c]{8} \\
\hline  \hline 
KO_n  &
\Z_g & \Z_2 \oplus \Z_g & \Z_2^3 & \Z_2^3 & 
\Z_2 \oplus \Z_g & \Z_g & 0 & 0 & 
\Z_g         \\
\hline  
KU_n   &
\Z_g & \Z_g & 
\Z_g & \Z_g & 
\Z_g & \Z_g & 
\Z_g & \Z_g & 
\Z_g         \\
\hline  
KT_n   &
\Z_g^2 & \Z_2 \oplus \Z_g & \Z_2^2 & \Z_2 \oplus \Z_g &  
\Z_g^2 & \Z_2 \oplus \Z_g & \Z_2^2 & \Z_2 \oplus \Z_g &
\Z_g^2         \\
\hline \hline
c_n & 1 & \smh{0}{1} & \smaaa{\tfrac{k}{2}}{0}{\tfrac{l}{2}} &
					\smaaa{0}{0}{\tfrac{d}{2}} & 
	\smh{0}{2} & 2 & 0 & 0 & 1 \\
\hline
r_n & 2 & \smv{0}{2} & \smabc{0}{1}{0} & \smabc{l'}{k'}{0} & 
	\smv{0}{1} & 1 & 0 & 0 & 2 \\
\hline
\varepsilon_n & \smv{0}{1} & \sm{1}{0}{0}{1} & \smaaabbb{1}{0}{0}{0}{0}{1} &
			\mmaaabbb{0}{0}{1}{\tfrac{k}{2}}{\tfrac{l}{2}}{0} &
	\mm{\tfrac{g}{2}}{0}{0}{2} & \smv{0}{2} & 0 & 0 & \smv{0}{1}  \\
\hline
\zeta_n & \smh{0}{1} & \smh{0}{1} & \smh{\tfrac{k}{2}}{\tfrac{l}{2}} &
			\smh{\tfrac{g}{2}}{0} &
	\smh{0}{1} & \smh{0}{1} & \smh{\tfrac{k}{2}}{\tfrac{l}{2}} &
			\smh{\tfrac{g}{2}}{0}  & \smh{0}{1}  \\
\hline
(\psi\su)_n & 1 & 1 & -1 & -1 & 1 & 1 & -1 & -1 & 1 \\
\hline
(\psi\st)_n & \sm{-1}{0}{0}{1} & \sm{1}{0}{0}{1} & 
			\sm{1}{0}{0}{1} & \sm{1}{0}{0}{-1} & 
	\sm{-1}{0}{0}{1} & \sm{1}{0}{0}{1} & 
		\sm{1}{0}{0}{1} & \sm{1}{0}{0}{-1} & \sm{-1}{0}{0}{1} \\
\hline
\gamma_n & \smv{0}{1} & \smv{1}{0} & \smv{1}{0} & \smv{l'}{k'} &   
	\smv{0}{1} & \smv{1}{0} & \smv{1}{0} & \smv{l'}{k'} & \smv{0}{1} \\
\hline
\tau_n & \sm{0}{1}{2}{0} & \smaabbcc{0}{l'}{1}{0}{0}{k'} & 
			\smaabbcc{1}{0}{0}{1}{0}{0} & \sm{1}{0}{0}{1} &
	 \smh{1}{0} & 0 & 0 & \smh{0}{2} & \sm{0}{1}{2}{0} \\
\hline \hline
\end{array}$$
\end{table}

\begin{table}  
\caption{$K\crt(\mathcal{O}\pr_{k+1} \otimes \mathcal{O}\pr_{l+1})$ for $k \equiv 2$ 
and $l \equiv 0 \pmod 4$} \label{cuntzproductD}
$$\begin{array}{|c|c|c|c|c|c|c|c|c|c|}  
\hline  \hline  
n & \makebox[1cm][c]{0} & \makebox[1cm][c]{1} & 
\makebox[1cm][c]{2} & \makebox[1cm][c]{3} 
& \makebox[1cm][c]{4} & \makebox[1cm][c]{5} 
& \makebox[1cm][c]{6} & \makebox[1cm][c]{7} 
& \makebox[1cm][c]{8} \\
\hline  \hline 
KO_n  &
\Z_g & \Z_2 \oplus \Z_g & \Z_4 \oplus \Z_2 & \Z_2 \oplus \Z_4 & 
\Z_2 \oplus \Z_g & \Z_g & 0 & 0 & 
\Z_g         \\
\hline  
KU_n   &
\Z_g & \Z_g & 
\Z_g & \Z_g & 
\Z_g & \Z_g & 
\Z_g & \Z_g & 
\Z_g         \\
\hline  
KT_n   &
\Z_g^2 & \Z_2 \oplus \Z_g & \Z_2^2 & \Z_2 \oplus \Z_g &  
\Z_g^2 & \Z_2 \oplus \Z_g & \Z_2^2 & \Z_2 \oplus \Z_g &
\Z_g^2         \\
\hline \hline
c_n & 1 & \smh{0}{1} & \smh{\tfrac{g}{2}}{0} & \smh{0}{\tfrac{g}{2}} & 
	\smh{0}{2} & 2 & 0 & 0 & 1    \\
\hline
r_n & 2 & \smv{0}{2} & \smv{2}{0} & \smv{0}{2} & 
	\smv{0}{1} & 1 & 0 & 0 & 2  \\
\hline
\varepsilon_n & \smv{0}{1} & \sm{1}{0}{0}{1} & 
		\sm{1}{0}{0}{1} & \mm{0}{1}{\tfrac{g}{2}}{0} &
	\mm{\tfrac{g}{2}}{0}{0}{2} & \smv{0}{2} & 0 & 0 & \smv{0}{1}  \\
\hline
\zeta_n & \smh{0}{1} & \smh{0}{1} & 
		\smh{\tfrac{g}{2}}{0} & \smh{\tfrac{g}{2}}{0} &
	\smh{0}{1} & \smh{0}{1} & 
		\smh{\tfrac{g}{2}}{0} & \smh{\tfrac{g}{2}}{0} &
	\smh{0}{1} \\
\hline
(\psi\su)_n & 1 & 1 & -1 & -1 & 1 & 1 & -1 & -1 & 1  \\
\hline
(\psi\st)_n & \sm{-1}{0}{0}{1} & \sm{1}{0}{0}{1} & 
			\sm{1}{0}{0}{1} & \sm{1}{0}{0}{-1} & 
	\sm{-1}{0}{0}{1} & \sm{1}{0}{0}{1} & 
			\sm{1}{0}{0}{1} & \sm{1}{0}{0}{-1} &
	\sm{-1}{0}{0}{1}   \\
\hline
\gamma_n & \smv{0}{1} & \smv{1}{0} & \smv{1}{0} & \smv{0}{1} &  
	\smv{0}{1} & \smv{1}{0} & \smv{1}{0} & \smv{0}{1} & 
	\smv{0}{1}     \\
\hline
\tau_n & \sm{0}{1}{2}{0} & \sm{2}{0}{0}{1} & 
			\sm{1}{0}{0}{2} & \sm{1}{0}{0}{1} & 
	\smh{1}{0} & 0 & 0 & \smh{0}{2} & \sm{0}{1}{2}{0}    \\
\hline \hline
\end{array}$$
\end{table}

Notice that, unlike in the complex case, it is not enough to know the greatest common 
divisor of $k$ and $l$.  Suppose that $k$ and $l$ are both congruent to 2 modulo 4.
 Then we find the united $K$-theory in Table~\ref{cuntzproductB}.  If we replace $l$
by $2l$, then we find the united $K$-theory in Table~\ref{cuntzproductD}.  This
change does not change the greatest common divisor, but it does change the
$K$-theory and therefore the isomorphism class of the tensor product.

To be concrete, let $A = \mathcal{O}\pr_3 \otimes\r \mathcal{O}\pr_3$ and 
$B = \mathcal{O}\pr_3 \otimes\r \mathcal{O}\pr_5$.  Then $A$ and $B$ are simple,
real C*-algebras that are not isomorphic, distinguished by $K$-theory.  However,
$A\c$ and $B\c$ are simple complex C*-algebra that have the same $K$-theory and
hence are isomorphic by the classification theorems of Phillips and Kirchberg (\cite{Phillips} and
\cite{Kirchberg}). 

\begin{proof}[Computation of Table~\ref{cuntzproductA}]
The complex Cuntz algebras are stably isomorphic to a crossed product of an AF-algebra 
by $\Z$ (see \cite{Cun1}).  Hence all of complex Cuntz algebras are in Schochet's
bootstrap category $\mathcal{N}$.  We are therefore entitled to use our K\"unneth
sequence for the united $K$-theory of the product.  

Let $k$ and $l$ be integers with $k$ odd.  We must first compute 
$K\crt(\mathcal{O}\pr_{k+1}) \otimes\scrt K\crt(\mathcal{O}\pr_{l+1})$ and 
$\tor\scrt(K\crt(\mathcal{O}\pr_{k+1}), K\crt(\mathcal{O}\pr_{l+1}))$.  
For this we must build a free resolution of $K\crt(\mathcal{O}\pr_{k+1})$.  In the 
case that $k$ is odd, this is easy since $K\crt(\mathcal{O}\pr_{k+1})$ is generated
by the class of the unit in $KO_0(\mathcal{O}\pr_{k+1})$.  Let $F(b, 0, \R)$ be the
free ${\it CRT}$-module with a single generator in the real part in degree 0.  Our
free resolution for $k$ odd is: $$0 \rightarrow  F(b, 0, \R) \xrightarrow{\mu_1} 
F(b, 0, \R) \xrightarrow{\mu_0}  K\crt(\mathcal{O}\pr_{k+1}) \rightarrow 0 $$
where $\mu_1$ is multiplication by $k$ and $\mu_0$ is the homomorphism which sends the 
generator $b$ in $F(b,0,\R)$ to the class of the unit $[1]$ in
$KO_0(\mathcal{O}\pr_{k+1})$.  When we tensor this resolution on the right by
$K\crt(\mathcal{O}\pr_{l+1})$, we obtain the homomorphism
$$K\crt(\mathcal{O}\pr_{l+1}) \xrightarrow{\mu_1 \otimes 1}
K\crt(\mathcal{O}\pr_{l+1}) \; ,$$ 
referring to Proposition~\ref{tensorR}.  Now
$K\crt(\mathcal{O}\pr_{k+1}) \otimes\scrt K\crt(\mathcal{O}\pr_{l+1})$ is isomorphic
to the cokernel of $\mu_1 \otimes 1$ and $\tor\scrt(K\crt(\mathcal{O}\pr_{k+1}),
K\crt(\mathcal{O}\pr_{l+1}))$ is isomorphic to the kernel of $\mu_1 \otimes 1$.  But
$\mu_1 \otimes 1$ is multiplication by $k$.  Therefore, the ${\it CRT}$-module
$K\crt(\mathcal{O}\pr_{k+1}) \otimes\scrt K\crt(\mathcal{O}\pr_{l+1})$ is easily
computed and is shown in Table~\ref{tensorA}.  The $\tor$ group turns out to be the
same in this case so we will not reproduce it in a separate table.

\begin{table}  
\caption{$M = K\crt(\mathcal{O}\pr_{k+1}) \otimes\scrt K\crt(\mathcal{O}\pr_{l+1})$ for 
$k$ odd} \label{tensorA}
$$\begin{array}{|c|c|c|c|c|c|c|c|c|c|}  
\hline  \hline 
n & \makebox[1cm][c]{0} & \makebox[1cm][c]{1} & 
\makebox[1cm][c]{2} & \makebox[1cm][c]{3} 
& \makebox[1cm][c]{4} & \makebox[1cm][c]{5} 
& \makebox[1cm][c]{6} & \makebox[1cm][c]{7} 
& \makebox[1cm][c]{8} \\
\hline  \hline 
M\so  &
\Z_g & 0 & 0 & 0 & 
\Z_g & 0 & 0 & 0 & 
\Z_g         \\
\hline  
M\su   &
\Z_g & 0 & 
\Z_g & 0 & 
\Z_g & 0 & 
\Z_g & 0 & 
\Z_g         \\
\hline  
M\st   &
\Z_g & 0 & 0 & \Z_g &  
\Z_g & 0 & 0 & \Z_g &
\Z_g         \\
\hline \hline
c_n & 1 & 0 & 0 & 0 & 2 & 0 & 0 & 0 & 1 \\
\hline
r_n & 2 & 0 & 0 & 0 & 1 & 0 & 0 & 0 & 2  \\
\hline
\varepsilon_n & 1 & 0 & 0 & 0 & 2 & 0 & 0 & 0 & 1 \\
\hline
\zeta_n & 1 & 0 & 0 & 0 & 1 & 0 & 0 & 0 & 1 \\
\hline
(\psi\su)_n & 1 & 0 & -1 & 0 & 1 & 0 & -1 & 0 & 1 \\
\hline
(\psi\st)_n & 1 & 0 & 0 & -1 & 1 & 0 & 0 & -1 & 1 \\
\hline
\gamma_n & 1 & 0 & 0 & 0 & 1 & 0 & 0 & 0 & 1 \\
\hline
\tau_n & 0 & 0 & 0 & 1 & 0 & 0 & 0 & 2 & 0  \\
\hline \hline
\end{array}$$
\end{table}

It remains to solve the extension problem given by the K\"unneth Sequence 
\begin{multline*}
0 \longrightarrow 
K\crt(\mathcal{O}\pr_{k+1}) \otimes\scrt K\crt(\mathcal{O}\pr_{l+1}) \mapis{\alpha}
K\crt(\mathcal{O}\pr_{k+1} \otimes \mathcal{O}\pr_{l+1})    \\
\mapis{\beta}
\tor\scrt(K\crt(\mathcal{O}\pr_{k+1}), 
K\crt(\mathcal{O}\pr_{l+1})) \longrightarrow 0 \;. 
\end{multline*}
We show that it splits in this case.  Indeed, the only graded degree in
which there is a possibility of a non-trivial extension is that for
$KT_0(\mathcal{O}\pr_{k+1} \otimes \mathcal{O}\pr_{l+1})$ which is an extension of
$\Z_g$ by $\Z_g$.  However, since $\gamma_0$ is an isomorphism in the ${\it
CRT}$-module $\tor\scrt(K\crt(\mathcal{O}\pr_{k+1}),
K\crt(\mathcal{O}\pr_{l+1}))$ and since $\beta$ is an isomorphism in degree 1, we
have a splitting $s = \gamma_1 \circ \beta^{-1} \circ \gamma_0^{-1}$ for $\beta$ in
the self-conjugate part in degree 0. \end{proof}

The other three cases are somewhat more complicated than the case when $k$ is 
odd.  We will sketch the computation of Table~\ref{cuntzproductC} as a
representative example.

\begin{proof}[Sketch of Computation of Table~\ref{cuntzproductC}]
Let $k$ and $l$ be integers which are both multiples of 4.  The united $K$-theory of 
$\mathcal{O}\pr_{k+1}$ is generated by two elements:  one in
$KO_0(\mathcal{O}\pr_{k+1})$ and one in $KO_2(\mathcal{O}\pr_{k+1})$.  Thus we must
start our resolution with a surjective homomorphism $\mu_0$ from $F(b, 0, \R) \oplus
F(b, 2, \R)$ to $K\crt(\mathcal{O}\pr_{k+1})$.  The kernel of $\mu_0$ turns out to
be isomorphic to $F(b, 0, \C)$.  Hence our resolution is 
$$0 \rightarrow  F(b, 0, \C) \xrightarrow{\mu_1}  
		 F(b, 0, \R) \oplus F(b, 2, \R) \xrightarrow{\mu_0} 
		 K\crt(\mathcal{O}\pr_{k+1}) \rightarrow 0 \; $$
where $\mu_1$ sends the generator $b \in F(b,0,\C)_0 \pu$ to 
$b \oplus \beta^{-1} c b 
	\in \left( F(b, 0, \R) \oplus F(b, 2, \R) \right)_0 \pu$.

Using Propositions~\ref{tensorR} and \ref{tensorC}, we study the tensor product of 
this resolution with $K\crt(\mathcal{O}\pr_{l+1})$.  Then we obtain
Tables~\ref{tensorB} and \ref{torB} by finding the cokernel and kernel
(respectively) of  $$F(b, 0, \C) \otimes\scrt K\crt(\mathcal{O}\pr_{l+1}) 
\xrightarrow{\mu_1 \otimes 1} (F(b, 0, \R) \oplus F(b, 2, \R)) \otimes\scrt
K\crt(\mathcal{O}\pr_{l+1}) \; .$$

\begin{table}  
\caption{$M = K\crt(\mathcal{O}\pr_{k+1}) \otimes\scrt K\crt(\mathcal{O}\pr_{l+1})$ 
for $k \equiv l \equiv 0 \pmod 4$} \label{tensorB}
$$\begin{array}{|c|c|c|c|c|c|c|c|c|c|}  
\hline  \hline 
n & \makebox[1cm][c]{0} & \makebox[1cm][c]{1} & 
\makebox[1cm][c]{2} & \makebox[1cm][c]{3} 
& \makebox[1cm][c]{4} & \makebox[1cm][c]{5} 
& \makebox[1cm][c]{6} & \makebox[1cm][c]{7} 
& \makebox[1cm][c]{8} \\
\hline  \hline 
M\so  &
\Z_g & \Z_2 & \Z_2^3 & \Z_2^2 & 
\Z_2 \oplus \Z_g & \Z_2 & 0 & 0 & 
\Z_g         \\
\hline  
M\su   &
\Z_g & 0 & 
\Z_g & 0 & 
\Z_g & 0 & 
\Z_g & 0 & 
\Z_g         \\
\hline  
M\st   &
\Z_2 \oplus \Z_g & \Z_2 & \Z_2^2 & \Z_g &  
\Z_2 \oplus \Z_g & \Z_2 & \Z_2^2 & \Z_g &
\Z_2 \oplus \Z_g         \\
\hline \hline
c_n & 1 & 0 & \smaaa{\tfrac{k}{2}}{0}{\tfrac{l}{2}} & 0 & \smh{0}{2} & 
0 & 0 & 0 & 1 \\
\hline
r_n & 2 & 0 & \smabc{0}{1}{0} & 0 & \smv{0}{1} & 0 & 0 & 0 & 2  \\
\hline
\varepsilon_n & \smv{0}{1} & 1 & \smaaabbb{1}{0}{0}{0}{0}{1} & 
\smh{\tfrac{k}{2}}{\tfrac{l}{2}} & 
\sm{1}{0}{0}{2} & 0 & 0 & 0 & \smv{0}{1} \\
\hline
\zeta_n & \smh{0}{1} & 0 & \smh{\tfrac{k}{2}}{\tfrac{l}{2}} & 0 &
 \smh{0}{1} & 0 & \smh{\tfrac{k}{2}}{\tfrac{l}{2}} & 0 &
\smh{0}{1}   \\
\hline
(\psi\su)_n & 1 & 0 & -1 & 0 & 1 & 0 & -1 & 0 & 1 \\
\hline
(\psi\st)_n & 1 & 1 & \sm{1}{0}{0}{1} & -1 & 1 & 1 & \sm{1}{0}{0}{1} & -1 & 1\\
\hline
\gamma_n & 1 & 0 & 1 & 0 & 1 & 0 & 1 & 0 & 1 \\
\hline
\tau_n & \smh{0}{1} & \smabc{0}{1}{0} & \sm{1}{0}{0}{1} & \smv{0}{1} & 
\smh{1}{0} & 0 & 0 & 2 & \smh{0}{1} \\
\hline \hline
\end{array}$$
\end{table}

\begin{table}   
\caption{$M = \tor\scrt(K\crt(\mathcal{O}\pr_{k+1}), K\crt(\mathcal{O}\pr_{l+1}))$ for 
$k \equiv l \equiv 0 \pmod 4$} \label{torB}
$$\begin{array}{|c|c|c|c|c|c|c|c|c|c|}  
\hline  \hline 
n & \makebox[1cm][c]{0} & \makebox[1cm][c]{1} & 
\makebox[1cm][c]{2} & \makebox[1cm][c]{3} 
& \makebox[1cm][c]{4} & \makebox[1cm][c]{5} 
& \makebox[1cm][c]{6} & \makebox[1cm][c]{7} 
& \makebox[1cm][c]{8} \\
\hline  \hline 
M\so  &
\Z_g & 0 & \Z_2 & 0 & 
\Z_{\tfrac{g}{2}} & 0 & 0 & 0 & 
\Z_g         \\
\hline  
M\su   &
\Z_g & 0 & 
\Z_g & 0 & 
\Z_g & 0 & 
\Z_g & 0 & 
\Z_g         \\
\hline  
M\st   &
\Z_g & 0 & \Z_2 & \Z_{\tfrac{g}{2}} &  
\Z_g & 0 & \Z_2 & \Z_{\tfrac{g}{2}} &
\Z_g         \\
\hline \hline
c_n & 1 & 0 & \tfrac{g}{2} & 0 & 2 & 0 & 0 & 0 & 1 \\
\hline
r_n & 2 & 0 & 0 & 0 & 1 & 0 & 0 & 0 & 2  \\
\hline
\varepsilon_n & 1 & 0 & 1 & 0 & 2 & 0 & 0 & 0 & 1 \\
\hline
\zeta_n & 1 & 0 & \tfrac{g}{2} & 0 & 1 & 0 & \tfrac{g}{2} & 0 & 1 \\
\hline
(\psi\su)_n & 1 & 0 & -1 & 0 & 1 & 0 & -1 & 0 & 1 \\
\hline
(\psi\st)_n & 1 & 0 & 1 & -1 & 1 & 0 & 1 & -1 & 1 \\
\hline
\gamma_n & 1 & 0 & 0 & 0 & 1 & 0 & 0 & 0 & 1 \\
\hline
\tau_n & 0 & 0 & 0 & 1 & 0 & 0 & 0 & 2 & 0  \\
\hline \hline
\end{array}$$
\end{table}

It remains to determine the extension
$$0 \longrightarrow 
K\crt(\mathcal{O}\pr_{k+1}) \otimes\scrt K\crt(\mathcal{O}\pr_{l+1}) \mapis{\alpha}
K\crt(\mathcal{O}\pr_{k+1} \otimes \mathcal{O}\pr_{l+1}) \hspace{3cm}
$$
$$
\hspace{6cm} \mapis{\beta}
\tor\scrt(K\crt(\mathcal{O}\pr_{k+1}), K\crt(\mathcal{O}\pr_{l+1})) 
\longrightarrow 0 \; . $$  Fortunately, the structure imposed by the ${\it
CRT}$-relations and the acyclicity of united $K$-theory will be enough to determine
$K\crt(\mathcal{O}\pr_{k+1} \otimes\r \mathcal{O}\pr_{l+1})$.

Let $A$ denote $\mathcal{O}\pr_{k+1} \otimes \mathcal{O}\pr_{l+1}$.  We can immediately 
establish $KO_*(A)$ in degrees 0, 2, 4, 6, and 7; $KU_*(A)$ in all degrees; and
$KT_*(A)$ in degrees 2 and 6.  The realification and complexification operations are
also established in the even degrees.  

In the following, we will make frequent use of the long exact sequence 
$$
\cdots \longrightarrow  KO_n(A) \mapis{\eta\so} 
KO_{n+1}(A) \mapis{c} 
KU_{n+1}(A) \mapis{r \beta\su^{-1}}
KO_{n-1}(A) \longrightarrow \cdots \; . 
$$

First, we show that $KO_1(A)$ is isomorphic to $\Z_2 \oplus \Z_g$, rather than 
$\Z_{2g}$.  Indeed, the kernel of $c_2 \colon \Z_2^3 \rightarrow \Z_g$ must have at
least two generators.  Since the kernel of $c_2$ is the image of $(\eta\so)_1$, it
follows that $KO_1(A)$ must have two generators.

Second, we show that $KO_5(A)$ is isomorphic to $\Z_g$.  Indeed, $KO_6(A) = KO_7(A) = 0$,
 so the long exact sequence above shows that $r_5$ is an isomorphism.  It also
follows from this and the relation $rc = 2$ that the image of $c_5$ is $2 \Z_g$.

Third, we show that $KO_3(A)$ is isomorphic to $\Z_2^3$ rather than $\Z_2 \oplus \Z_4$.  
Assume that $KO_3(A) \cong \Z_2 \oplus \Z_4$.  Because the image of $c_5$ is $2
\Z_g$, the realification map $r_3$ must take the generator of $KU_3(A)$ to an
element of degree 2.  On the other hand, the image of $(\eta\so)_3$ (which is the
kernel of $c_4$) has two generators.  Therefore, the image of $r_3$ cannot be the
$\Z_2$ summand of $\Z_2 \oplus \Z_4$.  By re-choosing the basis, we may assume $r_3
= \smv{0}{2}$.

The image of $(\eta\so)_2$ is a copy of $\Z_2 \oplus \Z_2$.  As there is only one such 
subgroup of $\Z_2 \oplus \Z_4$, this determines $c_3 = \smh{0}{\tfrac{g}{2}}$. 
However, then we have $rc = \sm{0}{0}{0}{g}$ which is not multiplication by 2.  By
this contradiction, we conclude that $KO_3(A) \cong \Z_2^3$.  This completes the
computation of the groups $KO_*(A)$.

Because $KO_6(A) \cong KO_7(A) \cong 0$, the long exact sequence 
$$
\cdots \longrightarrow  KO_n(A) \mapis{\eta\so^2} 
KO_{n+2}(A) \mapis{\varepsilon} 
KT_{n+2}(A) \mapis{\tau \beta\st^{-1}}
KO_{n-1}(A) \longrightarrow \cdots \;  
$$
shows that $\varepsilon_1$ and $\tau_3$ are isomorphisms.  Thus we know $KT_*(A)$ in 
degrees 1 and 3, and hence in degrees 5 and 7.  It remains to find $KT_0(A) \cong
KT_4(A)$.  The sequence $$KU_1(A) \xrightarrow{1 - \psi\su}
KU_1(A) \xrightarrow{\gamma}
KT_0(A) \xrightarrow {\zeta}
KU_0(A) \xrightarrow{1 - \psi\su}
KU_0(A)
$$
becomes
$$\Z_g \xrightarrow{0} \Z_g \xrightarrow{\gamma}
KT_0(A) \xrightarrow{\zeta} \Z_g \xrightarrow{0} \Z_g $$
so $KT_0(A)$ is an extension of $\Z_g$ by $\Z_g$.  We claim that it is isomorphic to 
$\Z_g \oplus \Z_g$.  For this we also need the sequence
$$KO_{-2}(A) \xrightarrow{\eta\so^2}
KO_0(A) \xrightarrow{\ve}
KT_0(A) \xrightarrow{\tau \beta\st^{-1}}
KO_{-3}(A) \xrightarrow{\eta\so^2}
KO_{-1}(A)
$$
which becomes
$$0 \rightarrow
\Z_g \xrightarrow{\ve}
KT_0(A) \rightarrow
\Z_g \rightarrow 0 \, .$$
Since $c_0 = \zeta_0 \ve_0 \colon \Z_g \rightarrow \Z_g$ is an isomorphism, the image of 
$\ve_0$ has a trivial intersection with $\ker \zeta_0 = \im \gamma_1$.  Therefore,
$KT_0(A) = \im \ve \oplus \im \gamma \cong \Z_g \oplus \Z_g$.

Now that all of the groups of $K\crt(A)$ have been established, the behavior of the 
operations can be computed using the exact sequences and the {\it CRT}-relations.
\end{proof}

This example shows that the K\"unneth sequence for united $K$-theory does not split 
in general, even on the level of abelian groups.  




%
%
\end{document}